\documentclass[12pt]{amsart}
\usepackage{amscd,amsmath,amssymb,amsfonts}
\usepackage[cmtip, all]{xy}

\newtheorem{thm}{Theorem}[section]
\newtheorem{prop}[thm]{Proposition}
\newtheorem{lem}[thm]{Lemma}
\newtheorem{cor}[thm]{Corollary}

\theoremstyle{definition}

\newtheorem{defn}[thm]{Definition}
\theoremstyle{remark}
\newtheorem{remk}[thm]{Remark}
\newtheorem{remks}[thm]{Remarks}
\newtheorem{exm}[thm]{Example}
\newtheorem{exms}[thm]{Examples}
\newtheorem{notat}[thm]{Notation}
\numberwithin{equation}{section}

{\hfill$\square$\end{defn}}

\newenvironment{rem}{\begin{remk}}%
{\hfill$\square$\end{remk}}

{\hfill$\square$\end{remks}}

\newenvironment{ex}{\begin{exm}}%
{\hfill$\square$\end{exm}}

{\hfill$\square$\end{exms}}

{\hfill$\square$\end{notat}}

\newcommand{\sA}{{\mathcal A}}

\newcommand{\sC}{{\mathcal C}}
\newcommand{\sD}{{\mathcal D}}

\newcommand{\sO}{{\mathcal O}}

\newcommand{\sT}{{\mathcal T}}
\newcommand{\sU}{{\mathcal U}}

\newcommand{\sW}{{\mathcal W}}

\newcommand{\sZ}{{\mathcal Z}}
\newcommand{\A}{{\mathbb A}}

\newcommand{\G}{{\mathbb G}}

\renewcommand{\P}{{\mathbb P}}
\newcommand{\Q}{{\mathbb Q}}

\newcommand{\W}{{\mathbb W}}

\newcommand{\Z}{{\mathbb Z}}

\newcommand{\CH}{{\rm CH}}

\newcommand{\surj}{\twoheadrightarrow}
\newcommand{\red}{{\rm red}}
\newcommand{\codim}{{\rm codim}}

\newcommand{\Hom}{{\rm Hom}}

\newcommand{\Spec}{{\rm Spec \,}}

\newcommand{\supp}{{\rm supp}\,}
\newcommand{\0}{\emptyset}
\newcommand{\sHom}{{\mathcal{H}{om}}}

\newcommand{\id}{{\operatorname{id}}}

\newcommand{\Sch}{{\operatorname{\mathbf{Sch}}}}

\newcommand{\op}{{\text{\rm op}}}

\newcommand{\Sets}{{\mathbf{Sets}}}
\newcommand{\del}{\partial}

\renewcommand{\max}{{\operatorname{\rm max}}}

\newcommand{\Sm}{{\mathbf{Sm}}}
\newcommand{\SmProj}{{\mathbf{SmProj}}}

\newcommand{\GL}{{\operatorname{\rm GL}}}
\newcommand{\SL}{{\operatorname{\rm SL}}}

 \newcommand{\Ab}{{\mathbf{Ab}}}
\newcommand{\Tot}{{\operatorname{\rm Tot}}}

\newcommand{\End}{{\operatorname{\text{End}}}}

\newcommand{\ds}{{/\kern-3pt/}}

\renewcommand{\log}{{\operatorname{log}}}

 \newcommand{\gen}{{\text{gen}}}
\newcommand{\Gr}{{\text{\rm Gr}}}

\newcommand{\Proj}{{\operatorname{Proj}}}

\newcommand{\colim}{\mathop{\text{colim}}}

\newcommand{\Cube}{{\mathbf{Cube}}}

\newcommand{\Mot}{{\operatorname{Mot}}}

\newcommand{\Cor}{{\operatorname{Cor}}}
\newcommand{\DMH}{\operatorname{DM}^H}

\newcommand{\TZ}{{\operatorname{Tz}}}
\renewcommand{\TH}{{\operatorname{TCH}}}
\newcommand{\un}{\underline}
\newcommand{\ov}{\overline}
\renewcommand{\hom}{\text{hom}}
\newcommand{\dgn}{{\operatorname{degn}}}
\renewcommand{\dim}{\text{dim}}

\begin{document}
\title{Additive higher Chow groups of schemes}
\author{Amalendu Krishna, Marc Levine}
\address{}
\email{}

\date{January 10, 2007}


\thanks{The second author gratefully acknowledges the support of the Humboldt Foundation, and support of the NSF via the grant DMS-0457195}

\begin{abstract}
We show how to make the additive Chow groups of Bloch-Esnault, R\"ulling and Park into a graded module for Bloch's higher Chow groups, in the case of a smooth projective variety over a field. This yields a  a projective bundle formula as well as a blow-up formula for  the additive Chow groups of a smooth projective variety.

In case the base-field admits resolution of singularieties, these properties allow us to apply the technique of Guill\'en and Navarro Aznar to define the additive 
Chow groups ``with log poles at infinity" for an arbitrary finite-type 
$k$-scheme $X$. This theory has all the usual properties of a Borel-Moore theory on finite type $k$-schemes: it is covariantly functorial for projective morphisms, contravariantly functorial for morphisms of smooth schemes,  and has a projective bundle 
formula, homotopy property, Mayer-Vietoris and localization sequences.

Finally, we show that the regulator map defined by Park from the additive Chow groups of $1$-cycles to the  modules of absolute K\"ahler differentials of an algebraically closed field of
characteristic zero is surjective, giving an evidence of conjectured isomorphism between these two groups.
\end{abstract} 

\maketitle   

\tableofcontents

\section*{Introduction}
Since Bloch's introduction of his higher Chow groups over last twenty years 
ago, Voevodsky's construction of a good triangulated category of motives, 
with the higher Chow groups as the corresponding {\em motivic cohomology} has 
given what was at first a provisional definition a solid foundation. However, 
when compared with algebraic $K$-theory, this theory of motivic cohomology 
still has at least one serious deficiency: it fails to see the difference 
between a non-reduced scheme $X$ and its associated reduced closed subscheme 
$X_\red$. 

One important case of a non-reduced scheme is of course that of the 
infinitesimal thickening of a reduced scheme.    S. Bloch and H. Esnault have 
introduced in \cite{BlochEsnault} the first attempt at a cycle-theoretic 
description of this motivic cohomology `with modulus", defining groups which 
are supposed to describe the part given by zero-cycles for $k[t]/t^2$. They 
showed that  these groups  coincide
with the absolute differential forms $\Omega^n_k$. In \cite{BlochEsnault2} 
they defined 
a cubical version of these  {\em additive} Chow groups, extending the 
definition to zero-cycles for $k[t]/t^{m+1}$, and showed that for modulus 
$m=1$ these groups agreed with the ones defined earlier. 

K. R\"ulling  \cite{Ruelling} generalized these computations to arbitrary 
modulus $m$, showing that the
zero-cycles with modulus $m$ leads to a cycle theoretic 
explanation the generalized deRham-Witt complex of 
Hasselholt-Madsen. Park \cite{Park} extended the definition of zero-cycles 
with modulus $m$ to a full cycle complex, and made some computations for 
1-cycles.

Our aim in this paper is to understand the structural properties of Park's 
additive cycle complexes for arbitrary modulus. For a general smooth 
quasi-projective variety, these complexes do not have the nice behavior of 
Bloch's cycle complexes. The homotopy property quite obviously fails, but 
more seriously, these complexes do not have any reasonable localization 
property with respect the closed subschemes (even smooth ones). 
 However, for smooth projective varieties, they seem to have 
interesting properties. Our first main technical result is the construction 
of  a module structure for the additive Chow groups (with arbitrary modulus) 
of  $X$ over Chow ring of $X$ leading immediately to an 
action of correspondences (for smooth projective schemes over a a smooth $k$-scheme $S$) on the additive Chow groups. This gives us by 
standard constructions a projective bundle formula for arbitrary smooth $X$ and a blow-up formula for 
the additive Chow groups of a smooth projective $X$. 

If we now assume that our base field $k$ admits resolution of singularities, 
we may apply the method  of Guill\'en and Navarro Aznar to define the additive 
Chow groups ``with log poles at infinity" for an arbitrary finite-type 
$k$-scheme $X$. This construction requires a detour through an integral 
version of Hanamura's triangulated category of motives. The resulting theory 
is a Borel-Moore theory, in that it is covariantly functorial for projective 
morphisms in $\Sch_k$. We add to the general theory of Guill\'en and Navarro Aznar and its application by Hanamura in that we construct in addition 
pull-back (Gysin) maps for morphisms of smooth quasi-projective varieties. 
This theory of additive Chow groups (for arbitrary modulus) with log poles 
has all nice structural properties one could hope for: projective bundle 
formula, homotopy property, Mayer-Vietoris and localization sequences. 

As an application
of this action of 'higher correspondences', we show that the regulator
maps \cite{Park} from the additive Chow groups of $1$-cycles to the 
modules of absolute K\"ahler differentials of an algebraically closed field of
characteristic zero is surjective, giving an evidence of conjectured
isomorphism between these two groups.

Our hope is that the Atiyah-Hirzebruch spectral sequence for the algebraic $K$-theory of a smooth $k$-scheme $X$:
\[
E^{p,q}_2:=\CH^{-q}(X,-p-q)\Longrightarrow K_{-p-q}(X)
\]
(see \cite{BL,FriedSus,LevineAH}) has an analog with modulus and log poles:
\[
E^{p,q}_2:=\TH_\log^{-q}(X,-p-q;m)\Longrightarrow K^\log_{-p-q}(X;m);
\]
here $\TH_\log^{-q}(X,-p-q;m)$ is the logarithmic additive Chow group with modulus $m$  (see \S~\ref{subsec:LogChow} for details). For $X$ smooth and projective over $k$, $K^\log(X;m)$ is defined as the (-1 connected) homotopy fiber of the restriction map
\[
K(X\times\Spec k[t])\to K(X\times\Spec k[t]/(t^{m+1});
\]
in general, the logarithmic $K$-groups are at present not even defined, so the ``correct" definition is part of the conjecture. At the very least, the groups $K^\log_n(X;m)$ should be the homotopy groups of a spectrum $K^\log(X;m)$, the assignment $X\mapsto K^\log(X;m)$ should extend to a functor on at least smooth quasi-projective $k$-schemes, covariant for projective maps and contravariant for arbitrary maps (where we consider $K^\log(X;m)$ as an object in the stable homotopy category of spectra). In addition, the spectra $K^\log(X;m)$ should satisfy an $\A^1$ homotopy invariance property, Nisnevich excision, a projective bundle formula and a blow-up formula. 
 
Our paper is orgainzed as follows. We begin with a general discussion of cubical objects in an addtive or abelian category, then recall the definition and basic properties of the cubical version of Bloch's cycle complex. In section~\ref{sec:MovingLemmas} we recall the ``Chow's moving lemma" for Bloch's cycle complex, as the techniques used to prove this result will form the core of the technical results here. Next, in section~\ref{sec:NewMovingLemma} we show how to refine the moving lemma to take into account behavior on the closure of the algebraic $n$-cubes, as is needed for applications to the additive cycle complexes. We recall the classical category of Chow motives and Hanamura's construction of a triangulated category of motives in section~\ref{sec:Motives}; we also adapt Hanamura's construction to our purposes in this section. We also show at this point  how the constructions of Guill\'en and Navarro Aznar allow one to define the (Borel-Moore) motive of an arbitrary finite-type $k$-scheme, assuming one has resolution of singularities, and we finish the section by construction of the pull-back maps for morphisms in the category $\Sm/k$ of smooth quasi-projective varieties over 
$k$.

In section~\ref{sec:AdditiveChow}, we introduce our main object of study, the additive cycle complexes for arbitrary modulus, together with their elementrary properties. We construct the action of higher Chow groups in section~\ref{sec:ChowAction}, and extend this in section~\ref{sec:AdditiveMotive} to define the additive higher Chow groups of a Chow $S$-motive. This leads to the projective bundle formula and the blow-up formula. In section~\ref{sec:log}, we use the method of Guill\'en and Navarro Aznar to define the logarithmic additive higher Chow groups of a scheme of finite type over $k$, assuming $k$ admits resolution of singularities; the theory is covariant for projective maps and contravariant for arbitrary maps of smooth varieties, satisfies a projective bundle formula, homotopy property, localization and Mayer-Vietoris and a blow-up formula.
Finally, in section~\ref{sec:1cycle} we study the regulator maps from the 
additive Chow group of $1$-dimensional cycles to the modules of absolute
K\"ahler differentials of a characteristic zero field, as constructed by Park
\cite{Park}. We establish some properties of the regulator maps and then use
the results of section~\ref{sec:ChowAction} to show the surjectivity of these
regulator maps. 

\section{Cycle complexes}

\subsection{Cubical objects}
We introduce the ``cubical category" $\Cube$. This is the subcategory of $\Sets$ with objects  $\un{n}:=\{0,\infty\}^n$, $n=0,1,2,\ldots$, and morphisms generated by\\
\\
1. Inclusions: $\eta_{n,i\epsilon}:\un{n-1}\to \un{n}$, $\epsilon=0,\infty$, $i=1,\ldots, n$
\[
\eta_{n,i,\epsilon}(y_1,\ldots, y_{n-1})=(y_1,\ldots, y_{i-1},\epsilon,y_i,\ldots, y_{n-1})
\]
2. Projections: $p_i:\un{n}\to\un{n-1}$, $i=1,\ldots, n$.\\
3. Permutations of factors: 
$(\epsilon_1,\ldots, \epsilon_n)\mapsto (\epsilon_{\sigma(1)},\ldots, \epsilon_{\sigma(n)})$ for $\sigma\in S_n$.
4. Involutions: $\tau_{n,i}$ exchanging $0$ and $\infty$ in the $i$th factor of $\un{n}$.\\
\\
Clearly all the Hom-sets in $\Cube$ are finite. For a category $\sA$, we call a functor $F:\Cube^\op\to \sA$ a {\em cubical object} of $\sA$.

\begin{ex}  Let $k$ be a field and set   $\P^1:=\Proj k[Y_0,Y_1]$, and let $y:=Y_1/Y_0$ be the standard
coordinate function on $\P^1$. We set $\square^n:=(\P^1\setminus\{1\})^n$. 

Letting $p_i:  ({\P}^{1})^n\to({\P}^{1})^{n-1}$ be the projection , we use the (rational) coordinate system 
$(y_1, \cdots , y_{n-1})$ on $\square^n$, with $y_i:=y\circ p_i$.

A {\em face} of $\square^n$ is  a subscheme $F$ defined by equations of the form
 \[
 y_{i_1}=\epsilon_1, \ldots,  y_{i_s}=\epsilon_s;\ \epsilon_j\in\{0,\infty\}.
 \]

Let $\eta_{n,i,\epsilon}:\square^{n-1}\to \square^n$ be the inclusion
\[
\eta_{n,i,\epsilon}(y_1,\ldots, y_{n-1})=(y_1,\ldots, y_{i-1},\epsilon,y_i,\ldots, y_{n-1})
\]
Thus $n\mapsto \square^n$ is a functor $\square:\Cube\to \Sm/k$.
\end{ex}

Let $\sA$ be an abelian category, $\un{A}:\Cube\to \sA$ be a cubical object. Let $\pi_{n,i}:\un{A}(\un{n})\to \un{A}(\un{n})$ be the endomorphism $p_i^*\circ \eta_{n,i,\infty}^*$, and set
\[
\pi_n:=(\id-\pi_{n,n})\circ\ldots\circ (\id-\pi_{n,1}).
\]
Let 
\[
\un{A}(\un{n})_0:=\cap_{i=1}^n\ker \eta_{n,i,\infty}^*\subset \un{A}(\un{n})
\]
and
\[
\un{A}(\un{n})_\dgn:=\sum_{i=1}^np_i^*(\un{A}(\un{n-1}))\subset \un{A}(\un{n}).
\]
Finally let $(\un{A}_*,d)$ be the complex with $\un{A}_n:=\un{A}(\un{n})$ and with
\[
d_n:=\sum_{i=1}^n(-1)^i(\eta_{n,i,\infty}^*-\eta_{n,i,0}^*).
\]

The following result is the basis of all ``cubical" constructions; the proof is elementary and is left to the reader.
\begin{lem}
 Let $\un{A}:\Cube\to\sA$ be a cubical object in an abelian category $\sA$. Then\\
\\
1. For each $n$, $\pi_n$ maps $\un{A}(\un{n})$ to $\un{A}(\un{n})_0$ and defines a splitting
\[
\un{A}(\un{n})=\un{A}(\un{n})_\dgn\oplus \un{A}(\un{n})_0.
\]
2. $d_n(\un{A}(\un{n})_\dgn)=0$, $d_n(\un{A}(\un{n})_0)\subset \un{A}(\un{n-1})_0$
\end{lem}

\begin{defn}
 Let $\un{A}:\Cube\to\sA$ be a cubical object in an abelian category $\sA$. Define the complex $(A_*,d)$ to be the quotient complex  $n\mapsto (\un{A}(\un{n})/\un{A}(\un{n})_\dgn,d_n)$ of $\un{A}_*,d)$.
\end{defn}

The lemma shows that $(A_*,d)$ is well-defined and is isomorphic to the subcomplex $n\mapsto (\un{A}(\un{n})_0,d_n)$ of $\un{A}_*$.

\subsection{Bloch's cycle complex}

We recall the cubical version of Bloch's cycle complexes.

\begin{defn}
[{\it cf.} \cite{BlochAlgCyc}] Let $X$ be $k$-scheme of finite type. For $r\in \Z$, $n \ge 0$, we  let $\un{z}_r(X, n)$ be the free abelian group 
on integral closed subschemes $Z$ of $X \times \square^n$ of dimension $r+n$ over $k$
such that: \\ \\
(Good position) For each face $F$ of $\square^n$, $Z$ intersects  $X \times F$ properly:
\[
\dim_k Z\cap X\times F\le r+\dim_kF
\]
\end{defn}

it follows directly from the definition of $\un{z}_r(X, n)
$ that for $\epsilon=0,\infty$,   $\partial^\epsilon_{n,i}:=\eta_{n,i,\epsilon}^*$  gives a well-defined map
\[
\partial^\epsilon_{n,i}:\un{z}_r(X, n)
\to \un{z}_r(X, n-1),
\]
so we have a cubical abelian group $\un{n}\mapsto \un{z}_r(X, n)
$.

\begin{defn}
 Let
\[
z_r(X, n):=\un{z}_r(X, n)/\un{z}_r(X, n)_\dgn.
\]
The complex $(z_r(X,*),d)$ associated to the cubical abelian group 
\[
n\mapsto \un{z}_r(X, n)
\]
 is the {\em cubical cycle complex} of dimension $r$ cycles on $X$.
Set
\[
\CH_r(X,n):=H_n(z_r(X,*)).
\]
If $X$ is equi-dimensional of dimension $d$ over $k$, set
\[
z^i(X,*):=z_{i-r}(X,*)
\]
and
\[
\CH^i(X,n):=H_n(z^i(X,*)).
\]
\end{defn}

Recall from \cite{BlochAlgCyc} the original simplicial version of this construction, Bloch's cycle complex $\sZ_r(X,*)$ and the {\em higher Chow groups} 
\[
\CH_r(X,n):=H_n(\sZ_r(X,*)).
\]

\begin{thm}[\cite{LevineChowRev}] \label{thm:CubicalComparison}There is a natural isomorphism in $D^-(\Ab)$
\[
z_r(X,*)\cong \sZ_r(X,*),
\]
hence a natural isomorphism
\[
\CH_r(X,n)\cong H_n(\sZ_r(X,*)).
\]
\end{thm}

\begin{rem} Both cycle complexes $\sZ_r(X,*)$ and $z_r(X,*)$ are covariantly functorial for proper morphisms, and contravariantly functorial for flat morphisms (with a shift in $r$). The naturality in
the above comparison theorem is with respect to these two functorialities.
\end{rem}

\subsection{Products}
The main advantage the cubical complexes have over the simplicial version is the existence of an associative external product on the level of complexes.

For $X,Y\in\Sch_k$, let 
\[
\tau:X\times\square^n\times Y\times\square^m\to 
X\times Y\times\square^n\times\square^m
\]
be the exchange-of-factors isomorphism. For integral cycles $Z$ on $X\times\square^n$, $W$ on $Y\times\square^m$, define 
\[
Z\boxtimes W:=\tau_*(Z\times W)
\]
where $Z\times W$ the cycle associated to the subscheme $Z\times_kW$. Extending by linearity gives us
\[
\boxtimes_{n,m}:\un{z}_r(X,n)\otimes\un{z}_s(Y,m)\to \un{z}_{r+s}(X\times_kY,n+m).
\]
\begin{lem}[cf. \cite{LevineChowRev}] The maps $\boxtimes_{n,m}$ define a map of complexes
\[
\boxtimes_{X,Y}:z_r(X,*)\otimes z_s(Y,*)\to z_{r+s}(X\times_kY,*)
\]
called {\em external product}. The external product is associative. Up to a natural homotpy, the external product is commutative.
\end{lem}

\begin{proof} It is clear that 
\begin{multline*}
\boxtimes_{n,m}\left(\un{z}_r(X,n)_\dgn\otimes\un{z}_s(Y,m)+\un{z}_r(X,n)_\dgn\otimes\un{z}_s(Y,m)\right)\\
\subset  \un{z}_{r+s}(X\times_kY,n+m)_\dgn. 
\end{multline*}
The fact that $\boxtimes$ is a map of complexes follows directly from the definition of the boundary maps. the associativity is obvious from the definition.

The homotopy commutativity is proved in \cite{LevineChowRev}.
\end{proof}

\subsection{Moving lemmas}\label{sec:MovingLemmas}

We recall  the application of the classical moving lemma of Chow to Bloch's cycle complex. For $X\in\Sch_k$, let $X_{sm}$ be the largest open subscheme of $X$ which is smooth over $k$.

\begin{defn}
 For $X\in\Sch_k$ equi-dimensional over $k$, let $\sC$ be a finite set of locally closed irreducible subsets of $X_{sm}$ and let $e:\sC\to \Z_+$ be a function. We define
$\un{z}^r(X,n)_{\sC, e}\subset \un{z}^r(X,n)$ to be the subgroup generated by integral closed subschemes $Z\subset X\times\square^n$ such that, for each $C\in\sC$ and each face $F$ of $\square^n$, 
\[
\codim_{C\times F}Z\cap (C\times F)\ge r-e(C).
\]
\end{defn}

The $\un{z}^r(X,n)_{\sC, e}$ form a subcomplex $\un{z}^r(X,*)_{\sC, e}$ of 
$\un{z}^r(X,n)$. In addition, for $z_1,\ldots, z_n\in \un{z}^r(X,n-1)$,
\[
\sum_i(\id\times p_i)^*(z_i)\in \un{z}^r(X,n)_{\sC, e}\Leftrightarrow
z_i\in \un{z}^r(X,n-1)_{\sC, e}\text{ for all }i=1,\ldots, n.
\]
Thus, with the evident definition of $\un{z}^r(X,*)_{\sC, e,\dgn}\subset \un{z}^r(X,*)_{\sC, e}$, the map 
\[
\frac{\un{z}^r(X,*)_{\sC, e}}{\un{z}^r(X,*)_{\sC, e,\dgn}}\to
\frac{\un{z}^r(X,*)}{\un{z}^r(X,*)_\dgn}=z^r(X,*)
\]
is injective, and we have the subcomplex $z^r(X,*)_{\sC, e}$ of $z^r(X,*)$.

The second moving lemma is

\begin{thm}\label{thm:ML2} Let $X$ be a quasi-projective $k$-scheme, equi-di\-men\-sional over $k$, $\sC$ a finite set of  locally closed subsets of $X_{sm}$ and $e:\sC\to \Z_+$ a function. Then the inclusion
$z^r(X,*)_{\sC, e}\to z^r(X,*)$ is a quasi-isomorphism.
\end{thm}

As we will be using these techniques later on, we recall the idea of the proof, which goes in three steps:\\
\\
{\bf Step 1: $X=\P^n$.} In this case, one use the action of $G:=\SL_{n+1}$ on $\P^n$. Let $K=k(G)$, a pure transcendental extension of $k$. Since $K$ is a field, we can express the generic matrix $\eta\in\SL_{n+1}$ over $k$ as a product of elementary matrices with  $K$-coefficients. 

Thus, there is a map
\[
\phi:\A^1_K\to G
\]
with $\phi(0)=\id$, $\phi(1)=\eta$. Via the action of $G$ on $\P^n$, we have the map
\[
\mu_\phi:\P^n\times\A^1_K\to \P^n.
\]

Identifying $(\A^1,0,1)$ with $(\square^1,0,\infty)$, pull-back by the maps 
\[
\mu_\phi\times\id:\P^n\times\square^{n+1}_K\to \P^n\times\square^n
\]
define a degree 1 map of complexes
\[
H_\phi:z^r(\P^n,*)\to z^r(\P^n_K,*)
\]
that gives a homotopy between the base-extension $\pi:z^r(\P^n,*)\to z^r(\P^n_K,*)$ and $T_\eta^*$, where $T_\eta:\P^n\to \P^n_K$ is translation by $\eta$.

Both $H_\phi$ and $T_\eta^*$ map $z^r(\P^n,*)_{\sC,e}$ to $z^r(\P^n_K,*)_{\sC,e}$. Since $\eta$ is generic over $k$, $T_\eta^*$ maps $z^r(\P^n,*)$ to $z^r(\P^n_K,*)_{\sC,e}$. As in the proof of the localization theorem, this implies that base-extension
\[
\pi:\frac{z^r(\P^n,*)}{z^r(\P^n,*)_{\sC,e}}\to 
\frac{z^r(\P^n_K,*)}{z^r(\P^n_K,*)_{\sC,e}}
\]
is both injective and 0 on homology,whence the result in this case.\\
\\
{\bf Step 2. $X\subset \P^N$ projective.} This uses the method of the projecting cone. Suppose $X$ has dimension $n$ and let $L\subset \P^N$ be a dimension $N-n-1$ linear subspace with $L\cap X=\0$. Projection from $L$ defines a finite morphism
\[
\pi_L:X\to \P^n.
\]

For $Z$ a cycle on $X$, define
\[
\tilde{L}(Z):=\pi_L^*(\pi_{L*}(Z))-Z.
\]
Since $\pi_L$ is finite and $\P^n$ is smooth, $\tilde{L}$ is well-defined. In addition, $\tilde{L}(Z)\ge0$ if $Z\ge 0$.

For a closed subset $W\subset X_{sm}$, define $\codim_{X_{sm}}W$ to be the minimum of 
$\codim_{X_{sm}}Z$ as $Z$ runs over the irreducible components of $W$.

\begin{defn}
 For $A, B$ closed subsets of $X$ of pure codimension $a, b$, respectively, define the {\em excess} 
\[
e(A,B):= a+b-\codim_{X_{sm}}(A\cap B\cap X_{sm})
\]
if $A\cap B\cap X_{sm}=\0$, define the excess to be 0. If $Z$ and $W$ are cycles on $X$, define
\[
e(Z,W):=e(\text{supp}(Z), \text{supp}(W)).
\]
\end{defn}

Let $\Gr(N-n-1,N)$ be the Grassmannian of $\P^{N-n-1}$'s in $\P^N$, and let $\Gr(N-n-1,N)_X$ be the open subscheme of $L$ with $L\cap X=\0$. The classical moving lemma of Chow is

\begin{lem}\label{lem:ChowML1}  Let $Z$ and $W$ be cycles on $X$. Then there is a non-empty open subscheme $U_{Z,W}\subset \Gr(N-n-1,N)_X$ such that, for each field $K\supset k$ and each $K$-point $L$ of $U_{Z,W}$, 
\[
e(\tilde{L}(Z), W))\le \max(e(Z,W)-1, 0).
\]
\end{lem}

For each $L\in \Gr(N-n-1, N)_X$, we have the maps
\[
\pi_L\times\id:X\times\square^n\to \P^n\times\square^n.
\]
Define
\[
\tilde{L}_n:z^r(X,*)\to z^r(X,*)
\]
by
\[
\tilde{L}_n(Z):=\pi_L\times\id^*(\pi_L\times\id_*(Z))-Z
\]
The maps $\tilde{L}_n$ define the map of complexes
\[
\tilde{L}_*:z^r(X,*)\to z^r(X,*).
\]

With a little bit of work (see e.g. \cite{BlochWeb} or \cite{MixMot} for details) Chow's moving lemma implies

\begin{lem}\label{lem:ChowML2} Fix a finite set $\sC$ of locally closed subsets of $X_{sm}$ and a function $e:\sC\to\Z_+$. Define $e-1:\sC\to \Z_+$ by
\[
(e-1)(C):=\max(e(C)-1, 0).
\]
Let $K=k( \Gr(N-n-1,N))$ and let $L_\gen\in \Gr(N-n-1,N)(K)$ be the generic point of $\Gr(N-n-1,N)$. Then
\[
\tilde{L}_{\gen*}:z^r(X,*)\to z^r(X_K,*).
\]
maps $z^r(X,*)_{\sC,e}$ to $z^r(X_K,*)_{\sC,e-1}$.
\end{lem}

\begin{proof}[Sketch of proof] Let $|Z|\subset X\times\square^n$ be  the support of a cycle $Z$ in 
$z^r(X,n)$. For each  face $F$ of $\square^n$, let
\[
Z(F,d)=\{x\in X\ | \dim_k x\times F\cap |Z|\ge d\}
\]
where the inequality is satisfied if there is an irreducible component of $x\times F\cap Z$ of dimension $\ge d$ over $k$. $Z(F,d)$ is a constructible set, and $Z$ is in $z^r(X,n)_{\sC,e}$ if and only if
\begin{equation}\label{eqn:MLIequal}
\codim_CC\cap [Z(F,d)\setminus Z(F,d-1)]\ge r-e(C)-\dim_kF+d
\end{equation}
for each $C$, $F$ and $d$. We have a similar condition for membership in $z^r(X,n)_{\sC,e-1}$.

One has the operation $A\mapsto \tilde{L}(A)$ on constructible subsets $A$ of $X$ by setting
\[
\tilde{L}(A):=\ov{(\pi_L^{-1}(\pi_L(A))\setminus A)}
\]
where the closure is in $\pi_L^{-1}(\pi_L(A))$. One shows that
\[
\tilde{L}(Z)(F,d)\subset \tilde{L}(Z(F,d)).
\]
Using the moving lemma~\ref{lem:ChowML1}, one shows that the inequality \eqref{eqn:MLIequal}
 for the subsets $Z(F,d)$ implies the inequality  \eqref{eqn:MLIequal}
 for the subsets $ \tilde{L_\gen}(Z(F,d))$ with $e$ replaced by $e-1$, which proves the result.
 \end{proof}

To use this lemma, we proceed as follows: For each $C\in \sC$, we have the constructible subset $\pi_{L_\gen}(C)$. Write 
\[
\pi_{L_\gen}(C)=C'_1\cup\ldots\cup C'_r
\]
with each $C'_i$ locally closed in $\P^n_K$. Let $d_i=\codim_{\P^n}(C'_i)-\codim_XC$ and define
\[
e'(C'_i):=e(C)+d_i.
\]
Let $\sC'=\{C'_j\ |\ C\in \sC\}$.

Then $\pi_{L_\gen*}:z^r(X,*)\to z^r(\P^n_K,*)$ maps $z^r(X,*)_{\sC,e}$ to $z^r(\P^n_K,*)_{\sC',e'}$ and  $\pi_{L_\gen}^*:z^r(\P^n_K,*)\to z^r(X_K,*)$ maps $z^r(\P^n_K,*)_{\sC',e'}$
to  $z^r(X_K,*)_{\sC,e}$.

Lemma~\ref{lem:ChowML2} implies that
\[
\tilde{L}_*=
\pi_{L_\gen}^*\circ \pi_{L_\gen*}-\pi:\frac{z^r(X,*)_{\sC,e}}{z^r(X,*)_{\sC,e-1}}\to
\frac{z^r(X_K,*)_{\sC,e}}{z^r(X_K,*)_{\sC,e-1}}
\]
is zero, where $\pi$ is the base-extension map. Thus $\pi_{L_\gen}^*\circ \pi_{L_\gen*}$ is the same as base-extension on this quotient complex. However, we can factor $\pi_{L_\gen}^*\circ \pi_{L_\gen*}$  as
\[
\frac{z^r(X,*)_{\sC,e}}{z^r(X,*)_{\sC,e-1}}\xrightarrow{\pi_{L_\gen*}}
\frac{z^r(\P^n_K,*)_{\sC,e}}{z^r(\P^n_K,*)_{\sC,e-1}}
\xrightarrow{\pi_{L_\gen}^*}
\frac{z^r(X_K,*)_{\sC,e}}{z^r(X_K,*)_{\sC,e-1}}
\]

By Step 1, the complex in the middle is acyclic, so $\pi$ induces zero on homology. Since $K$ is a pure transcendental extension of $k$, the quotient complex ${z^r(X,*)_{\sC,e}}/{z^r(X,*)_{\sC,e-1}}$ is acyclic. Since $z^r(X,*)_{\sC,e-r}=z^r(X,*)$, the result follows by induction.

{\bf Step 3:} $X$ quasi-projective. This uses Bloch's moving lemma \cite{BlochMovLem} involving blowing up faces of $\square^n$.  We omit the details.

\subsection{A new moving lemma}\label{sec:NewMovingLemma}
To allow for an action of ``higher correspondences" on the additive cycle complexes with modulus, we will need an extension of Chow's moving lemma that takes into account the closures of the cycles. 

We have the open immersion $\square^n\subset {(\P^1)}^n$; write $\ov{\square}^n$ for ${(\P^1)}^n$. For an  integral subscheme or closed subset $Z$ of $X\times\square^n$, we let $\ov{Z}$ denote the closure in $X\times \ov{\square}^n$. Extending linearly, we have the operation
\[
\ov{(-)}:z_*(X\times\square^n)\to z_*(X\times\ov{\square}^n).
\]

\begin{defn}

For a set $\sC$ of locally closed subsets of $X$ and a function $e:\sC\to\Z_+$, let $\un{z}^r(X,n)_{\ov\sC,e}\subset \un{z}^r(X,n)_{\sC,e}$ be the subgroup generated by integral $Z\in  \un{z}^r(X,n)$ such that, for each face $F$ of $\square^n$, 
\[
\codim_{C\times \ov{F}}[\ov{Z\cap X\times F}]\cap C\times\ov{F}\ge r-e(C).
\]
\end{defn}

The $\un{z}^r(X,n)_{\ov\sC,e}$ form a cubical abelian group, hence we have the subcomplex
$z^r(X,*)_{\ov\sC,e}$ of $z^r(X,n)_{\sC,e}$.

\begin{lem}\label{lem:ML3} Let $X$ be a quasi-projective $k$-scheme, $\sC$ a finite set of locally closed subsets of $X_{sm}$, $e:\sC\to\Z_+$ a function. \\
\\
(1) If $X$ is projective,  the inclusion $z_r(X,*)_{\ov\sC,e}\to z_r(X,*)$ is a quasi-isomorphism.\\
\\
(2) In general, the inclusion $z_r(X,*)_{\ov\sC,e}\to z_r(X,*)$ induces an isomorphism on $H_0$.
\end{lem}

\begin{proof} The proof is essentially the same as theorem~\ref{thm:ML2}. We indicate the modifications.

\noindent
{\bf Step 1. $X=\P^n$} In modifying our earlier Step 1, since we need to keep track of closures, we use a finite sequence of ``generic translations" rather than writing the generic matrix as a product of elementary matrices.

In detail: Let $H\subset \P^n$ be the generic hyperplane, defined over a pure transcendental extension $K_1$ of $k$. Let $X_0,\ldots, X_n$ be homogeneous coordinates on $\P^n_{K_1}$ such that $H$ is defined by $X_0=0$. Let $x:=(x_1,\ldots, x_n)$ be the generic point of the affine space $\A^n_{K_1}:=\P^n\setminus H$, let $K=K_1(x_1,\ldots, x_n)$ and let $\phi(t)$ be the $K[t]$-valued point of $\GL_{n+1}$ defined linear transformation
\[
(X_0,X_1,\ldots, X_n)\mapsto (X_0, X_1+tx_1X_0,\ldots, X_n+tx_nX_0).
\]
This gives us the action $\phi:\A^1\times\P^n_K\to \P^n$.

For $Z$ a cycle on $\P^n\times\ov{\square}^n$ we have the cycle  $\phi^*(Z)$ on 
$\A^1\times\P^n_K\times\ov{\square}^n$. Let $\ov{\phi^*(Z)}$ be the closure of $\phi^*(Z)$ on
$\P^1\times\P^n_K\times\ov{\square}^n$. One can easily compute the intersection
$\ov{\phi^*(Z)}\cap \infty\times\P^n_K\times\ov{\square}^n$. This cycle is the join of $(x_1:\ldots:x_n)\times \ov\square^n\in H\times\ov\square^n$ with the the intersection $Z\cap H\times\ov\square^n$.

Since $H$ is generic over $k$, it is easy to check that sending $Z\in \un{z}_r(X,n)$ to $\phi^*(Z)\in 
\un{z}^r(X_K,n+1)$ defines a homotopy $H_\phi$ of $\phi(1)^*$ with the base-change map on the complex $z_r(X,*)$, and in addition $H_\phi$ restricts to give a degree one map
\[
H_\phi:\un{z}_r(X,n)_{\ov\sC,e}\to \un{z}_r(X,n+1)_{\ov\sC,e}
\]
giving a homotopy between the base-change on $\phi(1)^*$ on the subcomplex $z_r(X,*)_{\ov\sC,e}$.

Since the generic projective linear transformation is a composition of $n+1$ such translations, the argument we used for Step 1 above goes through to show prove the theorem for $X=\P^n$.

\noindent
{\bf Step 2. $X$ projective.}  We use exactly the same argument, replacing lemma~\ref{lem:ChowML2} with

\begin{lem}\label{lem:ChowML2.5} Fix a finite set $\sC$ of locally closed subsets of $X_{sm}$ and a function $e:\sC\to\Z_+$. Define $e-1:\sC\to \Z_+$ by
\[
(e-1)(C):=\max(e(C)-1, 0).
\]
Let $K=k( \Gr(N-n-1,N))$ and let $L_\gen\in \Gr(N-n-1,N)(K)$ be the generic point of $\Gr(N-n-1,N)$. Then
\[
\widetilde{L_\gen}_*:z_r(X,*)\to z_r(X_K,*).
\]
maps $z_r(X,*)_{\ov\sC,e}$ to $z_r(X_K,*)_{\ov\sC,e-1}$.
\end{lem}

The proof is exactly the same, where we replace the sets $Z(F,d)$ with
\[
\ov{Z}(F,d)=\{x\in X\ | \dim_k x\times \ov{F}\cap \ov{|Z|\cap X\times F}\ge d\}
\]

\noindent
{\bf Step 3.} $X$ quasi-projective.  Here  the ``moving by blow-ups" technique does not pass well to the closures $\ov\square^n$.

However, if $j:X\to \bar{X}$ is a projective closure of $X$, then let $\un{z}^r(X_{\bar{X}},*)_{\ov\sC,e}\subset \un{z}^r(X,*)_{\ov\sC,e}$ be the image of $j^*:\un{z}^r({\bar{X}},*)_{\ov\sC,e}\to
\un{z}^r(X,*)_{\ov\sC,e}$. If $Z$ is in $\un{z}^r(X,n)_{\ov\sC,e}$, then the closure $\ov{Z}$ is in
$\un{z}^r({\bar{X}},*)_{\ov\sC,e}$ for $n=0,1$, 

Thus, letting $i:W\to \bar{X}$ be the complement of $X$, we have a commutative diagram with exact rows
\[
\xymatrix{
\CH_r(W,0)\ar[r]^{i_*}\ar@{=}[d]&H_0(z_r({\bar{X}},*)_{\ov\sC,e})\ar[d]_\alpha\ar[r]^{j^*}&
H_0(z_r(X,*)_{\ov\sC,e})\ar[d]^\beta\ar[r]&0\\
\CH_{\dim X-r}(W,0)\ar[r]^{i_*}&\CH_r(\bar{X},0)\ar[r]_{j^*}&\CH_r(X,0)\ar[r]&0
}
\]
The map $\alpha$ is an isomorphism by Step 2, hence $\beta$ is an isomorphism.
\end{proof}

\section{Categories of motives}\label{sec:Motives}
\subsection{Chow motives over $S$}
Let $S$ be in $\Sm/k$. We recall the construction of a category of (homological) Chow motives over $S$.

We first assume that $S$ is irreducible. Let $\SmProj
/S$ denote the full subcategory of $\Sm/S$ consisting of the (smooth) projective $S$ schemes. For irreducible $X,Y\in\SmProj
/S$ set
\[
\Cor_S^n(X,Y):=\CH_{\dim_SX+n}(X\times_SY)
\]
Extend the definition to arbitrary $X,Y$ by taking the direct sum over the pairs of irreducible components of $X$ and $Y$; proceed similarly if $S$ is not irreducible.

Define the category $\Cor_S$ with objects $(X,n)$, $X\in\SmProj
/k$, $n\in\Z$, morphisms
\[
\Hom_{\Cor_S}((X,n),(Y,m)):=\Cor_S^{m-n}(X,Y)
\]
and composition law
 \[
 \beta\circ \alpha:=p^{XYZ}_{XZ*}(p^{XYZ*}_{YZ}(\beta)\cup p^{XYZ*}_{XY}(\alpha)).
 \]
 for $\alpha\in \Cor_S^*(X,Y)$, $\beta\in\Cor_S^*(Y,Z)$ and where, for instance, $p^{XYZ}_{XZ}:X\times_SY\times_SZ\to X\times_SZ$ is the projection.
 
The product over $S$ makes $\Cor_S$ a tensor category, with
\[
(X,n)\otimes(Y,m):=(X\times_SY,n+m)
\]
and, for $\alpha:(X,n)\to (Y,m)$, $\alpha':(X',n')\to(Y',m')$, 
\[
\alpha\otimes\alpha':= p_{XY}^{XX'YY'*}(\alpha)\cdot p_{X'Y'}^{XX'YY'*}(\alpha').
\]
The unit is $1:=(S,0)$. Note that
\[
\Hom_{\Cor_S}(1,(X,n))=\CH_n(X),
\]
so sending $(X,n)$ to $\CH_n(X)$ defines a functor
\[
\CH:\Cor_S\to\Ab.
\]
 
For each $r\in\Z$, we have the functor
\[
m_\Cor(r):\SmProj/S\to \Cor_S 
\]
sending $X$ to $(X,r)$ and $f:X\to Y$ to the class of the graph $\Gamma_f\subset X\times_SY$. We write $m_\Cor$ for $m_\Cor(0)$.

 The category of Chow motives over $S$, $\Mot(S)$, is the pseudo-abelian hull of $\Cor_S$, i.e., objects are $(X,n,\alpha)$ with $\alpha\in \End_{\Cor_S}(X)$ idempotent: $\alpha^2=\alpha$, and with
 \[
 \Hom_{\Mot(S)}((X,n,\alpha),(Y,m,\beta)):=\beta\Hom_{\Cor_S}((X,n),(Y,m))
\alpha
 \]
 
 $\Mot(S)$ is a tensor category with $(X,n,\alpha)\otimes(Y,m,\beta):=((X\times_SY,n+m,\alpha\otimes\beta)$. $\Mot(S)$ is a rigid tensor category with $(X,n,\alpha)^\vee:=(X,\dim_SX-n,{}^t\alpha)$; the unit $\delta:1\to (X,n,\alpha)\otimes (X,n,\alpha)^\vee$ and the co-unit $\epsilon: (X,n,\alpha)\otimes (X,n,\alpha)^\vee\to 1$ are both the image of the class of the diagonal
 $[\Delta]\in\CH^{\dim_kX}(X\times_SX)$.
 
 Sending $(X,n)$ to $(X,n,\id)$ defines a full tensor embedding $\Cor_S\to\Mot(S)$; we let
 \[
 m(r):\SmProj/S\to\Mot(S)
 \]
 be the composition of $m_\Cor(r)$ with this embedding. We write $m(X)(n)$ for the object $(X,n,\id)$.

\subsection{Highly distinguished complexes}

Central to Hanamura's construction of motives is the use of ``distinguished subcomplexes" of various cycle complexes. These are defined by properness conditions on various intersections with faces on $\square^n$. We refine these conditions to pass to the closed $n$-cube $\ov{\square}^n$, leading to the ``highly distinguished" subcomplexes.
 
\begin{defn}  Fix $X\in\SmProj/k$ and let  $\sW$ be a finite set of irreducible closed subsets $W_n\subset X\times T_n$, $n=1,\ldots, N$ for schemes $T_n\in \Sm/k$. For $F\subset\square^n$ a face, let
\[
p_{F,n}:X\times\ov{F}\times T_n\to X\times T_n
\]
be the projection. Let $\un{z}^r(X,*)_{\ov\sW}\subset \un{z}^r(X,*)$ be the subcomplex generated by irreducible $Z\subset X\times \square^n$ such that 
\begin{enumerate}
\item $Z$ is in $\un{z}^r(X,n)$
\item For each face $F$ of $\square^n$, $p^{-1}_{F,n}(W_n)$ and
$\ov{(Z\cap X\times F)}\times T_n$ intersect properly on $X\times \ov{F}\times T_n$.
\end{enumerate}
A subcomplex  $z^r(X,*)'\subset  z^r(X,*)$ that is the image of a subcomplex of the form
 $\un{z}^r(X,*)_{\ov\sW}$ is called {\em highly distinguished}.
\end{defn}

\begin{prop}\label{prop:disting1} The inclusion of a highly distinguished subcomplex
$z^r(X,*)'\subset  z^r(X,*)$ is a quasi-isomorphism.
\end{prop}

\begin{proof} For each $W_n\subset X\times T_n$, form the contructible subsets 
\[
C_{n,d}=\{x\in X\ |\  x\times T_n\cap W_n\text{ contains a component of dimension}\ge d\}
\]
Write each $C_{n,d}\setminus C_{n,d-1}$ as a union of irreducible locally closed subsets $C_{n,d}^j$. Set
\[
e(C_{n,d}^j):=-\dim W_n-d-\dim C_{n,d}^j
\]
Note that $e(C_{n,d}^j)\ge0$. Let $\sC:=\{C_{n,d}^j\}$. It is an easy exercise to show that
\[
\un{z}^r(X,*)_{\ov\sC,e}=\un{z}^r(X,*)_{\ov\sW}
\]
The proposition thus follows from lemma~\ref{lem:ML3}.
\end{proof}

\begin{ex}\label{ex:dist} Let $X, Y$ be in $\SmProj/k$. Suppose we are given irreducible closed subsets $W_i\subset X\times Y\times\square^{\ell_i}$, $i=1,\ldots, N$, which intersect all faces properly, i.e., each $W_i$ is the support of some cycle in $\un{z}^*(X\times Y,\ell_i)$. Let $\sW$ be the set of all irreducible components of all intersections $W_i\cap X\times Y\times F\subset X\times Y\times F$, as the $F$ run over all faces of $\square^{\ell_i}$. This gives us the distinguished subcomplex
$z^r(X,*)_{\ov{\sW}}$. Letting $W=\{W_1,\ldots, W_N\}$, Hanamura has defined the {\em distinguished subcomplex} $z^r_W(X,*)\subset z^r(X,*)$ (see \cite{HanamuraII}, p. 6); it follows directly from the definitions that
\[
z^r(X,*)_{\ov\sW}\subset z^r_W(X,*).
\]
From our proposition~\ref{prop:disting1} and \cite{HanamuraII}, Proposition 1.3, this inclusion is a quasi-isomorphism.
\end{ex}

In short:
\begin{lem}\label{lem:DistComp} Let $X$ be in $\SmProj/k$. Then each distinguished subcomplex of $z^r(X,*)$, in the sense of \cite{HanamuraII}, contains a highly distinguished subcomplex. \end{lem}

\subsection{The construction of $\DMH(k)$}

Hanamura considers three partially defined operations:
\begin{enumerate}
\item Take $X\in \SmProj/k$. The partially defined {\em cup product}
\[
\cup_X:z^r(X,*)\otimes z^s(X,*)\dasharrow z^{r+s}(X,*)
\]
is $\cup_X:=\delta^*_X\circ\boxtimes$, where $\delta_X:X\to X\times X$ is the diagonal.
\item Take $X, Y\in\SmProj/k$ and $f\in z^s(X\times Y,\ell)$. The partially defined  {\em push-forward by a correspondence $f$} is
\begin{align*}
f_*:z^r(X,n)&\dasharrow z^{r+s}(Y,n+\ell)\\
Z&\mapsto p_{Y*}(f\cup_{X\times Y}p_X^*(Z))
\end{align*}
\item Take $X,Y, Z\in\SmProj/k$. The partially defined {\em composition of correspondences} is
\begin{align*}
\circ:z^r(X\times Y,*)\otimes z^s(Y\times Z,*)&\dasharrow z^{r+s}(X\times Z,*)\\
u\otimes v&\mapsto v\circ u;\\
v\circ u:=p_{XZ*}(p_{YZ}^*(v)&\cup_{XYZ}p_{XY}^*(u))
\end{align*} 
where the maps $p_{XZ}$, etc., are the projections, $p_{XZ}:X\times Y\times Z\to X\times Z$, etc.
\end{enumerate}

In \cite{HanamuraII}, propositions 1.4, 1.5, Hanamura shows that, given a distinguished subcomplex $z(-)'$ for the target of one of these operations, there is a distinguished subcomplex of the domain for which the operation is well-defined and the image lands in the given distinguished subcomplex $z(-)'$. Exactly the same argument, with the help of the dictionary given in example~\ref{ex:dist}, proves the analog for the highly distinguished subcomplexes:

\begin{prop}[analog of \hbox{\cite{HanamuraII}, propositions 1.4, 1.5}] \label{prop:Composition} Take $X$, $Y$, $Z$ and  $W$ in $\SmProj/k$.\\
\\
(1a) Take $f\in z^s(X\times Y,\ell)$. There is a highly distinguished subcomplex 
$z^r(X,*)'$ such that $f_*$ is defined on $z^r(X,*)'$.\\
\\
(1b) Given further an element $g\in z^t(Y\times Z,m)$ such that $g\circ f$ is defined, there are highly distinguished subcomplexes $z^r(X,*)'$  and $z^{r+s}(Y,*)'$ such that $f_*$ and $(g\circ f)_*$ are defined on  $z^r(X,*)'$, $g_*$ is defined on  $z^{r+s}(Y,*)'$, $f_*(z^r(X,*)')\subset z^{r+s}(Y,*)'$ and
$g_*\circ f_*=(g\circ f)_*$ on $z^r(X,*)'$.\\
\\
(2) Take $w\in z^s(X,\ell)$. There is a highly distinguished subcomplex $z^r(X,*)'$ such that $w\cup_X(-)$ is defined on $z^r(X,*)'$.\\
\\
(3a) Take $v\in z^s(Y\times Z,\ell)$. There is a highly distinguished subcomplex $z^r(X\times Y,*)'$ such that $v\circ(-)$ is defined on $z^r(X\times Y,*)'$\\
\\
(3b) Given further an element $w\in z^t(Z\times W,m)$ such that $w\circ v$ is defined, there is a highly distinguished subcomplexes $z^r(X\times Y,*)'$ and $z^{r+s}(X\times Z,*)'$
 such that $v\circ(-)$ and $(w\circ v)\circ (-)$ are defined on $z^r(X\times Y,*)'$, $w\circ(-)$ is defined on 
 $z^{r+s}(X\times Z,*)'$, $v\circ z^r(X\times Y,*)'\subset z^{r+s}(X\times Z,*)'$ and
 $w\circ(v\circ (-))=(w\circ v)\circ (-)$ on $z^r(X\times Y,*)'$.\\
 \\
(4) These results generalize in the evident manner to any finite number of compositions of the
  operations $f_*$, $w\cup_X(-)$ and $v\circ(-)$.
  \end{prop}
  
  We follow Hanamura's construction of $\sD_{finite}(k)$ with three modifications:
 \begin{enumerate}
 \item We construct a homological category of motives rather than a cohomological one, in that $\DMH(k)$ comes with a functor
 \[
 m:\SmProj/k\to \DMH(k)
 \]
 rather than a functor from $\SmProj/k^\op$.
  \item We replace ``distinguished subcomplex" with ``highly distinguished subcomplex"
  \item We replace Hanamura's use of complexes of $\Q$-vector spaces $z^*(-,*)_\Q$ with the integral complexes $z^*(-,*)$. Structurally, the only price we pay is that the external product is now only commutative up to natural homotopy, rather than being strictly commutative. This means we need to
be a bit more careful in making sure that we use a consistent order in our products. This also prevents us from giving the integral construction a tensor structure.
\end{enumerate}
With these changes, one has the triangulated category  $\DMH(k)$. We briefly recall  the construction from \cite{HanamuraII} and show how to adapt it to our situation.

For integral $X, Y\in \SmProj/k$ and integers $r,s$, we set
\[
C^n((X,r), (Y,s)):=z^{r-s+\dim Y}(X\times Y, -n):=z_{s-r+\dim X}(X\times Y,-n).
\]
Clearly the differential in $z^{r-s+\dim Y}(X\times Y, *)$ makes $C^*((X,r), (Y,s))$ a cohomological complex; we write the differential in $C^*((X,r), (Y,s))$ as $\del$. We extend to finite formal sums of 
symbols $(X,r)$, $(Y,s)$ by making $C^*$ additive in each variable.

The objects of $\DMH(k)$ are built from {\em finite diagrams}: A finite diagram $(K^m, f^{n,m})$ consists of
\begin{enumerate}
\item finite formal symbols $K^m:=\sum_i(X_{im},r_{im})$, where each $X_{im}$ is in $\SmProj/k$, integral, and each $r_{im}$ is an integer. The sum is finite, and $K^m=0$ except for finitely many $m$.
\item correspondences $f^{n,m}=\sum_{ij}f_{ij}^{n,m}$ with $m<n$ and with
\[
f_{ij}^{n,m}\in C^{1+m-n}((X_{im}, f_{im}), (X_{jn}, r_{jn}))
\]
In addition, we assume that all compositions
\[
f^{m_n,m_{n-1}}\circ \ldots\circ f^{m_2, m_1}
\]
are defined.
\item The identity
\[
\del f^{n,m}+\sum_{m<\ell<n}f^{n,\ell}\circ f^{\ell,m}=0
\]
is satisfied for all $m<n$.
\end{enumerate}

Given two objects $(K,f)$, $(L,g)$, one defines a filtered system of quasi-isomorphic {\em function complexes} $\sHom(K,L)'$ as follows: Choose highly distinguished subcomplexes
\[
C^*(K^m,L^n)'\subset C^*(K^m,L^n)
\]
such that, for each $n,m$, $n',m'$ the maps
\begin{align*}
&(-)\circ f^{m,n}_{m'}:C^*(K^m,L^{m'})'[m-m']\to C^*(K^n,L^{m'})'[1+n-m']\\
&g^{m',n'}_n\circ(-):C^*(K^n, L^{n'})'[n-n']\to C^*(K^n,L^{m'})'[1+n-m']
\end{align*}
are defined. Let $C(K,M)'_p:=\oplus_{n-m=p}C^*(K^m,L^n)'$ and let
\[
h^{q,p}:C(K,M)'_p[-p]\to C(K,M)'_q[-q]
\]
be the degree +1 map $\sum_{m,n,m',n'}g^{m',m}_n\circ(-)-\epsilon\cdot(-)\circ f^{n,n'}_{m}$,
where the factor  $\epsilon$ is $(-1)^N$ on $(C(K,M)'_p[-p])^N$. We have the total complex $\Tot(C(K,M)')$ with
\[
\Tot(C(K,M)'):=\oplus_{p}C(K,M)'_p[-p]
\]
and differential $\sum_{p,q}\del_{C(K,M)'_p[-p]}+h^{q,p}$.  We set
\[
\sHom(K,L)':=\Tot(C(K,M)').
\]
This construction follows the formula for the category $C^b(\sA)$ of complexes over a d.g. category $\sA$, defined by Kapranov \cite{Kapranov}.

The complex $\sHom(K,L)'$ depends on the choice of the highly distinguished subcomplexes
$C^*(K^m,L^n)'$, but if we make another choice $C^*(K^m,L^n)''\subset C^*(K^m,L^n)'$, the induced inclusion map 
\[
\sHom(K,L)''\to \sHom(K,L)'
\]
 is a quasi-isomorphism. Since the intersection of two highly distinguised subcomplexes is again highly distinguished, this gives us a well-defined filtered system of quasi-isomorphic complexes $\sHom(K,L)'$. In particular, the cohomology $H^0(\sHom(K,L)')$ is canonically defined. We set
\[
\Hom(K,L):=H^0(\sHom(K,L)').
\]

The composition law $\Hom(K,L)\otimes\Hom(L,M)\to \Hom(K,M)$ is defined as follows: First fix an element $h_1\in \Hom(K,L)$. Choose a collection of highly distinguished subcomplexes $C^*(K^m,L^n)'$ and a representative $(h_1^{n,m})$ in $Z^0(\sHom(K,L)')$ for $h_1$. By 
 proposition~\ref{prop:Composition}, there is a collection of highly distinguished subcomplexes  $C^*(L^{m',}M^{n'})'$ and  $C^*(K^{m,}M^{n'})'$ for which $\sHom(L,M)'$ and $\sHom(K,M)'$ are defined and for which the operation
 \[
 (-)\circ (h_1^{n,m}):\sHom(L,M)'\to \sHom(K,M)'
 \]
 is defined. This gives us the operation
\[
 (-)\circ h_1:\Hom(L,M)\to \Hom(K,M),
 \]
which is independent of the choice of highly distinguished subcomplexes  $C^*(L^{m',}M^{n'})'$ and  $C^*(K^{m,}M^{n'})'$.

If we now change the representative $(h_1^{n,m})$ to $(h_1^{n,m})'$, there is an element $(H_1^{n,m})\in C^{-1}(\sHom(K,L)')$ with
\[
d(H_1^{n,m})=(h_1^{n,m})-(h_1^{n,m})'
\]
By replacing our choices  $C^*(L^{m',}M^{n'})'$ and  $C^*(K^{m,}M^{n'})'$ with smaller highly distinguished subcomplexes if necessary, we may assume that 
\[
(-)\circ *:\sHom(L,M)'\to \sHom(K,M)'
\]
is defined on the subcomplex of $\sHom(K,L)'$ generated by all irreducible components of
all the representatives $(h_1^{n,m})$, $(h_1^{n,m})'$ and $(H_1^{n,m})$. Thus the maps 
$(-)\circ (h_1^{n,m})$ and $(-)\circ (h_1^{n,m})'$ are homotopic, hence the composition $(-)\circ h_1$ is independent of the choice of representative for $h_1$ in $Z^0(\sHom(K,L)')$.

The proof that $(-) \circ h_1$ is independent of the choice of highly distinguished subcomplexes $C^*(K^m,L^n)'$ is similar as is the argument that the composition law is associative.  This gives us the additive category $\widetilde{\DMH}(k)$.

Using the construction of \cite{Kapranov} as extended by \cite{HanamuraII} one defines the cone of a morphism $h:K\to L$ in $\widetilde{\DMH}(k)$, or more precisely, the cone of a choice of representative of $h$ in some $Z^0(\sHom(K,L)')$. The  collection of cone sequences gives $\widetilde{\DMH}(k)$ the structure of a triangulated category; the argument of Hanamura copied word for word proves this result.

Finally, we let $\DMH(k)$ denote the pseudo-abelian hull of $\widetilde{\DMH}(k)$. The main theorem of 
\cite{BalmerSchlicht} says that  $\DMH(k)$ inherits a canonical structure of a triangulated category from 
$\widetilde{\DMH}(k)$.

Fix an integer $r$. We have the functor
\[
m(r):\SmProj/k\to \DMH(k)
\]
sending $X$ to the complex $K^0=(X,r)$, $K^m=0$ for $m\neq0$, and sending $f:Y\to X$ to the graph $\Gamma_f$ in $z^{\dim Y}(X\times Y)$.

\begin{lem} The functor $m(r)$ extends to a full embedding
\[
i:\Mot(k)\to  \DMH(k).
\]
\end{lem}

\begin{proof}  Since $\sHom((X,r), (Y,s))=z_{\dim X+s-r}(X\times Y,*)$ we have
\begin{align*}
\Hom_{\DMH(k)}((X,r), (Y,s))&=\CH_{\dim X+s-r}(X\times Y)\\
&=\Hom_{\Mot(k)}((X,r),(Y,s)),
\end{align*}
giving the full embedding $i:\Cor(k)\to \DMH(k)$. The extension to  $\Mot(k)$ follows since $\DMH(k)$ is pseudo-abelian.
\end{proof}

\subsection{The homological motive} 
Let $\Z\SmProj/k$ be the additive category generated by $\SmProj/k$: for $X, Y$ integral define
\[
\Hom_{\Z\SmProj/k}(X,Y):=\Z[\Hom_{\SmProj/k}(X,Y)]
\]
and extend to arbitrary $X$ and $Y$ by taking direct sums over the irreducible components. The composition law in $\Z\SmProj/k$ is induced from $\SmProj/k$.

Form the category of bounded complexes $C^b(\Z\SmProj/k)$ and the homotopy category $K^b(\Z\SmProj/k)$. We denote the complex concentrated in degree 0 associated to an $X\in\SmProj/k$ by $[X]$. Sending $X$ to $[X]$ defines the functor
\[
[-]:\SmProj/k\to C^b(\Z\SmProj/k)
\]

Let $i:Z\to X$ be a closed immersion in $\SmProj/k$, $\mu:X_Z\to X$ the blow-up of $X$ along $Z$ and $i_E:E\to X_Z$ the exceptional divisor with structure morphism $q:E\to Z$. Let $C(\mu)$ be the complex
\[
[E]\xrightarrow{(i_E,-q)}[X_Z]\oplus [Z]\xrightarrow{\mu+i}[X]
\]
with $[X]$ in degree 0. 
\begin{defn} The category $\sD_\hom(k)$ is the localization of the triangulated category $K^b(\Z\SmProj/k)$ with respect to the thick subcategory generated by the complexes $C(\mu)$.

Let 
\[
m_\hom:\SmProj/k\to \sD_\hom(k)
\]
be the composition
\[
\SmProj/k\xrightarrow{[-]}C^b(\SmProj/k)\to K^b(\SmProj/k)\to \sD_\hom(k).
\]
\end{defn}

\begin{lem}\label{lem:Descent} $\sD_\hom(k)$ with functor $m_\hom:\SmProj/k\to \sD_\hom(k)$ is a category of homological descent, in the sense of Guill\'en and Navarro Aznar \cite{GNA}.
\end{lem}

\begin{proof}  We use the notations of  \cite{GNA}. We first need to extend $m_\hom$ to a functor from the category of cubical objects of $\SmProj/k$ to $\sD_\hom(k)$. In fact, we can use the evident total complex functor to extend $[-]$ to a functor from cubical objects of $\SmProj/k$ to $C^b(\Z\SmProj/k)$. The remaining conditions are direct consequences of the definition of $\sD_\hom(k)$.
\end{proof}

We let $\Sch_k'$ denote the subcategory of proper morphisms in $\Sch_k$.

\begin{thm}\label{thm:Extension} Suppose that $k$ admits resolution of singularities. Then the functor $m_\hom$ extends to a functor
\[
M_\hom:\Sch_k'\to  \sD_\hom(k)
\]
such that\\
\\
1. Let $\mu:Y\to X$ be a proper morphism in $\Sch_k$, $i:Z\to X$ a closed immersion. Suppose that $\mu:\mu^{-1}(X\setminus Z)\to X\setminus Z$ is an isomorphism. There is a canonical extension of
the commutative square
\[
\xymatrix{
M_\hom(\mu^{-1}(Z))\ar[r]\ar[d]&M_\hom(Y)\ar[d]\\
M_\hom(Z)\ar[r]&M_\hom(X)
}
\]
to a map of distinguished triangles in $\sD_\hom(k)$
\[
\xymatrix{
M_\hom(\mu^{-1}(Z))\ar[r]\ar[d]&M_\hom(Y)\ar[d]\ar[r]&C_1\ar[r]\ar[d]^\alpha&M_\hom(\mu^{-1}(Z))[1]\ar[d]\\
M_\hom(Z)\ar[r]&M_\hom(X)\ar[r]&C_2\ar[r]&M_\hom(Z)[1]
}
\]
such that $\alpha$ is an isomorphism.\\
\\
2. Let $j:U\to X$ be an open immersion in $\Sch_k$ with closed complement $i:Z\to X$. Then there is a canonical distinguished triangle
\[
M_\hom(Z)\xrightarrow{i_*}M_\hom(X)\xrightarrow{j^*}M_\hom(U)\to M_\hom(Z)[1],
\]
natural with respect to proper morphisms of pairs $f:(X,U)\to (X',U')$.
\end{thm}

\begin{proof} It follows from lemma~\ref{lem:Descent} that the category $\sD_\hom(k)^\op$ with functor
\[
m_\hom^\op:\SmProj/k^\op\to \sD_\hom(k)^\op
\]
is a category of cohomological descent, in the sense of \cite{GNA}. By \cite{GNA}, th\'eor\`em 2.2.2. $m_\hom^\op$ extends to a functor (``cohomology with compact supports")
\[
M_{\hom}^c:\Sch_k^{\prime\op}\to  \sD_\hom(k)^\op
\]
satisfying (1) and  (2) with all arrows reversed. We take $M_\hom:=(M_{\hom}^c)^\op$.
\end{proof}

Fix an integer $r$. We have the functor
\[
m(r):\SmProj/k\to \DMH(k)
\]
Since the composition (as correspondences) of the graphs of composable morphisms is always defined, $m$ extends canonically to the exact functor
\[
K^b(m(r)):K^b(\Z\SmProj/k)\to \DMH(k)
\]
\begin{lem} $K^b(m)$ extends canonically to an exact functor
\[
\sD_\hom(m(r)):\sD_\hom(k)\to \DMH(k)
\]
\end{lem}

\begin{proof} We have already seen that the functors $m(r)$ extends to a functor
\[
i:\Mot_\CH(k)\to  \DMH(k).
\]
Let 
\[
\xymatrix{
E\ar[r]^{i_E}\ar[d]_q&X_Z\ar[d]^\mu\\
Z\ar[r]_i&X
}
\]
be a blow-up square in $\SmProj/k$. It is well-known that
\[
m(E)(r)\xrightarrow{(i_{E*},-q_*)}m(X_Z)(r)\oplus m(Z)(r)\xrightarrow{\mu_*+i_*}m(X)(r)
\]
is a split exact sequence in $\Mot_\CH(k)$. Thus $K^b(m(r))(C(\mu))\cong0$ in 
$\DMH(k)$, giving the desired extension to $\sD_\hom(k)$.
\end{proof}

\begin{defn} Let $bm(r):\Sch_k'\to  \DMH(k)$ be the composition
\[
\Sch_k'\xrightarrow{M_\hom}  \sD_\hom(k)\xrightarrow{\sD(m(r))}\DMH(k).
\]
We write $bm(X)(r)$ for $bm(r)(X)$, and call $bm(X)(r)$ the {\em twisted Borel-Moore motive of $X$}.
\end{defn}

\subsection{The category of extended motives}

\begin{defn} 
For $r\in\Z$, we have the full subcategory $\sC(r)\subset \DMH(k)$ with objects of the form
$\sD_\hom(m(r))(X)$, with $X\in \sD(k)$ and $r\in\Z$. Clearly the pseudo-abelian hull of $\sC(r)$ defines a full additve subcategory $\hat\Mot(k)_r$ of $\DMH(k)$. Let $\hat\Mot(k)$ be the full additive subcategory of $\DMH(k)$ generated by the $\hat\Mot(k)_r$, $r\in\Z$, and let $\hat\Mot(k)^*$ denote the smallest full subcategory of $\DMH(k)$ containing $\hat\Mot(k)$ and closed under translation.  We call  $\hat\Mot(k)^*$ the category of {\em extended motives over $k$}.
\end{defn}
Clearly the full embedding $i:\Mot(k)\to\DMH(k)$ factors through $\hat\Mot(k)$, and the functor $bm(r):
\Sch_k'\to \DMH(k)$ factors through the inclusion  $\hat\Mot(k)_r\to \DMH(k)$.

\subsection{Contravariant functoriality}

For $X\in\Sm/k$ equi-dimensional over $k$, set $h(X):=bm(X)(\dim X)$. Extend to arbitrary $X\in\Sm/k$ by taking the direct sum over the components of $X$.

We have the functor
\[
h_\SmProj:\SmProj/k^\op\to \Mot(k)
\]
with $h(X)$ defined as above and with $h(f:Y\to X)=[{}^t\Gamma_f]_*$. 

Let $f:X\to Y$ be a morphsim in $\Sm/k$. We proceed to define a morphism
\[
f^*:h(Y)\to h(X).
\]

By resolution of singularities, we may find projective completions $j:X\to\bar X$, $j':Y\to\bar Y$ such that\\
\\
1. the complements $D:=\bar X\setminus X$, $E:=\bar Y\setminus Y$ are strict normal crossing divisors.

Let $\bar\Gamma_f\subset \bar X\times\bar Y$ be the closure of the graph of $f$. Write $E=\sum_{j=1}^mE_j$ with each $E_j$ irreducible, and let $E_I=\cap_{i\in I}E_i$ for each subset $I\subset\{1,\ldots, m\}$. By \cite{BlochMovLem}, applied to the projection $\bar\Gamma_f\to\bar Y$, we assume that\\
\\
2. $\bar\Gamma_f$ intersects $\bar X\times E_I$ properly for each $I$. \\
\\
We say that a pair of projective completions  $j:X\to\bar X$, $j':Y\to\bar Y$ satisfying (1) and (2) are {\em good for $f$}.

Write $D=\sum_{i=1}^nD_i$ with each $D_i$ smooth , and let $D_I:=\cap_{i\in I}D_i$ for each subset $I\subset \{1,\ldots, n\}$. Then $m_\hom(X)$ is represented by the complex $D_*$:
\[
D_{1,\ldots,n}\xrightarrow{d_{n-1}}\oplus_{|I|=n-1}D_I\xrightarrow{d_{n-2}}\ldots \xrightarrow{d_1}\oplus_iD_i\xrightarrow{d_0} \bar X
\]
with $\bar X$ in degree 0, and where the differential is the signed sum of the inclusion maps
$D_I\to D_{I\setminus i}$, with sign $(-1)^{j-1}$ if $I=i_1<\ldots i_s$ and $i=i_j$.  We have a similar description of a complex $E_*$ representing $m_\hom(Y)$. We write $D_*(d_X)$  for the finite diagram \[
(\oplus_{|I|=j}(D_I,d_X), f^{j+1,j}), 
\]
where $f^{j+1,j}$ is the signed sum of the graphs of the corresponding inclusion maps. We define the finite diagram $E_*(d_Y)$ similarly. By definition of $h(X)$ and $h(Y)$, we have canonical isomorphisms in $\DMH(k)$
\[
D_*(d_X)\cong h(X),\ E_*(d_Y)\cong h(Y).
\]

We have the highly distinguished subcomplexes \[
C^*((E_J,d_Y), (D_I,d_X))'\subset C^*((E_J,d_Y),(D_I,d_X))
\]
defined by
\[
C^*((E_J,d_Y), (D_I,d_X))':=z_{d_X-|J|}(E_J\times D_I,-*)_{\ov\sC_{J,I}},
\]
where $\sC_{J,I}$ is the set of closed subsets $E_{J'}\times D_I$, $J\subset J'$; we include the cases $J=\0$ or $I=\0$, $E_\0=\bar Y$, $D_\0=\bar X$. The $C^*((E_J,d_Y), (D_I,d_X))'$ form a triple complex: the first differential is given by the differential in $C^*((E_J,d_Y), (D_I,d_X))'$, the second is the alternating sum of the push-forward maps
\[
C^*((E_J,d_Y), (D_I,d_X))'\xrightarrow{i_{J,I\setminus\{i\} i\subset I*}} C^*((E_J,d_Y), (D_{I\setminus \{i\}},d_X))'
\]
and the third is the alternating sum of the restriction maps
\[
C^*((E_J,d_Y), (D_I,d_X))'\xrightarrow{i^*_{J\cup\{j\},I\setminus i i\subset I}} C^*((E_{J\cup\{j\}},d_Y), (D_I,d_X))'
\]
The total complex of this triple complex is the internal Hom complex $\sHom(E_*(d_Y),D_*(d_X))'$ with respect to our choices of highly distinguished subcomplexes.

For each $j=1,\ldots, m$, the cycle $\bar\Gamma_f\cdot (\bar X\times E_j)$ is supported in $D\times E_j$, so we may write
\[
\bar\Gamma_f\cdot (\bar X\times E_j)=\sum_{i=1}^n \bar\Gamma_f(i,j).
\]
with $\bar\Gamma_f(i,j)$ supported in $D_i\times E_j$. This decomposition is not in general unique, as some components of $\bar\Gamma_f\cdot \bar X\times E_j$ may lie in $D_{ii'}\times E_j$.

It follows from our definitions that 
\begin{align*}
&{}^t\bar\Gamma_f\text{ is in } C^0((E_\0,d_Y),(D_\0,d_X),\\
&{}^t\bar\Gamma_f(i,j)\text{ is in } C^0((E_j,d_Y), (D_i,d_X)),\\
&{}^t\bar\Gamma_f+\sum_{i,j}{}^t\bar\Gamma_f(i,j)\text{ is in }Z^0(\sHom(E_*(d_Y),D_*(d_X))'). 
\end{align*}
We let $f^*_{\bar X,\bar Y}$ denote the class  
\[
\left[{}^t\bar\Gamma_f+\sum_{i,j}{}^t\bar\Gamma_f(i,j)\right]\in  H^0(\sHom(E_*(d_Y),D_*(d_X))'). 
\]
Having fixed $\bar X$ and $\bar Y$, the only choice we have made is the decomposition of the $\bar\Gamma_f\cdot (\bar X\times E_j)$ into pieces 
$ \bar\Gamma_f(i,j)$; as the difference of two such choices comes from the image of the differential
\begin{multline*}
z_{d_X-1}(E_j\times D_{ii'},-*)_{\ov\sC_{j,\{ii'\}}}\\\to
z_{d_X-1}(E_j\times D_i,-*)_{\ov\sC_{j,i}}\oplus
z_{d_X-1}(E_j\times D_i,-*)_{\ov\sC_{j,i'}}
\end{multline*}
the cohomology class $f^*_{\bar X,\bar Y}$ depends only on the choice of compactifications
$\bar X,\bar Y$.

By resolution of singularities, the collection of projective completions of $X$ and $Y$ that are good for $f$ form a filtered system. Suppose we have morphisms $a:\bar X'\to\bar X$, $b:\bar Y'\to\bar Y$ of completions. Let $\bar\Gamma_f'\subset \bar X'\times\bar Y'$ denote the closure of $\Gamma_f$. Then
\[
(a\times\id)_*(\bar\Gamma_f')=(\id\times b)^*(\bar\Gamma_f)
\]
as cycles on $\bar X\times\bar Y'$; this implies that, via the canonical isomorphisms
\[
D_*(d_X)\cong h(X)\cong D'_*(d_X),\ E_*(d_Y)\cong h(Y)\cong E'_*(d_Y),
\]
 the morphisms $f^*_{\bar X,\bar Y}$ and $f^*_{\bar X',\bar Y'}$ induce the same morphism $h(Y)\to h(X)$.
 
 Thus we have for each morphism $f:X\to Y$ in $\Sm/k$ a well-defined morphism
 \[
 f^*:h(Y)\to h(X)
 \]
 in $\DMH(k)$. Clearly $f^*$ is the morphism defined by ${}^t[\Gamma_f]$ in case $X$ and $Y$ are in $\SmProj/k$.
 
 Thus, it remains to prove the functoriality $(fg)^*=g^*f^*$. This does not follows by a direct computation, because the composition of the representatives we have used to define $f^*$ and $g^*$ may not be defined. We therefore proceed in stages. We require first a result from resolution of singularities:

\begin{lem}\label{lem:BlowUpTrans} Let $D$ be a reduced strict normal crossing divisor on some $X\in\Sm/k$ and let $i:W\to X$ be a closed immersion. Suppose that $W$ is irreducible, not contained in $D$ and that $W\setminus D$ is smooth. Then there is a sequence of blow-ups with irreducible smooth centers lying over $W\cap D$,
\[
X_N\to\ldots\to X_1\to X_0=X,
\]
such that
\begin{enumerate}
\item  the reduced pull-back  $D_i$ of $D$ to  $X_i$ is a strict normal crossing divisor for each $i$.
\item Let $Z_i\subset X_i$ be the smooth  irreducible closed subscheme blown up to form $X_{i+1}$, and let $F\subset D_i$ be the minimal face of $D_i$ containing $Z_i$. Then for each face $F'$ of $D_i$, $Z_i$ intersects $F\cap F'$ properly and transversely on $F$.
\item letting $W_N\subset X_N$ be the proper transform of $W$ and $D_N\subset X_N$ the reduced inverse image of $D$,  $W_N$ intersects each face of $D_N$ properly and transversely.
\end{enumerate}
\end{lem}

 We can now handle the case of the composition of closed immersions.
 
 \begin{lem}\label{lem:IncFunct} Let $i_1:X\to Y$, $i_2:Y\to Z$ be closed immersions. Then $(i_2i_1)^*=i_1^*i_2^*$.
 \end{lem}
 
 \begin{proof} We use the previous lemma to find projective completions $j':Y\to \bar Y$, $j'':Z\to\bar Z$ with strict normal crossing complements $D'\subset \bar Y$, $D''\subset \bar Z$ such that
 \begin{enumerate}
 \item $i_2$ extends to a closed immersion $\bar i_2:\bar Y\to \bar Z$
 \item $\bar Y$ intersects each face $D''_I$ properly and transversely
 \end{enumerate}
 Now let $X'\subset\bar Y$ be the closure of $X$. Then $X=X'\setminus D'$, and we may apply the previous lemma to the closed subscheme $X'$ of $\bar Y$. Let $T_i\subset \bar Y_i$ be the center of the $i$th blow-up $\bar Y_{i+1}\to \bar Y_i$ used to achieve the conclusion of the lemma.  We blow-up $\bar Z$, starting with the blow-up along $T_0$. Since each face of $D''$ intersects $\bar Y$ transversely, it follows that $T_0\subset D''$ satisfies the same transversality in lemma~\ref{lem:BlowUpTrans}(2) that $T_0$ satisfies with respect to $D'$. Thus $D''$ pulls back to a strict normal crossing divisor $D''_1$ on $\bar Z_1:=\bar Z_{T_0}$, intersecting the proper transform $\bar Y_1\subset \bar Z_1$ transversely. Thus we may continue the process, eventually replacing the closed immersion $\bar Y\to \bar Z$ with the closed immersion $\bar Y_n\to \bar Z_n$.
 
 Replacing $\bar Z$ and $\bar Y$ with $\bar Z_n$ and $\bar Y_n$, we may thus assume that we have a projective completion $j:X\to\bar X$ such that $i_1$ extends to $\bar i_1:\bar X\to \bar Y$, with 
 $\bar X$ transverse to each face of $D'$. Let $D=\bar X\setminus X$.
 
 We have the graphs $\Gamma_{\bar i_1}\subset \bar X\times\bar Y$, $\Gamma_{\bar i_2}\subset \bar Y\times\bar Z$ and $\Gamma_{\bar i_2\circ \bar i_1}$. By our conditions on $\bar i_1$ and $\bar i_2$, it follows that 
 $\Gamma_{\bar i_1}\cap \bar X\times D'_J$ is a disjoint union of components $D_I$ with $|I|=|J|$, and similarly for $\Gamma_{\bar i_2}\subset \bar Y\times\bar Z$ and $\Gamma_{\bar i_2\circ \bar i_1}$. This easily implies that the composition of correspondences $i^*_{1\bar X,\bar Y}\circ i^*_{2\bar Y,\bar Z}$ is defined and
 \[
 i^*_{1\bar X,\bar Y}\circ i^*_{2\bar Y,\bar Z}=
 (i_2\circ i_1)^*_{\bar X,\bar Z}
 \]
 As we are free to choose our projective completions to define the maps $i_1^*$, $i_2^*$ and $(i_2\circ i_1)^*$, this  gives $i_1^*i_2^*=(i_2i_1)^*$.
 \end{proof}
 
 \begin{lem}\label{lem:OpenExt} Let $g:U\to X$ be an open immersion with complement $X\setminus U$ a strict normal crossing divisor, $f:X\to Y$, $h:V\to U$ morphisms in $\Sm/k$. Then
 $(fg)^*=f^*g^*$ and $(gh)^*=h^*g^*$.
 \end{lem}
 
 \begin{proof} Let $j:X\to \bar X$ be a completion such that $\bar X\setminus U$ is a strict normal crossing divisor, and let $j':Y\to \bar Y$ be a completion such that $f^*_{\bar X,\bar Y}$ is defined. Taking $\bar U=\bar X$, $(fg)^*_{\bar X, \bar Y}$ is defined and $g^*_{\bar U,\bar X}$ is defined. As $\bar \Gamma_g$ is the diagonal in $\bar X$, the composition $g^*_{\bar U,\bar X}\circ f^*_{\bar X,\bar Y}$ is defined and 
 \[
 g^*_{\bar U,\bar X}\circ f^*_{\bar X,\bar Y}= 
 (fg)^*_{\bar X, \bar Y}.
 \]
 
 Similarly, let $V\to\bar V$ and $U\to \bar U$ be projective completions that are good for $h$. Set $\bar X=\bar U$. Then $\bar \Gamma_g$ is the diagonal in $\bar X$ and $(\bar V,\bar X)$ is good for $gh$, so $h^*_{\bar V,\bar U}$, $(gh)^*_{\bar V,\bar X}$ and $g^*_{\bar U,\bar X}$ are all defined. Since $\bar g=\id_{\bar U}$, the composition
 $g^*_{\bar U,\bar X}\circ h^*_{\bar V,\bar U}$ is defined and
 \[
 g^*_{\bar U,\bar X}\circ h^*_{\bar V,\bar U}=(gh)^*_{\bar V,\bar X}.
 \]
 \end{proof}

\begin{lem}\label{lem:ImProj} Given a  morphism $f:X\to Y$ in $\Sm/k$, factor $f$ as $p_2\circ i$, with
$p_2:X\times Y\to Y$ the projection and $i:X\to X\times Y$ the closed immersion $(\id_X,f)$. Then
$f^*=i^*\circ p_2^*$.
\end{lem}

\begin{proof} First choose completions $X\to\bar X$, $Y\to \bar Y$ such that $f$ extends to a morphism 
$\bar f:\bar X\to\bar Y$. Let $X':=\bar f^{-1}(Y)$. Then the restriction of $\bar f$ to $f':X'\to Y$ is a proper morphism, and we have the commutative diagram
\[
\xymatrix{
X\ar[r]^j\ar[d]_i\ar@/_30pt/[dd]|f&X'\ar[d]^{i'}\ar@/^30pt/[dd]|{f'}\\
X\times Y\ar[r]^{j\times\id}\ar[d]_{p_2}&X'\times Y\ar[d]^{p_2}\\
Y\ar@{=}[r]&Y}
\]
Suppose we have proven the result for $f':X'\to Y$. Using lemma~\ref{lem:OpenExt} we have
\begin{align*}
i^*\circ p_2^*&=i^*\circ(j\times\id)^*\circ p_2^*\\
&=(j\times\id\circ i)^*\circ p_2^*\\
&=(i'\circ j)^*\circ p_2^*\\
&=j^*\circ i^{\prime*}\circ p_2^*\\
&=j^*\circ f^{\prime*}\\
&=f^*
\end{align*}
Thus we may assume that $f$ is proper.

Choose completions $j:X\to \bar X$, $j':Y\to\bar Y$ for which the morphism $f_{\bar X,\bar Y}^*$ is defined. Let $\bar{i}:X\to \bar X\times Y$ be the composition  $(j\times\id)\circ i$ and let $\bar p_2:\bar X\times Y\to Y$ be the projection. Since $f$ is proper, $\bar i$ is a closed immersion. Using lemma~\ref{lem:OpenExt} again, it suffices to show
\[
f^*=\bar i^*\circ \bar p_2^*.
\]

The fact that $f^*_{\bar X,\bar Y}$ is defined implies that $\bar i^*_{\bar X,\bar X\times\bar Y}$ is also defined. Since $\bar X$ is smooth and projective, $\bar p^*_{2,\bar X\times\bar Y,\bar Y}$ is also defined, and the graph closure $\bar\Gamma_{\bar p_2}\subset \bar X\times\bar Y\times\bar Y$ is $\bar X\times\Delta_{\bar Y}$. The intersection 
\[
\bar\Gamma_{\bar i}\times\bar Y\cap \bar X\times\bar\Gamma_{\bar p_2}\subset
\bar X\times\bar X\times\bar Y\times\bar Y
\]
is $\delta_{\bar X\times\bar Y}(\bar\Gamma_f)$ (after rearranging the factors), hence the composition 
$\bar i^*_{\bar X,\bar X\times\bar Y}\circ \bar p^*_{2,\bar X\times\bar Y,\bar Y}$ is defined and  
\[
\bar i^*_{\bar X,\bar X\times\bar Y}\circ \bar p^*_{2,\bar X\times\bar Y,\bar Y}=
f^*_{\bar X,\bar Y}.
\]
\end{proof}

A slight variation:

\begin{lem}\label{lem:DoubleIm} Let $i:X\to Y\times T$ be a closed immersion, with $X$, $Y$ and $T$ in $\Sm/k$. Suppose that $p_2\circ i:X\to Y$ is a closed immersion. Then $(p_2\circ i)^*=i^*p_2^*$.
\end{lem}

\begin{proof} Using lemma~\ref{lem:OpenExt}, we may replace $T$ with a projective completion, i.e., we may assume $T$ is smooth and projective. If we take projective completions $\bar X$, $\bar Y$ good for 
$p_2\circ i$, then $\bar X$, $\bar Y\times T$ is good for $i$ and $\bar Y\times T$, $\bar Y$ is good for $p_2$. Also,  the composition $i^*_{\bar X,\bar Y\times T}\circ p^*_{2,\bar Y\times T,\bar Y}$ is defined and
\[
i^*_{\bar X,\bar Y\times T}\circ p^*_{2,\bar Y\times T,\bar Y}=
(p_2\circ i)^*_{\bar X,\bar Y}.
\]
\end{proof}

Next, we consider the case of a projection.
 \begin{lem}\label{lem:ProjFunct} For $X, Y\in \Sm/k$, let $p_2:X\times Y\to Y$ be the projection.  Let $f:V\to Z$ be a morphism in $\Sm/k$. Then $(fp_2)^*=p_2^*f^*$.
 \end{lem}
 
 \begin{proof} Let $j:X\to\bar X$ be a projective completion with strict normal crossing complement $\bar X\setminus X$. This gives us the commutative diagram
 \[
 \xymatrix{
 X\times Y\ar[d]_{p_2}\ar[dr]^{j\times\id}\\
 \bar X\times Y\ar[r]_{\bar p_2}&Y\ar[r]_f&Z}
 \]
 Using lemma~\ref{lem:OpenExt}, it suffices to show that $(f\bar p_2)^*=\bar p_2^*f^*$, so we may assume that $X$ is smooth and projective.
 
 Take projective completions $Y\to\bar Y$ and $Z\to\bar Z$ good for $f$. Since $p_2$ is smooth and projective, $X\times Y\to X\times\bar Y$, $Y\to\bar Y$ is good for $p_2$, and $X\times Y\to X\times\bar Y$, $Z\to\bar Z$ is good for $fp_2$. Similarly, the composition
 $p_{2 X\times\bar Y, \bar Y}^*\circ f^*_{\bar Y,\bar Z}$ is defined and
 \[
 p_{2 X\times\bar Y, \bar Y}^*\circ f^*_{\bar Y,\bar Z}=
 (fp_2)^*_{X\times\bar Y, \bar Z}.
 \]
 \end{proof}

Given an $X\in\Sm/k$ and a morphism $g:Y\to Z$ in $\Sm/k$, we have the closed immersions
\[
i_1:=\id_X\times (\id_Y,g):X\times Y\to X\times Y\times Z
\]
and
\[
i_2:=(\id_Y,g):Y\to Y\times Z
\]
giving the commutive diagram
\[
\xymatrix{
X\times Y\ar[r]^-{i_1}\ar[d]_{p_2}&X\times Y\times Z\ar[d]_{p_{23}}\ar[rd]^{p_3}\\
Y\ar[r]_-{i_2}&Y\times Z\ar[r]_{p_2}&Z}
\]

\begin{lem}\label{lem:Exchange} $p_2^*\circ i_2^*=i_1^*\circ p_{23}^*$ and $p_3^*=p_{23}^*\circ p_2^*$.
\end{lem}

\begin{proof} We first show that  $p_2^*\circ i_2^*=i_1^*\circ p_{23}^*$. Let $h:X\times Y\to Y\times Z$ be the composition $p_{23}\circ i_1$. By lemma~\ref{lem:ImProj} we have
\[
h^*=i_1^*\circ p_{23}^*,
\]
so we need only show that $h^*=p_2^*\circ i_2^*$. This follows from lemma~\ref{lem:ProjFunct}.

The identity $p_3^*=p_{23}^*\circ p_2^*$ also follows from   lemma~\ref{lem:ProjFunct}.
\end{proof}

\begin{thm}\label{thm:MotFunct} The assignment sending $X$ to $h(X)$ and $f:X\to Y$ to  $f^*:h(Y)\to h(X)$ defines a functor
\[
h:\Sm/k^\op\to \DMH(k)
\]
with $h_{|\SmProj}=h_\SmProj$.
\end{thm}

\begin{proof} We need to show that $(gf)^*=g^*f^*$ for morphisms $f:X\to Y$, $g:Y\to Z$ in $\Sm/k$. We have the commutative diagram
\[
\xymatrix{
X\ar[r]^{i_1}\ar[dr]_f\ar@/^20pt/[rr]|{i_3}&X\times Y\ar[d]^{p_2}\ar[r]^{\id\times i_2}&X\times Y\times Z\ar[d]^{p_{23}}\ar@/^30pt/[dd]|{p_3}\\
&Y\ar[r]^{i'}\ar[dr]_{g}&Y\times Z\ar[d]^{p_2}\\
&&Z}
\]
with $i_1=(\id_X,f)$, $i_2=(\id_Y,g)$, $i_3:=\id_X\times i_2\circ i_1$. Using lemmas~\ref{lem:IncFunct}, \ref{lem:ImProj}, \ref{lem:ProjFunct} and \ref{lem:Exchange}, we have
\[
f^*g^*=i_3^*\circ p_3^*
\]
Factoring $p_3$ as
\[
X\times Y\times Z\xrightarrow{p_{31}}X\times Z\xrightarrow{p_2} Z
\]
and using lemma~\ref{lem:ProjFunct} gives $p_3^*=p_{13}^*\circ p_2^*$. Thus
\[
f^*g^*=i_3^*\circ p_{13}^*\circ p_2^*.
\]
We have the commutative diagram
\[
\xymatrix{
&X\times Y\times Z\ar[d]^{p_{13}}\\
X\ar[ru]^-{i_3}\ar[r]_{i_4}\ar[rd]_{gf}&X\times Z\ar[d]^{p_2}\\
&Z}
\]
with $i_4:=(\id_X,gf)$. By lemma~\ref{lem:DoubleIm}
$i_4^*=i_3^*p_{13}^*$,
so by lemma~\ref{lem:ImProj} again
\[
f^*g^*=i_4^*p_2^*=(gf)^*.
\]
\end{proof}     
    
\section{Additive Chow groups} \label{sec:AdditiveChow}
\subsection{Normal schemes}

We prove some elementary results on normal varieties over a field.

\begin{lem} \label{lem:ProductNormality} 
Let $k$ be a perfect field.  and let $X$ and $Y$ be reduced normal finite type $k$-schemes. Then $X\times_kY$ is normal.
\end{lem}

\begin{proof}

 Since $k$   is perfect, the regular locus of $X$ or $Y$ is smooth over $k$, hence  the singular locus of  $X\times_kY$ is
\[
(X\times_kY)_{sing} =X_{sing}\times_kY\cup X\times_k Y_{sing}.
\]
In particular, $X\times_kY$ is smooth in codimension one.

Take a point $z\in (X\times_kY)_{sing}$ and let $x=p_1(z)$, $y=p_2(z)$.  Then either $x$   has codimension at least two on $X$ or  $y$   has codimension at least two on $Y$; we may suppose the former.

We recall that a local domain $\sO$ is normal if and only if $\Spec \sO$ is regular in codimension one and $\sO$ has depth at least two ({\it cf.} \cite{Matsumura}).

Since $X$ is normal, the local ring $\sO_{X,x}$ has depth at least two. Let $s_1, s_2$ be a regular sequence in the maximal ideal of $\sO_{X,x}$. Since $Y\to\Spec k$ is flat, $p_1^*(s_1), p_2^*(s_2)$ form a regular sequence in the maximal ideal of $\sO_{X\times_kY,z}$, hence $\sO_{X\times_kY,z}$ has depth at least two. Thus $X\times_kY$ is normal.
\end{proof}

\begin{lem}\label{lem:Effective} Let $k$ be a field, $f:Y\to X$ a projective,
surjective map of normal $k$-schemes of finite type. Let $D $ be  a 
Cartier divisor on $X$ such that $f^*D\ge 0$ on $X$. Then $D\ge0$ on $Y$. 
\end{lem}

\begin{proof}
Localizing at the set of generic points of $\text{supp}(D)$, we may assume 
that  $X=\Spec \sO$  with $\sO$ a DVR. Thus $Y\to X$ factors as a closed 
embedding $Y\to \P^N_\sO$ followed by the projection $\P^N_\sO\to X$. 

By repeatedly intersecting $Y$ with a hypersurface $H\subset \P^N_\sO$ of 
large degree and normalizing, we may assume that $f:Y\to X$ is generically 
finite. 

Let $Y\to Y'\to X$ be the Stein factorization of $f$. As $Y'$ is normal and 
$Y'\to X$ is finite and surjective, we must have $Y'=\Spec\sO'$, with $\sO'$ 
a semi-local PID. But $Y\to Y'$ is projective and birational, hence an 
isomorphism: $Y=\Spec\sO'$.

Let $t\in\sO$ be a generator of the maximal ideal $m$ of $\sO$, let 
$m'\subset \sO'$ be a maximal ideal  and let $s$ be a generator of $m'$. 
Since $\sO\to\sO'$ is finite, $m\sO'\subset m'$, hence $t=us^n$ for some unit 
$u$ in $\sO'_{m'}$ and some integer $n>0$. 

But $D$ is defined by $t^a$ for some $a\in\Z$, hence the restriction of 
$f^*D$ to $\Spec\sO'_{m'}$ is defined by $s^{an}$. Since $f^*D$ is effective, 
$an\ge 0$ hence $a\ge0$ and thus $D\ge0$.
\end{proof}

Let $X$ be a normal variety over $k$. Recall that a  Weil divisor on $D$ 
on $X$ is an element of the free abelian group $Div(X)$ on the set 
$PDiv(X)$ of prime divisors on $X$.
Thus $D$ is uniquely written as a finite sum $\sum a_Y [Y]$, where $Y$'s
are the prime divisors on $X$. Then the assignment $D \mapsto a_Y$ gives
a homomorphism $ord_Y : Div(X) \rightarrow {\Z}$, called the order function
of $Y$. A Weil divisor is called $effective$ if $ord_Y \ge 0$ for all
prime diviosrs $Y$ and one writes $D \ge 0$. The support of $D$ is the
set of all prime divisors $Y$ such that $ord_Y(D) \neq 0$. This is clearly
a finite set.

Let $Y_1, \cdots , Y_n$ be a set of Weil divisors on $X$. The $supremum$
of these divisors, denoted by ${sup}_{1 \le i \le n} Y_i$ is the Weil
divisor defined to be 
\begin{equation}\label{eqn:sup}
{sup}_{1 \le i \le n} Y_i : = \sum_{Y \in PDiv(X)} \left ({\rm max}_{1 \le i
\le n} \ ord_Y (Y_i) \right ) [Y].
\end{equation}  

\subsection{The additive cycle complex} \label{subsec:AdditiveCycleComplex}
We recall the definition of the additive cycle complexes with modulus from 
Park \cite{Park}. We fix a field $k$.

Set $\A^1:=\Spec k[t]$,  $\G_m:=\Spec k[t,t^{-1}]$, $\P^1:=\Proj k[Y_0,Y_1]$, 
and let $y:=Y_1/Y_0$ be the standard
coordinate function on $\P^1$. We set $\square^n:=(\P^1\setminus\{1\})^n$.

For $n \ge 1$, let $B_n = \A^1 \times \square^{n-1}$, 
$\ov{B}_n = \A^1 \times ({\P}^{1})^{n-1}\supset B_n$ and
$\widehat{B}_n = {\P}^{1} \times ({\P}^{1})^{n-1} \supset \ov{B}_n$. 
 We use the coordinate system 
$(t, y_1, \cdots , y_{n-1})$ on $\ov{B}_n$, with $y_i:=y\circ p_i$.

Let $F^1_{n,i}$, $i=1,\ldots, n-1$ be the Cartier  divisor
 on $\ov{B}_n$ defined by $y_i=1$ and $F_{n,0}\subset \ov{B}_n$
 the Cartier divisor defined by $t=0$.
 
 A {\em face} of $B_n$ is a subscheme $F$ defined by equations of the form
 \[
 y_{i_1}=\epsilon_1, \ldots,  y_{i_s}=\epsilon_s;\ \epsilon_j\in\{0,\infty\}.
 \]

For $\epsilon=0, \infty$,   $i = 1, 
\cdots, n-1$ let 
\[
\iota_{n,i,\epsilon}: B_{n-1}\to B_n
\]
be the inclusion
\[
\iota_{n,i,\epsilon}(t,y_1,\ldots, y_{n-2})=(t,y_1,\ldots, y_{i-1}, 
\epsilon,y_i,\ldots, y_{n-2}),
\]

\begin{defn}[Park\cite{Park}] \label{def:AdditiveComplex}
Let $X$ be a finite type $k$-scheme, $r,m$
nonnegative integers $m\ge1$. \\
\\
(0) $\un{\TZ}_r(X, 1; m)$  is the free abelian group 
on integral closed subschemes $Z$ of $X \times \A^1$ of dimension $r$ 
such that $Z\cap (X \times \{0\}) \time0=\0$.\\
\\
For $n>1$,  $\un{\TZ}_r(X, n; m)$ is the free abelian group 
on integral closed subschemes $Z$ of $X \times B_n$ of dimension $r+n-1$ 
such that: \\ \\
$(1)$ (Good position) For each face $F$ of $B_n$, $Z$ intersects  $X \times F$
properly:
\[
\dim_{X \times_k F} (Z\cap (X\times_k F)) \le r+\dim_k (F),
\] 
$(2)$ (Modulus condition) Let $\ov{Z}$ be the closure of $Z$ in $X \times 
\ov {B}_n$,   $\bar{p}:{\ov Z}^N \rightarrow {\ov Z}$ the normalization, and 
$p:
{\ov Z}^N \rightarrow X\times\ov{B}_n$ the evident map.Then 
\[
(m+1)\cdot p^*(X\times F_{n,0})\le 
\text{sup}_{i=1,\ldots n-1}p^*(X\times F_{n,i}^1) 
\]
(3) (Induction) For each $i=1,\ldots, n-1$, $\epsilon=0,\infty$, 
$(\id_X\times\iota_{n,i\epsilon})^*(Z)$ is in
$\un{\TZ}_r(X, n-1; m)$.
 \end{defn}
Note that the good position condition on $Z$ implies that the cycle 
$(\id_X\times\iota_{n,i\epsilon})^*(Z)$
is well-defined and each component satisfies the good position condition. If 
$X$ is equi-dimensional of dimension $d$ over $k$, set
\[
\un{\TZ}^r(X, n; m)=\un{\TZ}_{d+1-r}(X, n; m).
\]

We thus have the cubical abelian group $\un{n}\mapsto \un{\TZ}_{r}(X, n; m)$, 
and if $X$ is equi-dimensional over $k$, the cubical abelian group 
$\un{n}\mapsto \un{\TZ}^{r}(X, n; m)$.

\begin{rem} We have changed Park's notation a bit; our $m$ corresponds to his $m+1$. This fits with the convention in \cite{Ruelling}.
\end{rem}

\begin{defn}
$\TZ_r(X, *; m)$ is the {\em additive cycle complex} of
$X$ in dimension $r$ and with modulus $m$ which in degree $n$:
$$\TZ_r(X, n; m): = \frac{\un{\TZ}_{r}(X, n; m)}{\un{\TZ}_{r}(X, n; m)_{\dgn}}$$

\[
\TH_r(X, n; m): = H_n (\TZ_r(X, *; m)); \ n \ge 1
\]
is the 
{\it additive higher Chow group}  of $X$ with modulus $m$. 
\end{defn}

\subsection{Projective push-forward}

For a morphism $f:X\to Y$, we have the induced morphism
\[
f_n:=f\times\id_{B_n}:X\times B_n\to Y\times B_n.
\]
We write $\bar{f}_n$ for the extension $f\times\id_{\ov{B}_n}$.

\begin{lem}\label{lem:AddPushforward} Let $f : X \to Y$ be a projective morphism of finite type $k$-schemes. Then
\[
Z\in \un\TZ_r(X, n; m)\Longrightarrow
f_{n*}(Z)\in \un\TZ_r(Y, n; m).
\]
\end{lem}
\begin{proof}

We may assume that $Z$ is an integral cycle and
that   $Z\to f_n(Z)$ is generically finite. Since $f_n$ is proper, we have 
\[
f_n(Z)\cap Y\times F=f_n(Z\cap X\times F)
\]
(as closed subsets) for each face $F$. Thus
$f_{n*}(Z)$ is in good position. 

We now show that $f_{n*}(Z)$ satisfies 
the modulus condition. Let $\ov{Z}$ be the closure of $Z$ in $X \times 
{\ov{B}_n}$ and let $p : \ov{Z}^N \rightarrow \ov{Z} \hookrightarrow X \times {\ov{B}_n}$ be 
the normalization map for $\ov{Z}$. Put $V = f_n(W)$. Let $\ov{V}$ be the closure
of $V$ in $Y \times {\ov B_n}$ and let $q : \ov{V}^N \rightarrow \ov{V}\hookrightarrow Y 
\times {\ov B_n}$ be the normalization map. We have a commutative 
diagram
$$
\xymatrix{
\ov{Z}^N \ar[d]_{h} \ar[r]^{p} &  X \times {\ov B_n} \ar[d]^{f \times id} \\
\ov{V}^N \ar[r]_{q} & Y \times {\ov B_n} \\
}
$$ 
Thus
\begin{align*}
&h^*q^*(Y\times F_{n,0})=p^*(X\times F_{n,0})\\
&h^*q^*(Y\times F_{n,i}^1)=p^*(X\times F_{n,i}^1)
\end{align*}
Therefore the modulus condition for $Z$  and Lemma~\ref{lem:Effective} imply the modulus condition for $V$.

Suppose $n=0$. Since $f$ is proper, we have  
$f(X\times 0\cap \ov{Z})=\bar{f}_0(\ov{Z})\cap(Y\times0)$. Thus the condition (0) for $Z$ implies the same for $f_{0*}(Z)$.

For $n>0$, we have
\[
\partial_i^\epsilon(f_{n*}(Z))=f_{n*}(\partial_i^\epsilon(Z)).
\]
Thus the inductive condition for $Z$ (and induction) implies the same for $f_{n*}(Z)$
\end{proof}

Coming back to the construction of push-forward map, Lemma~\ref{lem:AddPushforward} shows 
that the maps $f_{n*}$ gives a well defined map $f_{n*}:\TZ
_r(X, n ; m) \rightarrow 
\TZ_r(Y, n ; m)$. These maps clearly commute with the boundary maps, 
as  in the case of higher Chow cycles. Thus we get a  map of 
complexes 
\[
f_* : \TZ_r(X, * ; m) \rightarrow \TZ_r(Y, * ; m)
\]
and hence   maps $f_ * : \TH_r(X, n ; m) \rightarrow \TH_r(Y, n ; m)$.

The functoriality
\[
(fg)_*=f_*g_*
\]
follows from the functoriality on the level of cycles.

\subsection{Flat pull-back}
Let $g:T\to S$ be a flat morphism of schemes. Recall that $f^*:z^i(S)\to z^i(T)$ is defined on generators by
\[
f^*(1\cdot Z):=[f^{-1}(Z)],
\]
where $f^{-1}(Z):=T\times_SZ$, considered as a closed subscheme of $T$ and for a closed subscheme $W\subset T$ of pure codimension $i$, $[W]\in z^i(T)$ is the associated cycle:
\[
[W]:=\sum_{w\in W_{gen}}\ell_{\sO_{T,w}}(\sO_{W,w})[w].
\]

Let $f : X \to Y$ be a flat morphism of  finite type $k$-schemes, of relative dimension $d$. The maps $f_n:X\times B_n\to Y\times B_n$ are therefore flat as well, hence we have a well-defined pull-back map
\[
f_n^*:z_r(Y\times B_n)\to z_{r+d}(X\times B_n).
\]

\begin{lem} For $V\in \un\TZ_r(Y, n ; m)$, $f_n^*(V)$ is in $\un\TZ_{r+d}(X, n ; m)$.
\end{lem}

\begin{proof}
 The good position condition follows immediately from the flatness of $f_n$.

To check the modulus condition, let $V\in \un\TZ_r(Y, n ; m)$ be an integral cycle, $\ov{V}$ the closure in $Y\times\ov{B}_n$, and let $Z=f^{-1}_n(V)_\red$.  We may assume that $Z$ is non-empty.

Since $\bar{f}_n$ is flat, we have (as closed subsets)
\[
\ov{Z}=\bar{f}_n^{-1}(\ov{V}).
\]
In particular, the morphism $\ov{Z}\to\ov{V}$ induced by $\bar{f}_n$ is equi-dimensional and dominant. Thus the map on the normalizations $\ov{Z}^N\to \ov{V}^N$ is also equi-dimensional and for each irreducible component $\ov{Z}^N_\alpha$, the restriction $q_\alpha:\ov{Z}^N_\alpha\to \ov{V}^N$ is  dominant.

We thus have the commutative diagram, where all the vertical maps are equi-dimensional and 
dominant.
$$
\xymatrix{
\ov{Z}_\alpha^N  \ar[d]_{q_\alpha} \ar[r]^-{p_{Z_\alpha}} & X \times \ov{B}_n \ar[d]^{\bar{f}_n} \\
 \ov{V}^N \ar[r]_-{p_V}& Y \times \ov{B}_n
}
$$
This gives the identities
\begin{align*}
&p_{Z_\alpha}^*(X\times F_{n,0})=q_\alpha^*(p_V^*(Y\times F_{n,0}))\\
&p_{Z_\alpha}^*(X\times F_{n,i}^1)=q_\alpha^*(p_V^*(Y\times F_{n,i}^1))
\end{align*}
which show that the modulus condition for $V$  implies the same for $Z$.

If $n=0$, the identity $\ov{Z}=\bar{f}_0^{-1}(\ov{V})$ shows that the condition (0) for $V$ implies
the condition (0) for $Z$.

The flatness of $f_n$ implies
\[
\partial^\epsilon_{n,i}(f_n^*V)=f_{n-1}^*(\partial^\epsilon_{n,i}(V))
\]
which yields the inductive condition.
\end{proof}

The maps $f_n^*$ thus gives us the map of complexes
\[
f^*: \TZ^i(Y, * ; m)\to  \TZ^i(X, * ; m).
\]
The functoriality of flat pull-bach for cycles yields the functoriality
\[
(fg)^*=g^*f^*
\] 
as maps of complexes, for composible flat morphisms $f$ and $g$.

The compatibility of projective push-forward and flat pull-back of cycles in cartesian squares yields:
\begin{lem}\label{pp}
Let 
$$
\xymatrix{
X' \ar[d]^{f'} \ar[r]^{g'} & X \ar[d]^{f} \\
Y' \ar[r]^{g} & Y \\
}
$$
be a cartesian square of maps of finite type $k$-schemes, where $f$ is 
flat and $g$ is projective. Then  
$$f^* g_* (\alpha) = g'_* f'^* (\alpha)$$
as maps of complexes $\TZ^i(Y', * ; m)\to  \TZ^i(X, * ; m)$.
\end{lem}

\section{Higher Chow groups and additive Chow groups}\label{sec:ChowAction}
Our aim in this section is to show that the higher Chow groups of smooth 
quasi-projective varieties act on the degenerate Chow groups in a natural
way. We suppose in this section that $k$ is a perfect field.

\subsection{External product}
We begin with the construction of an external product.

\begin{defn} Let $X$ and $Y$ be finite type $k$-schemes. Define
\[
\boxtimes_{a,b}:\un{z}_r(X,a)\otimes\un\TZ_s(Y,b;m)\to
z_{r+s}(X\times_kY\times B_{a+b})
\]
by
\[
Z\boxtimes W:=\tau_*(Z\times W),
\]
where 
\[
\tau:X\times\square^a\times Y\times\A^1\times \square^b\to
X\times Y\times\A^1\times\square^a\times\square^b
\]
is the exchange of factors.
\end{defn}

\begin{lem} For $Z\in \un{z}_r(X,a)$ and $W\in  \un\TZ_s(Y,b;m)$, $Z\boxtimes W$ is
in $\un\TZ_{r+s}(X\times_kY,a+b;m)$.
\end{lem}

\begin{proof} We may assume that $Z$ and $W$ are integral closed subschemes of $X\times\square^a$ and $Y\times B_b$, respectively.

Since 
\[
Z\times W\cap X\times F\times Y\times F'=(Z\cap X\times F)\times(W\cap Y\times F'),
\]
for faces $F$ of $\square^a$ and $F'$ of $B_b$, the good position condition is clear.

The condition (0) in case $a=b=0$ is checked the same way.

For the modulus condition, let $\ov{Z}\subset X\times{\P^1}^a$, $\ov{W}\subset Y\times\ov{B}_n$ be the respective closures of $Z$ and $W$, and let
\begin{align*}
p_Z:\ov{Z}^N\to X\times{\P^1}^a\\
p_W:\ov{W}^N\to Y\times\ov{B}_n
\end{align*}
be the normalizations of $\ov{Z}$ and $\ov{W}$. By lemma~\ref{lem:ProductNormality} 
\[
p_{Z\times W}:\ov{Z}^N\times_k\ov{W}^N\to X\times{\P^1}^a\times Y\times\ov{B}_n
\]
is the normalization of $\ov{Z}\times\ov{W}=\ov{Z\times W}$.

Suppose that $j$ is an index with
\[
m\cdot p_W^*(X\times F_{b,0})\le p_W^*(X\times F_{b,j}^1).
\]
Then
\[
(m+1)\cdot(\ov{Z}^N\times p_W^*(X\times F_{b,0}))\le \ov{Z}^N\times p_W^*(X\times F_{b,j}^1);
\]
since
\begin{align*}
&p_{Z\times W}^*(\tau^*(X\times Y\times F_{a+b,0}))=\ov{Z}^N\times p_W^*(X\times F_{b,0})\\
&p_{Z\times W}^*(\tau^*(X\times Y\times F^1_{a+b,j})=\ov{Z}^N\times p_W^*(X\times F_{b,j}^1),
\end{align*}
the modulus condition for $Z\boxtimes W$ is satisfied.

We finish by checking the inductive condition. If $b=0$, then $\ov{Z\boxtimes W}\cap X\times F_{a,0}=\0$, so the inductive condition is trivially satisfied. 

Suppose that $a+b>0$. For $1\le i\le a$, we have
\[
\partial_{a+b,i}^\epsilon(Z\boxtimes W)=\partial^\epsilon_{a,i}(Z)\boxtimes W
\]
and for $a+1\le i\le a+b$, we have
\[
\partial_{a+b,i}^\epsilon(Z\boxtimes W)=Z\boxtimes \partial^\epsilon_{b,i-a}(W).
\]
In the second case, the inductive condition for $W$, plus induction on $b$, shows that $Z\boxtimes \partial^\epsilon_{b,i-a}(W)$ is in $\un\TZ_{r+s}(X\times_kY,a+b-1;m)$. In particular, if $a=0$, then 
$Z\boxtimes W$ is in $\un\TZ_{r+s}(X\times_kY,a+b;m)$. By induction on $a$, 
$\partial^\epsilon_{a,i}(Z)\boxtimes W$ is in $\un\TZ_{r+s}(X\times_kY,a+b-1;m)$ as well, which completes the proof.
\end{proof}

\begin{thm}\label{thm:nolabel} Let $X$ and $Y$ be finite type $k$-schemes. Then the maps
\[
\boxtimes_{a,b}:\un{z}_r(X,a)\otimes\un\TZ_s(Y,b;m)\to
\un\TZ_{r+s}(X\times_kY,*;m)
\]
define a map of complexes
\[
\boxtimes:z_r(X,*)\otimes\TZ_s(Y,*;m)\to
\TZ_{r+s}(X\times_kY,*;m)
\]
\end{thm}

\begin{proof}  A formal computation shows that the product maps $\boxtimes_{a,b}$ are compatible with boundaries, and that
\begin{multline*}
\boxtimes_{a,b}(\un{z}_r(X,a)_\dgn\otimes\un\TZ_s(Y,b;m)+\un{z}_r(X,a)\otimes\un\TZ_s(Y,b;m)_\dgn)\\\subset
\un\TZ_{r+s}(X\times_kY,*;m)_\dgn.
\end{multline*}
\end{proof}

\begin{rem} It is easy to check that $\boxtimes$ is compatible with flat pull-back: for $f:X'\to X$, $g:Y'\to Y$ flat, we have
\[
f^*(Z)\boxtimes g^*(W)=(f\times g)^*(Z\boxtimes W)
\]
and projective push-forward:  if $f:X'\to X$  and $g:Y'\to Y$ are projective then
\[
f_*(Z)\boxtimes g_*(W)=(f\times g)_*(Z\boxtimes W).
\]
Furthermore we have the associativity
\[
Z\boxtimes(W\boxtimes V)=(Z\boxtimes W)\boxtimes V
\]
\end{rem}

Passing to homology, we thus have the external product
\[
\boxtimes:\CH_r(X,a)\otimes \TH_s(Y,b)\to \TH_{r+s}(X\times_kY,a+b)
\]

\subsection{Cap product}
We refine the external product to a ``cap product"
\[
\cap_X:z^r(X,*)\otimes \TZ^s(X,*;m)\to \TZ^{s+r}(X,*;m)
\]
for $X\in\Sm/k$. Let $\delta_X:X\to X\times_XX$ be the diagonal, and set 
\[
\delta_{X,n}:=\delta_X\times\id:X\times B_n\to X\times_kX\times B_n.
\]

\begin{defn} For $X\in\Sm/k$, let 
\[
\un\TZ^r(X\times_kX,n;m)_\Delta\subset\un\TZ^r(X\times_kX,n;m)
\]
be the subgroup generated by integral closed subschemes $Z\subset X\times X\times B_n$
such that\\
\\
1. $Z$ is in  $\un\TZ^r(X\times_kX,n;m)$\\
2. $\codim_{X\times F}\left (\delta^{-1}_{X,n}(Z)\cap X\times F \right )\ge r$  for all faces $F$ of $B_n$.\\
3. $\delta^*_{X,n}(Z)$ is in $\TZ^r(X,n;m)$.
\end{defn}
The subgroups $\un\TZ^r(X\times_kX,n;m)_\Delta$ form a cubical subgroup 
 \[
 \un\TZ^r(X\times_kX,-;m)_\Delta\subset \un\TZ^r(X\times_kX,-;m)
 \]
 and the maps $\delta^*_{X,n}$ give a well-defined map of complexes
 \[
 \delta_X^*:\TZ^r(X\times_kX,*;m)_\Delta\to
 \TZ^r(X,*;m).
 \]

Let $\sW$ be a finite set of closed subsets $W_n\subset X\times B_n$, $n=1,\ldots, N$. $\sW$ generates a ``cubical closed subset" $\sW^c$ by 
\[
\sW^c(\un{n}):=\cup_{\substack{g:\un{n}\to\un{m}\\ 1\le m\le N}}g^{-1}(W_m)\subset \square^n.
\]

\begin{defn} Let $\sW$ be a finite set of closed subsets $W_n\subset X\times B_n$, such that each $W_n$ is the support of a cycle in $\un\TZ_s(X,n;m)$. Let 
$\un\TZ_s^{\sW}(X,n;m)\subset \un\TZ_s(X,n;m)$ be the subgroup of cycles $Z$ with
\[
\text{supp}(Z)\subset \sW^c(\un{n}).
\]
Similarly, $\un{z}^r_\sW(X,n)\subset \un{z}^r(X,n)$ is the subgroup of cycles $Z$ such that
$\text{supp}(Z)\subset \sW^c(\un{n})$.
\end{defn}
From the construction of $\sW^c$, it is immediate that the $\un\TZ_s^\sW(X,n;m)$ form a cubical subgroup
$\un\TZ_s^\sW(X,-;m)$ of $\un\TZ_s(X,-;m)$, giving us the subcomplex
\[
\TZ_s^\sW(X,*;m):=\frac{\un\TZ^r_\sW(X,n;m)}{\un\TZ^r_\sW(X,n;m)_\dgn}\subset
\TZ_s(X,*;m).
\]
Similarly, we have the subcomplex $z^r_\sW(X,*)$ of $z^r(X,*)$.

The construction of the cap product is based on the following result:

\begin{lem}\label{lem:ProductML} Fix an integer $s$. Let $\sW$ be a finite set of closed subsets $W_n\subset X\times B_n$, $n=1,\ldots,N$, such that each $W_n$ is the support of a cycle in $\TZ_s(X,n;m)$. Let $f:X\to Y$ be a smooth morphism. Then there is a finite set of locally closed subsets $\sC$   such that $(A,B)\mapsto f^*(A)\boxtimes B$ restricts to a well-defined map of complexes
\[
f^*(-)\boxtimes:z^r(Y,*)_{\ov\sC}\otimes \TZ_s^\sW(X,*;m)\to
\TZ_{\dim_kX+s-r}(X\times_kX,*;m)_\Delta.
\]
for all $r$.
\end{lem}

\begin{proof} It is obvious that $\un\TZ_s^\sW(X,n;m)= \un\TZ_s^\sW(X,n;m)_\dgn$ for all $n>N$, hence 
$\TZ_s^\sW(X,n;m)=0$ for $n>N$.

If $F$ is a closed subset of some $T\in\Sch_k$, we write $\dim_kF\ge d$ if there is an irreducible component $F_\alpha$ of $F$ with $\dim_kF_\alpha\ge d$.

For each $n=0,\ldots, N$, define the constructible subsets $C_{n,d}$ of $Y$ by
\[
C_{n,d}=\{y\in Y\ | \dim_k f^{-1}(y)\times B_n\cap \sW^c(\un{n})\ge d\}
\]
Define the constructible subsets $C^0_{n,d}$ of $Y$ by
\[
C^0_{n,d}=\{y\in Y\ | \dim_k f^{-1}(y)\times F_{n,0} \cap \ov{\sW^c(\un{n})}\ge d\}
\]
and $C^1_{n,j,d}$ by
\[
C^1_{n,j,d}=\{y\in Y\ | \dim_k f^{-1}(y)\times F_{n,j}^1 \cap \ov{\sW^c(\un{n})}\ge d\}
\]
Write $C_{n,d}\setminus C_{n,d-1}$ and $C^0_{n,d} \setminus C^0_{n,d-1}$ as finite unions of locally closed subsets of $X$:
\begin{align*}
&C_{n,d}\setminus C_{n,d-1}=\cup_iC_{n,d,i}\\
&C^0_{n,d}\setminus C^0_{n,d-1}=\cup_iC^0_{n,d,i}\\
&C^1_{n,j,d}\setminus C^1_{n,j,d-1}=\cup_iC^1_{n,j,d,i}
\end{align*}

Let 
\begin{multline*}
\sC=\{C_{n,d,i}, C^0_{n,d,i}, C^1_{n,j,d,i}\ |\ n=1,\ldots, N, d=0,\ldots, n+s,\\ j=1,\ldots,n, i=1,2,\ldots\}.
\end{multline*}

For subschemes $Z\subset Y\times\square^a$, $W\subset X\times \A^1\times\square^b$, define $f^*(Z)\cap_XW$ as the image of $f_a^*(Z)\times_XW$ under the exchange of factors isomorphism
\[
(X\times\square^a)\times_X(X\times \A^1\times\square^b)\xrightarrow{\tau'}X\times B_{a+b}.
\]
Note that $f_a^*(Z)$ is well-defined since $f$ is smooth. Define $\ov{f^*(Z)}\cap_X\ov{W}\subset X\times\ov{B}_{a+b}$ similarly. Showing that $f^*(Z)\boxtimes W$ is in $\TZ_{\dim_kX+s-r}(X\times_kX,*;m)_\Delta$ is equivalent to showing that $f^*(Z)\cap_XW$ is in
 $\TZ_{s-r}(X,*;m)$, which we proceed to show for $Z\in \un{z}^r(Y,*)_{\ov\sC}$,
 $W\in \un\TZ_s^\sW(X,*;m)$.

For $Z, W$ integral cycles in  $\un{z}^r(Y,*)_{\sC}$, $\un\TZ_s^\sW(X,b;m)$ respectively, one easily computes:
\[
\dim_kf^{-1}_a(Z)\cap_XW\cap X\times F\le s-r+\dim_kF
\]
for each face $F$ of $B_{a+b}$. If moreover $Z$ is in $\un{z}^r(Y,a)_{\ov\sC}\subset \un{z}^r(Y,a)_{\sC}$, we have
\begin{align*}
&\dim_k\ov{f^{-1}_a(Z)}\cap_X(\ov{W}\cap X\times F_{b,0})\le s-r+a+b-1\\
&\dim_k\ov{f^{-1}(Z)}\cap_X(\ov{W}\cap X\times F^1_{b,j})\le s-r+a+b-1;\ j=1,\ldots, b.
\end{align*}
The first condition shows that each component of $f^{-1}(Z)\cap_XW$ has the correct dimension and is in good position; in particular, the cycle $f^*(Z)\cap_XW$ is defined. The second shows that 
$\ov{f^*(Z)}\times \ov{B}_b$ intersects $\tau'_*(\ov{\square}^a\times(\ov{W}\cap X\times F_{b,0}))$
and and the
 $\tau'_*(\ov{\square}^a\times(\ov{W}\cap X\times F^1_{b,j}))$ properly on $X\times\ov{B_{a+b}}$. This in turn shows that the modulus condition for $W$ implies the modulus condition for $f^*(Z) \cap_XW$.

The condition (0) in case $b=0$ trivially implies the induction condtion for $f^*(Z)\cap_XW$; in general, the inductive condition follows by induction on $b$.
\end{proof}

The same proof, ignoring the modulus and inductive conditions, yields 
\begin{lem}\label{lem:ProductML1} Fix an integer $r\ge0$. Let $\sW$ be a finite set of closed subsets $W_n\subset X\times B_n$, $n=1,\ldots,N$, such that each $W_n$ is the support of a cycle in $\un{z}^s(X,n)$. Let $f:X\to Y$ be a smooth morphism. Then there is a finite set of locally closed subsets $\sC$   such that the partially defined map $(A,B)\mapsto f^*(A)\cup_X B$ restricts to a well-defined map of complexes
\[
f^*(-)\cup_X(-):z^r(Y,*)_{\sC}\otimes \un{z}^s_\sW(X,*)\to
\un{z}^{r+s}(X,*).
\]
\end{lem}

\begin{lem}\label{lem:ProductML2} Let $s$,  $\sW$  and $f:X\to Y$ be as in lemma~\ref{lem:ProductML}, and let $\sC$ be a finite set of locally closed subsets of $X$ given by lemma~\ref{lem:ProductML}. Let $\sC_1$ be a  finite set of locally closed subsets of $X$ containing $\sC$. 

Fix an integer $r\ge0$, and let $\sT$ be a  finite set of closed subsets of $Y\times B_n$ of the form $\supp(T_n)$, where $T_n$ is in $\un{z}^r(Y,n)_{\ov{\sC_1}}$. Let $g:Y\to Z$ be a smooth morphism in $\Sm/k$. Then there is a finite set $\sC_2$ of locally closed subsets of $Z$ such that for each $W\in  \TZ_s^\sW(X,*;m)$,
$Z\in \un{z}^r_\sT(Y,*)$ and $V\in \un{z}^q(Z,*)_{\ov{\sC}_2}$:\\
\\
(1) The cycles $g^*(V)\cup_YZ$, $(gf)^*(V)\cap_X(f^*(Z)\cap_XW)$ and $f^*(g^*(V)\cup_YZ)\cap_XW$ are defined.\\
\\
(2)  $g^*(V)\cup_YZ$ is in $\un{z}^{q+r}(Y,*)$\\
\\
(3) $(gf)^*(V)\cap_X(f^*(Z)\cap_XW)$ and $f^*(g^*(V)\cup_YZ)\cap_XW$ are in 
 $\TZ_*(X,*;m)$ and
 \[
(gf)^*(V)\cap_X(f^*(Z)\cap_XW)=f^*(g^*(V)\cup_YZ)\cap_XW.
\]
 \end{lem}
 
 \begin{proof} Since the groups $\TZ_s^\sW(X,*;m)$ and  $\un{z}^r_\sT(Y,*)$ is finitely generated, it suffices to exhibit a finite set $\sC_2$ satisfying (1)-(3) for each integral $W\in\TZ_s^\sW(X,*;m)$ and integral $Z\in \un{z}^r_\sT(Y,*)$.
 
 Given $W\in\TZ_s^\sW(X,*;m)$ and integral $Z\in \un{z}^r_\sT(Y,*)$, the cycle 
 $f^*(Z)\cap_XW$ is in $\TZ_{s-r}(X,*;m)$. Letting $\sU$ be the set of intersections of 
 $\text{supp}(f^*(Z)\cap_XW)$ with faces $X\times F$, we may apply lemma~\ref{lem:ProductML}
to $\TZ_{s-r}^\sU(X,*;m)$ and morphism $gf:X\to Z$, to yield the set $\sC_2(fg)$, which is ``good" for the intersection $(gf)^*(V)\cap_X(f^*(Z)\cap_XW)$. Similarly, we may apply 
 lemma~\ref{lem:ProductML2} to the set of intersections $\sU'$ of $\text{supp}(Z)$ with faces $Y\times F$ to yield the set $\sC_2(g)$, which is ``good" for the intersection $g^*(V)\cup_YZ$.
 
Let $\sC_2= \sC_2(g)\cup \sC_2(fg)$. We claim that all $V\in  \un{z}^q(Z,*)_{\ov{\sC}_2}$ satisfy (1)-(3). For this, it suffices to show that for integral $V\in  \un{z}^q(Z,*)_{\ov{\sC}_2}$, the cycle 
$f^*(g^*(V)\cup_YZ)\cap_XW$ is well-defined (only as a cycle on $X\times B_n$ for appropriate $n$) and that this cycle is equal to  $(gf)^*(V)\cap_X(f^*(Z)\cap_XW)$. But since $f$ is smooth and the cycle  
$g^*(V)\cup_YZ$ is defined, so is the cycle
\[
f^*(g^*(V)\cup_YZ)=(gf)^*(V)\cup_Xf^*(Z).
\]
Since the cycle $(gf)^*(V)\cap_X(f^*(Z)\cap_XW)$ is also defined and all cycles involved are effective, this implies that $((gf)^*(V)\cup_Xf^*(Z))\cap f^*(Z)$ is defined. The associativity of cycle intersection implies that
 \[
(gf)^*(V)\cap_X(f^*(Z)\cap_XW)= ((gf)^*(V)\cup_Xf^*(Z))\cap f^*(Z),
\]
completing the proof.
\end{proof}

\begin{thm}\label{thm:product} Let $X$ be a smooth quasi-projective variety over $k$.\\
\\
(1) If $X$ is projective, there is a product
\[
\cap_X: \CH^r(X, p) \otimes \TH_s(X,q ; m)\to \TH_{s-r}(X, p+q ; m),
\] 
natural with respect to flat pull-back, and satisfying the projection formula 
\[
f_*( f^*(a)\cap_Xb)=a\cap_Y f_*(b)
\]
for $f:X\to Y$ a morphism of smooth projective varieties over $k$. If $f$ is a flat morphism, we have in addition the projection formula
\[
f_*(a\cap_X f^*(b))=f_*(a)\cap_Y b
\]
(2) In general, there is a product
\[
\cap_X:\CH^r(X) \otimes\TH_s(X,q ; m)  \to \TH_{s-r}(X, q ; m),
\] 
natural with respect to flat pull-back, and satisfying the projection formula 
\[
f_*( f^*(a)\cap_Xb)=a\cap_Yf_*(b)
\]
for $f:X\to Y$ a projective morphism of smooth quasi-projective varieties over $k$, and the projection formula
\[
f_*(a\cap_Xf^*(b))= f_*(a)\cap_Yb
\]
if $f$ is flat.\\
\\
In addition, all products are associative.
\end{thm}

\begin{proof} The existence of the product in (1) follows from lemma~\ref{lem:ProductML}:  the map
\[
\colim_{\sW} \TZ_s^\sW(X,*;m)\to  \TZ_s(X,*;m)
\]
is an isomorphism and by 
lemma~\ref{lem:ML3}(1), the inclusion $z^r(X,*)_{\ov\sC}\to z^r(X,*)$ is a quasi-isomorphism.

To prove the first projection formula, let $\Gamma\subset X\times_k Y$ be the 
graph of $f$. Let $\sW_X$ be a finite set  of closed subsets $W_n\subset 
X\times B_n$ as in the statement of lemma~\ref{lem:ProductML} and let $\sW$ be the finite set of subsets $(\id_X,f)\times\id_{B_n}(W_n)$ of $X\times_kY\times B_n$. We apply lemma~\ref{lem:ProductML} to the smooth morphism $p_2:X\times_kY\to Y$, the set $\sW$ and the integer $s$. This gives us the finite set of locally closed subsets $\sC$ of $Y$ and a well-defined map of complexes
\[
(-)\cap_{X\times_kY}p_2^*(-):z^r(Y,*)_{\ov\sC}\otimes \TZ_s^\sW(X\times_kY,*;m)\to
 \TZ_{s-r}(X\times_kY,*;m)
\]
Applying lemma~\ref{lem:ProductML1}  with $s=\dim_kY$ and $\sW=\{\Gamma\}$, we may add elements to $\sC$ so that we have a well-defined map of complexes
\[
p_2^*(-)\cap_{X\times_kY}(-): z^r(Y,*)_{\ov\sC}\otimes z^{\dim_kY}_{\Gamma}(X\times_kY,*)\to
 z^{s+r}(X\times_kY,*)
\]
For $Z\in z^r(Y,*)_{\ov\sC}$, the cycle $f^*(Z)$ is just $p_{1*}(\Gamma\cap p_2^*(Z))$, it follows that $f^*$ gives a well-defined map of complexes
\[
f^*:\un{z}^r(Y,*)_{\ov\sC}\to \un{z}^r(X,*).
\]
Similarly, for $W\in \un\TZ_s(X,*;m)$, 
\[
f^*(Z)\cap_XW=p_{1*}(p_2^*(Z)\cap_{X\times_kY}(\id,f)_*(W)), 
\]
so we have a well-defined map of complexes
\[
f^*(-)\cap_X(-):z^r(Y,*)_{\ov\sC}\otimes \un\TZ_s^{\sW_X}(X,*;m) \to
\un\TZ_{s-r}(X,*;m).
\]
The proof of the second projection formula for flat $f$ is similar, and is left to the reader.

Finally, if $Z$ is an integral cycle such that $f^*(Z)$ is defined and $W$ is an integral cycle such that $f^*(Z)\cap_XW$ is defined, then $Z\cap_Yf_*(W)$ is defined and
\[
f_*(f^*(Z)\cap_XW)=Z\cap_Yf_*(W).
\]
Thus sending $(Z,W)$ to $Z\cap_Yf_*(W)$ gives a well-defined map of complexes
\[
(-)\cap_Yf_*(-):z^r(Y,*)_{\ov\sC}\otimes \un\TZ_s^{\sW_X}(X,*;m) \to
\un\TZ_{s-r}(Y,*;m).
\]
and $f_*\circ(f^*(-)\cap_X(-))=(-)\cap_Yf_*(-)$. The projection formula thus follows from lemma~\ref{lem:ML3}(1).

The proof of the associativity $(z_1\cup_Xz_2)\cap_Xw=z_1\cap_X(z_2\cap_Xw)$ is proven similarly: Applying lemma~\ref{lem:ProductML2} to the projections $X^3\xrightarrow{p_{23}}X^2\xrightarrow{p_2}X$ and then using  lemma~\ref{lem:ML3}(1), we see that given classes $z_1\in \CH^r(X,a)$, $z_2\in \CH^t(X,b)$ and $w\in  \TH_s(X,q ; m)$, we can find representatives $Z_1, Z_2, W$ such that all intersections $Z_1\cup_XZ_2$, $Z_2\cap_XW$,
$Z_1\cap_X(Z_2\cap_XW)$ and $(Z_1\cup_XZ_2)\cap_XW$ are defined, are all in the appropriate complexes and satisfying
\[
(Z_1\cup_XZ_2)\cap_XW=Z_1\cap_X(Z_2\cap_XW)
\]
as cycles. 

The proof of (2) is the same, we just  use lemma~\ref{lem:ML3}(2) instead of lemma~\ref{lem:ML3}(1).
\end{proof}
 
\section{The additive Chow groups of an $S$-motive}\label{sec:AdditiveMotive}
\subsection{The action of correspondences}
Theorem~\ref{thm:product} gives us an action of correspondences on the groups $\TZ_s(X,*;m)$.
\begin{defn}
Take $X,Y\in\Sm/S$ with $Y$ projective and take $\alpha\in\Cor^r(X \times_S Y)=\CH_{r+\dim_SX}(X\times Y)$.
Define
\[
\alpha_*:\TZ_s(X,*;m)\to \TZ_{s+r}(X,*;m)
\]
to be the composition
\begin{multline*}
\TZ_s(X,*;m)\xrightarrow{p_1^*}\TZ_{s+\dim Y}(X\times_SY,*;m)\\
\xrightarrow{\alpha\cap-}
\TZ_{s+r}(X\times_SY,*;m)\xrightarrow{p_{1*}}
 \TZ_{s+r}(X,*;m)
 \end{multline*}
 \end{defn}
 
 \begin{prop}\label{prop:CorrFunct} Let $S$ be in $\Sm/k$. For $X,Y,Z\in\Sm/S$ with $Y$ and $Z$ projective, and for $\alpha\in\Cor^r(X, Y)$, $\beta\in \Cor^{r'}(Y, Z)$, we have
 \[
 (\beta\circ \alpha)_*=\beta_*\circ\alpha_*
 \]
 as maps from $\TZ_s(Z,*;m)$ to $\TZ_{s+r+r'}(X,*;m)$.
 \end{prop}
 
 \begin{proof} Relying on our previous results, especially theorem~\ref{thm:product}, the proof is standard,  using repeated application of the functoriality of flat pull-back, projective pushforward, associativity, compatiblity of projective push-forward and flat pull-back in transverse cartesian squares and the projection formula for a smooth projective morphism:
\begin{align*} 
(\beta\alpha)_*(w)&:=p_{X*}^{XZ}\left(p_{XZ*}^{XYZ}\left(p^{XYZ*}_{YZ}(\beta)\cup_{XYZ} p^{XYZ*}_{XY}(\alpha)\right)\cap_{XZ}p^{XZ*}_X(w) \right)\\
&=p^{XZ}_{Z*}\left(p^{XYZ}_{XZ*}\left[p^{XYZ*}_{YZ}(\beta)\cap_{XYZ}\left(p^{XYZ*}_{XY}(\alpha)\cap_{XYZ}p^{XYZ*}_{X}(w)\right)\right]\right)\\
&=p^{XYZ}_{Z*}\left(p^{XYZ*}_{YZ}(\beta)\cap_{XYZ}\left(p^{XYZ*}_{XY}(\alpha)\cap_{XYZ}p^{XYZ*}_X(w)\right)\right)\\
&=p^{YZ}_{Z*}\left(p^{XYZ}_{YZ*}\left[p^{XYZ*}_{YZ}(\beta) \cap_{XYZ}(p^{XYZ*}_{XY}(\alpha)\cap_{XYZ}p^{XYZ*}_X(w))\right]\right)\\
&=p^{YZ}_{Z*}\left(\beta\cap_{YZ}\left[p^{XYZ}_{YZ*}\left(p^{XYZ*}_{XY}\left(\alpha\cap_{XY}p^{XY*}_X(w)\right)\right)\right]\right)\\
&=p^{XY}_{X*}\left(\beta\cap_{YZ}\left[p^{YZ*}_Y\left(p^{XY}_{Y*}\left(\alpha\cap_{XY}p^{XY*}_X(w)\right)\right)\right]\right)\\
&=\beta_*(\alpha_*(w)).
\end{align*}
 \end{proof}
 
 The proposition immediately implies
 
 \begin{thm}\label{thm:addmotive} Let $\Gr\Ab$ denote the category of graded abelian groups. For each integer $p\ge1$ assignment
 \[
 (X,n)\mapsto  \TH_n(X,p;m) 
 \]
 extends to an additive functor
 \[
 \TH_*(p;m):\Mot(S)\to\Gr\Ab
 \]
 with
 \[
 \TH_*(p;m)(m(X)(n),\alpha)= \alpha_*(\TH_n(X,p;m))\subset \TH_n(X,p;m)).
 \]
 \end{thm}
 
 \begin{cor} For each $f:X\to Y$ in $\SmProj/S$, there is a well-defined pull-back map
 \[
 f^*:\TH^s(Y,*;m)\to \TH^s(X,*;m)
 \]
 with $(gf)^*=f^*g^*$. If $f$ is flat, then $f^*$ is the same as the flat pull-back. Finally,  the projection formula
 \[
 f_*(a\cap_Xf^*(b))=f_*(a)\cap_Yb
 \]
 is satisfied for $a\in \CH^r(X)$, $b\in \TH_s(Y,*;m)$.
 \end{cor}
 
 \begin{proof} Let $\Gamma_f\subset X\times_SY$ be the graph of $f$, and let ${}^t\Gamma_f\subset Y\times_SX$ be the transpose, giving the element $[{}^t\Gamma_f]\in \Cor_S^{\dim_SX-\dim_SY}(Y,X)$. Define $f^*=[{}^t\Gamma_f]_*$. The functoriality follows from the identity
 \[
 [{}^t\Gamma_f]\circ [{}^t\Gamma_g]=[{}^t\Gamma_{gf}]
\]
in $\Cor_S$ and proposition~\ref{prop:CorrFunct}.
 
 The fact that the new definition of $f^*$ agrees with the old one for flat $f$ follows from the identity
 \[
 (\id_X,f)_*(f^*_{old}(w))=[\Gamma_f]\cap p_2^*(w).
 \]
 
 The operations $a\cap_X(-)$, and $f_*(a)\cap_Y(-)$ can be also written as the action of correspondences, namely $\delta_{X*}(a)_*$ and $\delta_{Y*}(f_*(a))_*$, where $\delta_X:X\to X\times_SX$ and $\delta_Y:Y\to Y\times_SY$ are the diagonals.  Similarly, $f_*=[\Gamma_f]_*$.
 
The projection formula follows from proposition~\ref{prop:CorrFunct} and the identity of correspondences
 \[
 [\Gamma_f]\circ \delta_{X*}(a)\circ {}^t[\Gamma_f]=
 \delta_{Y*}(f_*(a)).
 \]
 \end{proof}
 
 \begin{rem}\label{rem:pushpull} Let $X\to S$ be a smooth projective morphism in $\Sm/k$ and let $i:Z\to X$ be a closed immersion in $\Sm/k$. Then we have the identity of operators
 on $\TH^*(X,*;m)$:
 \[
 i_*\circ i^*= [Z]\cap_X(-),
 \]
 where $[Z]\in\CH^*(X)$ denotes the class of $Z$. This follows from the projection formula:
 \begin{align*}
 i_*(i^*(w))&=i_*(1_Z\cap_Zi^*(w))\\
 &=i_*(1_Z)\cap_X w\\
 &=[Z]\cap_X w.
 \end{align*}
 \end{rem}

\subsection{Projective bundle formula}
In this section, we compute the additive higher Chow groups of the projective
bundles over smooth varieties. 

\begin{thm}\label{PBF}
Let $X$ be a smooth quasi-projective variety and 
let $E$ be a vector bundle on $X$ of rank $r+1$. Let 
$p:  {\P}(E) \to X$ be the associated projective bundle over 
$X$. Let ${\eta} \in CH^1(\P(E))$ be the class of the tautological 
line bundle $\sO(1)$.
Then for any $q , n \ge 1$ and $m \ge 2$, the map 
\[
\theta: \bigoplus_{i= 0}^r \TH^{q-i}(X, n; m)\to \TH^q(\P(E), n; m)
\]
 given by
\[
(a_o, \cdots , a_r) \mapsto {\sum}_{0 \le i \le r} 
\eta^i \cap_{\P(E)} p^*(a_i)
\]
 is an isomorphism. 
\end{thm} 
\begin{proof} Write $\P:=\P(E)$. As is well-known, the correspondences
\[
\alpha_i:=(\id_\P,p)_*(\eta^i)\in \Cor_S^i(\P,X)
\]
give an isomorphism
\[
\sum_{i=0}^r\alpha_i:m(\P) \cong \oplus_{i=0}^rm(X)(i)
\]
in $\Mot(S)$.

It is easy to see that
\[
\TH^*(*;m)(\alpha_i)(a)=\eta^i \cap_\P p^*(a),
\]
whence the result.
\end{proof}

\subsection{Blow-up formula}

 Let $i:Z\to X$ be a closed immersion in $\SmProj/S$. Denote the blow-up of $X$ along $Z$ by $\mu:X_Z
\to X$. Let $E:=\mu^{-1}(Z)$ be the exceptional divisor, giving us the commutative diagram
\[
\xymatrix{
E\ar[r]^{i_E}\ar[d]_q&X_Z\ar[d]^\mu\\
Z\ar[r]_i&X
}
\]
$E\to Z$ is the projective space bundle $p:\P(N_{Z/X})\to Z$, where $N_{Z/X}$ is the normal bundle. Applying the functor $m:\Sm\Proj/S\to \Mot(S)$ to this square gives us the sequence in $\Mot(S)$
\begin{equation}\label{eqn:Blowup}
m(E)\xrightarrow{(m(i_E),-m(q))}m(X_Z)\oplus m(Z)\xrightarrow{m(\mu)+m(i)}m(X)
\end{equation}

\begin{lem}\label{lem:Blowup} The maps $({}^t[\Gamma_\mu],0):m(X)\to m(X_Z)\oplus m(Z)$ and
$(m(i_E),-m(q)):m(E)\to m(X_Z)\oplus m(Z)$ gives an isomorphism
\[
\phi:m(E)\oplus m(X)\to m(X_Z)\oplus m(Z)
\]
exhibiting the sequence \eqref{eqn:Blowup} as a split exact sequence in $\Mot(S)$.
\end{lem}

\begin{proof} One computes that $m(\mu)\circ ({}^t[\Gamma_\mu])=\id_{m(X)}$, hence $({}^t[\Gamma_\mu],0)$ gives a splitting to $m(\mu)+m(i)$. Thus, it suffices to check that $\phi$ is an isomorphism.

By Manin's identity principle, it suffices to check that the map
\[
\id_T\otimes\phi:m(T)\otimes [m(E)\oplus m(X)]\to m(T)\otimes [m(X_Z)\oplus m(Z)]
\]
induces an isomophism
\[
(\id_T\otimes\phi)^*:\CH^*(T\times X_Z)\oplus \CH^*(T\times Z)\to
\CH^*(T\times E)\oplus \CH^*(T\times X)
\]
for all smooth projective $T\to S$. This follows easily from the known blow-up formula for the Chow groups.
\end{proof}

\begin{thm}\label{thm:Blowup} Let $i:Z\to X$ be a closed  immersion in $\SmProj/S$ and let $\mu:X_Z\to X$ be the blow-up of $X$ along $Z$, $i_E:E\to X_Z$ the exceptional divisor with morphism $q:E\to Z$. Then the sequences
\begin{multline*}
0\to \TH_s(E,n;m)\xrightarrow{(q_*,-i_{E*})} \TH_s(Z,n;m)\oplus\TH_s(X_Z,n;m)\\\xrightarrow{i_*+\mu_*}
 \TH_s(X,n;m)\to 0
 \end{multline*}
and
\begin{multline*}
0\to \TH^s(X,n;m)\xrightarrow{(i^*,\mu^*)}  \TH^s(Z,n;m)\oplus\TH^s(X_Z,;m)\\\xrightarrow{q^*-i_E^*}
\TH^s(E,n;m)\to0
\end{multline*}
are split exact.
\end{thm}
 
 \begin{proof} The result for the second sequence follows by applying the functor $\TH(n;m)$ to the sequence \eqref{eqn:Blowup} and using lemma~\ref{lem:Blowup}. For the first sequence, we take the dual of the sequence  \eqref{eqn:Blowup} and use the same argument.
\end{proof}

\subsection{Additive Chow groups of varieties with group actions}
Let $k$ be a perfect field which may not be algebraically closed.
As an application of our extension of the additive chow groups to a functor on 
the category of motives Mot$(S)$ (theorem~\ref{thm:addmotive}), we 
compute the additive Chow groups of smooth varieties with group actions, in
particular,  for projective homogeneous spaces
and Grassmann bundles over smooth projective varieties, generalising
theorem~\ref{PBF}.

\begin{thm}\label{thm:torusaction}
Let $X$ be a smooth projective variety of dimension $d$ over $k$ equipped with an action of
the multiplicative group ${\mathbb G}_m$. Let $\{Z_i , 0 \le i \le r \}$ be
the connected components of the fixed point locus. Then $Z_i$'s are all
smooth projective and one has for $n \ge 0, p \ge 1$ and $m \ge 1$ the
formula
\begin{equation}\label{torus}
\TH_n(X,p ; m) \cong \oplus_{i=0}^{r} \TH_{a_i}(Z_i,p;m)
\end{equation}
where $a_i$ is the dimension of the positive eigenspace of the action of
${\mathbb G}_m$ on the tangent space of $X$ at an arbitrary point $ z \in
Z_i$.
\end{thm}  

\begin{proof} By a theorem of Bialynicki-Birula \cite{Biyalinicki},
generalized by Hesselink \cite{Hesselink} to case of non-algebraically
closed fields, the fixed point locus $X^{{\mathbb G}_m}$ is smooth, closed
subscheme of $X$ with the decomposition into connected components as
stated. Furthermore, $X$ has a filtration
$$X = X_r \supset X_{r-1} \supset \cdots \supset X_0 \supset X_{-1} = \0$$
and $X_i\setminus X_{i-1}$ admits the structure of an  $\A^{a_i}$-bundle ${\phi}_i : X_i \setminus X_{i-1} \rightarrow Z_i$.

By a result of Karpenko and Chernousov-Gille-Merkurjev \cite{CGM}, theorem~7.1, the motive of $X$ has the following decomposition in the
category $\Mot(k)$: $m(X) = \oplus_{i=0}^{r} m(Z_i)(a_i)$, hence
\begin{equation}\label{tate}
m(X)(n) = \oplus_{i=0}^{r} m(Z_i)(n+a_i)
\end{equation}
The theorem now follows by applying theorem~\ref{thm:addmotive} (with
$S = \Spec(k)$) to the above decomposition. 
\end{proof}

\begin{remk}\label{PHS}
Theorem~\ref{thm:torusaction} applies to a smooth projective variety $X$ which has
an action of a $k$-split reductive group $G$, by restricting the action to a given $\G_m$ in $G$. Furthermore, if $X$ is
a quasi-homogeneous space of such a reductive group, then Brosnan 
\cite{Brosnan} has shown that the components $Z_i$'s are also 
quasi-homogeneous spaces of much smaller dimensions and hence their
additive Chow groups should be easier to compute.
\end{remk}

\begin{cor}\label{PHS1}
Let $G$ be a $k$-split reductive group and let $P$ be a parabolic subgroup
of $G$. Then one has for $n \ge 0, p \ge 1$ and $m \ge 1$, 
$$\TH_n(G/P, p;m) \cong \oplus_{w} \TH_{n-l(w)}(k,p;m)$$
where $w$ runs through the set of cells of $G/P$ in its Bruhat decomposition
and $l(w)$ is the dimension of the corresponding cell.
\end{cor}

\begin{proof}
 K\"ock 
has shown in \cite{Kock} that in this case one can choose  a $\G_m\subset G$ so that the varieties $Z_i$'s are
just rational points. The result thus
 follows directly from theorem~\ref{thm:torusaction}. \end{proof}

\begin{remk}\label{PHS2}
Corollary~\ref{PHS1} shows in particular that the additive Chow groups of
a Grassmann variety has a decomposition in terms of the additive Chow groups
of the ground field. One can write down the similar formula for the
additive Chow groups of any Grassmann bundle over a smooth projective
variety $X$ in terms of the additive Chow groups of $X$ because one
has the corresponding decomposition in Mot$(k)$ (see \cite{Kock}).
\end{remk}

\section{Logarithmic additive Chow groups}\label{sec:log}
The additive Chow groups quite clearly do not satisfy the homotopy invariance property; for this reason it is difficult to apply the existing technology to prove either a localization property with respect to a closed immersion, or a Mayer-Vietoris property for a Zariski open cover. However, since we have a blow-up exact sequence, we can use the machinery  of Guill\'en and Navarro Aznar to define ``additive Chow groups with log poles", at least if we assume $k$ admits resolution of singularities.

\subsection{The extension theorem}

\begin{thm}\label{thm:AdditiveChowExt} The functor
\[
\TZ_r(-;m):\SmProj/k\to C^-(\Ab)
\]
extends canonically to a functor
\[
\TZ_r(-,m):\sD_\hom(k)\to D^-(\Ab).
\]
\end{thm}

\begin{proof} Since $\TZ_r(-;m)$ is additive, there is a canonical extension to an additive functor
\[
\TZ_r(-;m):\Z\SmProj/k\to C^-(\Ab)
\]
Taking the total complex and passing to the homotopy categories gives the exact functor
\[
\TZ_r(-;m):K^b(\Z\SmProj/k)\to K^-(\Ab)
\]
For each blow-up $\mu:X_Z\to X$ in $\SmProj/k$, theorem~\ref{thm:Blowup} shows that
$\TZ_r(C(\mu);m)$ is acyclic. This gives us the desired extension
\[
\TH_r(-,m):\sD_\hom(k)\to D^-(\Ab).
\]
\end{proof}

\begin{cor}\label{cor:Ext} The functor $\TZ_r(-;m):\SmProj/k\to C^-(\Ab)$ extends to a functor
\[
\TZ_r^\log(-;m):\Sch_k'\to D^-(\Ab)
\]
satisfying:\\
\\
1. Let $\mu:Y\to X$ be a proper morphism in $\Sch_k$, $i:Z\to X$ a closed immersion. Suppose that $\mu:\mu^{-1}(X\setminus Z)\to X\setminus Z$ is an isomorphism. Set $E:=\mu^{-1}(X\setminus Z)$ with maps $i_E:E\to Y$, $q:E\to Z$. There is a canonical extension of
the sequence in $D^-(\Ab)$
\[
\TZ_r^\log(E;m)\xrightarrow{(i_{E*},-q_*)}\TZ^\log_r(Y;m)\oplus \TZ^\log_r(Z;m)\xrightarrow{\mu_*+i_*}
\TZ^\log_r(X;m)
\]
to a distinguished triangle in $D^-(\Ab)$.\\
\\
2. Let $i:Z\to X$ be a closed immersion in $\Sch_k$, $j:U\to X$ the open complement. Then there is a canonical distinguished triangle in $D^-(\Ab)$:
\[
\TZ^\log_r(Z;m)\xrightarrow{i_*}\TZ^\log_r(X;m)\xrightarrow{j^*}\TZ^\log_r(U;m)\to \TZ^\log_r(Z;m)[1],
\]
which is natural with respect to proper morphisms of pairs $(X,U)\to (X',U')$.
\end{cor}

\begin{proof} This follows from theorems~\ref{thm:Extension}  and \ref{thm:AdditiveChowExt}.
\end{proof}

\subsection{The logarithmic additive Chow groups of an extented motive} \label{subsec:LogChow}

For the Borel-Moore motive $bm(X)(r)$ of a finite-type $k$-scheme $X$, we have the associated object $\TZ^\log_r(X;m)$ in $D^-(\Ab)$; sending $(X,r)$ to $\TZ^\log_r(X;m)$ in $D^-(\Ab)$ defines a functor
\[
\TZ^\log(-;m):\Sch_k'\times\Z\to D^-(\Ab).
\]
We have as well the functor
\[
bm:\Sch_k'\times\Z\to \hat\Mot(k)\subset\DMH(k)
\]
which is $bm(r)$ on $\Sch_k'\times r$. 

For each $n\in\Z$, let $\TH^\log(-,n;m)$ be the functor $H_n(\TZ^\log(-;m))$.

\begin{thm}\label{thm:LogChowMot} For each $n\in\Z$, the functor $\TH^\log(-,n;m)$ extends via $bm$ to
\[
\TH_*^\log(-,n;m):\hat\Mot(k)^*\to \Ab
\]
such that the restriction of $\TH^\log(-,n;m)$ to $\Mot(k)\subset \hat\Mot(k)^*$ is the functor $\TH(-,n;m)$
\end{thm}

\begin{proof} Since $\hat\Mot(k)$ is the full pseudo-abelian subcategory of $\DMH(h)$ generated
by the objects $bm(X)(r)$ with $X\in\Sch_k$, it suffices to define an action  of higher correspondences
$\alpha\in\Hom(bm(X)(r), bm(Y)(s)[a])$,
\[
\alpha_*:\TH^\log_r(X,n;m)\to \TH^\log_s(Y,n+a;m)
\]
agreeing with the action  already defined for $X,Y\in\SmProj/k$, $a=0$, and satisfying the functoriality
\[
(\beta\circ\alpha)_*=\beta_*\circ\alpha_*.
\]
The construction is essentially the same as for $\Mot(k)$. Take $X\in\Sch_k$. Then there is a cubical object $\tilde X_*$ of $\SmProj/k$ such that the image under $C^b(m(r))$ of the total complex of 
$\tilde X_*$ in $C^b(\Z\SmProj/k)$ is a representative for $bm(X)(r)\in\DMH(k)$. In other words, the finite diagram $(K^m, f^{m+1,m})$ with 
\[
K^m=\oplus_{|I|=1-m}(\tilde X_I, r)
\]
and $f^{m+1,m}:K^m\to K^{m+1}$ the map induced by the graphs of the morphisms $f_{I,J}:\tilde X_I\to
\tilde X_J$ with $|I|=|J|+1$ represents $bm(X)(r)$. Note that 
\[
f_{I,J*}:\TZ_r(\tilde X_I,*;m)\to \TZ_r(\tilde X_J,*;m)
\]
is well-defined for every $I,J$,  yielding  the double complex $\TZ(K^*,*;m)$ with second differential given by the maps $f_{I,J*}$. The total complex $\Tot(\TZ(K^*,*;m))$ is thus by definition a representative in $C^-(\Ab)$ for $\TZ^\log_r(X,*;m)$.

Taking a second scheme $Y\in\Sch_k$, we have a choice of cubical object $\tilde Y_*$, giving us the finite diagram $L=(L^m, g^{m+1,m})$ and the identification of $\TZ^\log_r(X,*;m)$ with
$\Tot(\TZ(L^*,*;m))$. 

Choose an element $[x]\in H_n(\Tot(\TZ(K^*,*;m)))$. We may represent $[x]$ by some
$x\in Z_n(Tot(\TZ(K^*,*;m)))$, i.e., a collection of cycles $x_I\in \TZ_r(\tilde X_I, n-|I|;m)$ Letting $W_I$ be the set of all irreducible components of all cycles involved in $x_I$, we have for each $I$ the cubical subset $\tilde W_I$ generated by $W_I$. We enlarge this by adding all the images of the $\tilde W_J$ for all compositions of maps $X_J\to\ldots\to X_I$ involved in the diagram $K^*$, giving the finitely generated sub double complex $\TZ^\sW(K^*,*;m)$ of $\TZ(K^*,*;m)$. 

Applying lemma~\ref{lem:ProductML2}, we may construct highly distinguished subcomplexes
$C^*(K^m,L^n)'\subset C^*(K^m, L^n)$ such that the associated internal Hom complex $\sHom(K,L)'$ is 
defined, and that the action of correspondences for the individual schemes $\tilde X_I\times \tilde Y_J$ give rise to a well-defined action
\[
(-)_*(-):\sHom(K,L)'\otimes \Tot(\TZ^\sW(K^*,*;m))\to \Tot(\TZ(L^*,*;m)).
\]
For a given $\alpha\in\Hom(K,L[n])=\Hom(bm(X)(r),bm(Y)(s)[n])$, we define $\alpha_*([x])$ by 
\[
\alpha_*([x]):=\tilde\alpha_*(x)
\]
where $\tilde\alpha\in Z^a(\sHom(K,L))$ is a representative for $\alpha$. 

The argument we used to show that the composition law in $\DMH(k)$ is well-defined and associative also shows that $\alpha_*([x])$ is well-defined, independent of the choices we have made for representing elements and highly distinguished subcomplexes, and that we have the functoriality
\[
\beta_*(\alpha_*([x]))=(\beta\alpha)_*([x]),
\]
for $\beta\in \Hom(bm(Y)(s), bm(Z)(t)[b])$, once we are given a cubical object $\tilde Z_*$ of $\SmProj/k$ for $Z$,  yielding the resulting finite diagram $(L,h^{**})$ representing $bm(Z)(t)$.

Since different choices of cubical objects $\tilde X_*$, $\tilde Y_*$ representing $X$ and $Y$ give rise to isomorphic objects in $\DMH(k)$, the functoriality we have already proved shows that the map $\alpha_*$ is also independent of these choices.

In case $X$ and $Y$ are in $\SmProj/k$, we may take the $\tilde X_*=X$, $\tilde Y_*=Y$, so that $K=(X,r)$, $L=(Y,s)$ and we are back to our original definition of the action of the correspondence $\alpha\in \Cor^{s-r}(X,Y)$.
\end{proof}

For $X\in\Sch_k$ equi-dimensional over $k$, set $\TH^s_\log(X,n;m):=\TH_{\dim X-s}^\log(X,n;m)$ and extend this notation to $X$ locally equi-di\-men\-sional by taking the direct sum over the connected components of $X$.

\begin{cor}\label{cor:AddChowFunct}  Sending $X\in\Sm/k$ to $\TH^s_\log(X,n;m)$ extends to a functor
\[
\TH^s_\log(-,n;m):\Sm/k^\op\to \Ab.
\]
In addition, if $j:U\to X$ is an open immersion in $\Sm/k$, the map
\[
\TH^s_\log(j,n;m):\TH^s_\log(X,n;m)\to
\TH^s_\log(U,n;m)
\]
agrees with the map $j^*$ from corollary~\ref{cor:Ext}.
\end{cor}

\begin{proof}  The existence of the functor  $\TH^s_\log(-,n;m)$ follows from theorems~\ref{thm:MotFunct} and \ref{thm:LogChowMot}. 

For the assertion on the open immersion $j:U\to X$, the result follows essentially from  the definitions in case the complement of $U$ in $X$ is a strict normal crossing divisor. In general, temporarily write $[j]^*$ for the map in the distinguished triangle of corollary~\ref{cor:Ext}. Let $\mu:X'\to X$ be a projective birational morphism which is an isomorphism over $U$ such that the induced open immersion $j':U\to X'$ case strict normal crossing complement.  $\mu_*$ induces the map of distinguished triangles
\[
\xymatrix{
bm(X'\setminus U)\ar[r]^{i'_*}\ar[d]_{\mu_*}&
 bm(X')\ar[r]^{j^{\prime*}}\ar[d]_{\mu_*}&bm(U)\ar[r]\ar@{=}[d]&bm(X'\setminus U)[1]\ar[d]_{\mu_*[1]}\\
 bm(X\setminus U)\ar[r]_{i_*}&
 bm(X)\ar[r]_{[j]^{*}}&bm(U)\ar[r]& bm(X\setminus U)[1]
 }
 \]
i.e $j^{\prime*}=[j]^*\circ\mu_*$.   One easily calculates that 
$\mu_*\mu^*=\id_{bm(X)}$ and $j^{\prime*}=j^*\circ\mu_*$, hence $[j]^*=j^*$.
\end{proof}

\subsection{Codimension one} We conclude this section with a computation of  
$\TH^1_\log(X,1;m)$ for $X$ smooth and quasi-projective over $k$; we thank K. R\"ulling for suggesting that we make this computation. We continue to assume that $k$ admits resolution of singularities. 

\begin{prop}\label{prop:Computation} Suppose that $p:X\to \Spec k$ is geometrically irreducible of dimension $d$ over $k$. Then the pull-back map
\[
p^*:\TH_0^\log(\Spec k,n;m)=\TH_0(\Spec k,n;m)\to \TH_d^\log(X,n;m)
\]
is an isomorphism for all $n,m\ge1$. 
\end{prop}

\begin{proof} Let $j:X\to \bar X$ be an open immersion to a smooth projective geometrically irreducible $\bar X$ over $k$, with complement $i:Z\to \bar X$. We have the localization distinguished triangle 
\[
\TZ^\log_d(Z;m)\xrightarrow{i_*}\TZ^\log_d(\bar X;m)\xrightarrow{j^*}\TZ^\log_d(X;m)\to \TZ^\log_d(Z;m)[1]
\]
from corollary~\ref{cor:Ext}. Since $\bar X$ is projective, we have
\[
\TZ^\log_d(\bar X;m)=\TZ_d(\bar X;m)
\]
and similarly for $Z$. Since $d>\dim_kZ$, the complex $\TZ_d(Z;m)$ is the 0 complex and thus  the map $j^*$ is a quasi-isomorphism. This reduces us to the case of projective $X$.

Let $W\subset X \times B_n$ be a generator in $\un{\TZ}_d(X, n; m)$. In particular, $W$ is a codimension one integral closed subscheme of $X\times \A^1\times\square^{n-1}$, with
\[
W\cap X\times0\times\square^{n-1}=\0
\]
Since $X$ is projective, the projection $\bar W$ of $W$ to $B_n$ is closed and disjoint from $0\times\square^{n-1}$. Also, $\bar W$ is irreducible; since $X$ is geometrically irreducible over $k$, $X\times\bar W$ is irreducible. By reason of codimension, we therefore have
\[
W= p_{B_n}^{-1}(\bar W),
\]
and thus
\[
p^*:\TZ_0(\Spec k,*;m)\to \TZ_d(X,*;m)
\]
is an isomorphism of complexes. Since
\[
\TZ^\log_d(X,*;m)=\TZ_d(X,*;m),\ \TZ^\log_0(\Spec k,*;m)=\TZ_0(\Spec k,*;m)
\]
the proof is complete.
\end{proof}

\begin{rem} Using localization as above, if $X$ has irreducible components $X_1,\ldots, X_r$, such that $X_1,\ldots, X_n$ all have maximal dimension $d$ and $X_{n+1}, \ldots, X_r$ have dimension $<d$, and if $k_i$ is the field of constants in $k(X_i)$, then
\[
\TH_d^\log(X,n;m)\cong\oplus_{i=1}^n \TH_0(\Spec k_i,n;m).
\]

Combining this with R\"ulling's result \cite{Ruelling} that  $\TH^n(\Spec k,n;m)\cong \W_m\Omega^{n-1}_k$ gives the isomorphism
\[
\TH_d^\log(X,1;m)\cong\oplus_{i=1}^n \W_m(k_i).
\]
In particular, if $X$ is smooth and geometrically irreducible over $k$, then 
\[
\TH^1_\log(X,1;m)\cong  \W_m(k).
\]

Park has made an analog of the Beilinson-Soul\'e vanishing conjectures for the groups 
$\TH^q(\Spec k,p;m)$; for $q=1$, these say that $\TH^1(\Spec k,p;m)=0$ for $p\neq1$. 
\end{rem}

\section{Additive chow groups of 1-cycles on fields}\label{sec:1cycle}
The additive Chow groups of 0-cycles on a field $k$ were
studied by Bloch-Esnault \cite{BlochEsnault2} (for modulus $m=1$) and
 R\"ulling  \cite{Ruelling} (for general modulus $m$). Their main results are the construction of
natural regulator maps  
\begin{equation}\label{eqn:zcycle}
R^n_{0,m} : \TH^n(k,n;m)  \rightarrow {\W}_m {\Omega}^{n-1}_k.
\end{equation}
for $n \ge 1$ and $m \ge 1$ (the case $m=1$ is due to Bloch-Esnault), which
they  show are isomorphisms. Here the groups on the
right are the homology of the generalized deRham-Witt complex of 
Hasselholt-Madsen and ${\Omega}^{n-1}_k: ={\Omega}^{n-1}_{k/{\Z}}$
are absolute differentials of $k$. 

In particular, for $n \ge 1$ and $m \ge 1$, there is a natural surjection
\begin{equation}\label{eqn:zcycle1}
\TH^n(k,n;m)  \surj {\Omega}^{n-1}_k .
\end{equation}
These results motivate one to look for existence of such nontrivial 
regulator maps on the additive Chow groups of higher dimensional cycles 
for a field. If the field $k$ has characteristic zero, Park \cite{Park}, main theorem,
 has constructed regulator maps for $n \ge 1$ and 
$m \ge 1$,
\begin{equation}\label{eqn:1cycle1}
\TZ^{n-1}(k,n;m) \to {\Omega}^{n-3}_k 
\end{equation}
which induces the natural maps 
\begin{equation}\label{eqn:1cycle}
R^n_{1,m} :\TH^{n-1}(k,n;m)  \rightarrow {\Omega}^{n-3}_k .
\end{equation}
He has shown in \cite{Park1}, theorem~1.12, that this regulator map
is nontrivial for $n=3$. 

As an application of theorem~\ref{thm:product},
we prove in this section the following strengthening of the result of Park:
\begin{thm}\label{thm:reg}
Let $k$ be an algebraically closed field of characteristic zero. Then for
$n \ge 1$ and $m \ge 1$, the regulator map 
\[
R^n_{1,m} :\TH^{n-1}(k,n;m)  \rightarrow {\Omega}^{n-3}_k
\]
is surjective. In particular, all the additive Chow groups of 1-cycles
on $k$ are nontrivial for $n \ge 3$.
\end{thm}

\begin{rem} In view of the known calculations of relative $K$-theory of
nilpotent ideals and the expectation that the additive Chow groups 
would compute these relative $K$-groups, it is conjectured that the
regulator maps $R^n_{1,2}$ should be isomorphisms. In that context, the
above result gives some evidence for the conjecture.
\end{rem}

The proof of this theorem will be given after we establish some other 
preliminary results. In this section, we will assume our field $k$ to have characteristic zero
Recall the notation from section~\ref{subsec:AdditiveCycleComplex}: for $n \ge 1$, we have 
$B_n = \A^1 \times \square^{n-1}$, 
$\ov{B}_n = \A^1 \times ({\P}^{1})^{n-1}\supset B_n$ and
$\widehat{B}_n = {\P}^{1} \times ({\P}^{1})^{n-1} \supset \ov{B}_n$. 

The group of units $k^*$ acts on $B_n$ by
$$a * (x, t_1, \cdots , t_{n-1}) = (\frac {x}{a}, t_1, \cdots , t_{n-1})$$
This action clearly extends to an action of $k^*$ on $\ov {B}_n$ and 
$\widehat{B}_n$ by the same formula, where now $x$ could take value $\infty$.
It is easy to check that $k^*$ acts on the groups
$\TZ_r(k, n; m)$ such that all the boundary maps are $k^*$-equivariant.
In particular, the additive Chow groups are equipped with $k^*$-action.

For an irreducible curve $C \subset B_n$ such that $[C] \in \TZ^{n-1}
(k, n; m)$, let $\ov C$ and $\widehat C$ be the closures of $C$ in
$\ov {B}_n$ and $\widehat{B}_n$ respectively. Let 
\[
p :\widehat{B}_n \rightarrow {\P}^1
\]
 be the projection map. For $a \in k^*$, let 
${\sigma}_a$ denote the multiplication by $a$ on ${\P}^1$. We denote
the induced action on $\widehat{B}_n$ also by ${\sigma}_a$. We let $k^*$
act on $k$ by multiplication.

\begin{lem}\label{lem:acl}
For any $a \in k^*$, one has $\widehat {{\sigma}_a(C)}=
{\sigma}_a (\widehat C)$ and ${p}^{-1}(0) \cap {\widehat C} =
 {p}^{-1}(0) \cap {\sigma}_a (\widehat C)$.
\end{lem}  
\begin{proof} Let ${\sigma}_a :\widehat{B}_n \rightarrow  \widehat{B}_n$
be the multiplication map. This map is clearly closed which gives
$\widehat {{\sigma}_a(C)} \subset {\sigma}_a (\widehat C)$. To show the
other inclusion, let $z \in {\sigma}_a (\widehat C)$ and let $y \in
{\widehat C}$ such that $z = {\sigma}_a (y)$. Let $U$ be an open
neighborhood of $z$ and put $V = {\sigma}_a^{-1}(U)$. Then 
$V \cap C \neq \emptyset$ which gives ${\sigma}_a (C) \cap V \neq \emptyset$.
This proves the first part. The second part is immediate from the description
of the action.
\end{proof}

Fix integers $n \ge 1, m \ge 1$ and let $C \subset B_n$ be an irreducible 
curve such that 
$[C] \in \TZ^{n-1}(k, n; m)$. Recall that the modulus condition for
additive cycle $[C]$ (see section~\ref{sec:AdditiveChow}) is equivalent to the
following: for each closed point $P \in {p}^{-1}(t = 0)$ on  
the normalization ${\widehat C}^N$, there is an index $i, 1 \le i \le n-1$ 
such that
\begin{equation}\label{eqn:mod}
ord_P\left(p^*(F^1_{n,i}) - (m+1) p^*(F_{n,0})\right) \ge 0
\end{equation}
In such a case, we write $(C,P) \in {\mathcal M}^m(t_i)$.

Fix $a \in k^*$ and let $C_a$ denote ${\sigma}_a(C)$. We similarly
denote $\sigma_a(\widehat{C})$ by $\widehat{C}_a$.
 Consider the diagram  
$$ 
\xymatrix{
{\widehat C}^N \ar[r]^{\pi}\ar[d]_{{\sigma}^N_a} &
{\widehat C}\ar[r] \ar[d]_{{\sigma}_a} &
{\widehat{B}_n} \ar[d]^{{\sigma}_a} \\
{\widehat C}^N_a \ar[r]^{{\pi}_a} & {\widehat C}_a \ar[r] &
{\widehat{B}_n} 
}
$$
where ${\sigma}^N_a:{\widehat C}^N\to {\widehat C}^N_a$ is the map induced by $\sigma_a$ on the normalizations. Let 
\[
p_a:{\widehat C}^N_a\to \widehat{B}_n
\]
be the composition of the bottom arrows in the above diagram. 

Let $\{P_1, \cdots , P_r \}:={p}^{-1}(t=0)$. Then by 
lemma~\ref{lem:acl}, ${p_a}^{-1}(t=0) =
\{Q_1, \cdots , Q_r \}$, where $Q_i = {{\sigma}^N_a}(P_i)$. 

\begin{lem}\label{lem:mod1}
For each $i = 1, \cdots , r$, and $j = 1, \cdots , n-1$, 
$(C,P_i) \in {\mathcal M}^m(t_j)$ if and only if 
$(C_a,Q_i) \in {\mathcal M}^m(t_j)$.
\end{lem} 
\begin{proof} Let $\iota:F^1_{n,j}\to \widehat{B}_n$ be the inclusion. We have the cartesian diagram
\[
\xymatrix{
F^1_{n,j}\ar[r]^\iota\ar[d]_{\sigma_a}&\widehat{B}_n\ar[d]^{\sigma_a}\\
F^1_{n,j}\ar[r]^\iota&\widehat{B}_n
}
\]
from which it follows that
\[
\sigma^{N}_a:p^*(F^1_{n,j})\to p_a^*(F^1_{n,j})
\]
is an isomorphism. Similarly, 
\[
\sigma^{N}_a:p^*(F_{n,0})\to p_a^*(F_{n,0})
\]
is an isomorphism, from which the lemma follows directly.
\end{proof}

We recall from \cite{Park} the   rational absolute K\"ahler
differential $(n-1)$-forms ${\omega}^n_{l,m} \in
{\Gamma} \left(\widehat {B}_{n+1}, 
{\Omega}^{n-1}_{\widehat{B}_{n+1}/\Z}({\rm log} F_{n+1})(*\{x=0\})\right)$:
\begin{align*}
&\omega^n_{1,m} = \frac {1-t_1}{x^{m+1}} \frac {dt_2}{t_2}\wedge
\cdots \wedge \frac {dt_n}{t_n}\\
&\omega^n_{l,m} = \frac {1-t_l}{x^{m+1}}\frac {dt_{l+1}}{t_{l+1}}
\wedge \cdots \wedge \frac {dt_n}{t_n} \wedge \frac {dt_1}{t_1}
\wedge \cdots \wedge\frac {dt_{l-1}}{t_{l-1}}\ \ \ 
(1 < l < n)  \\
&\omega^n_{n,m} = \frac {1-t_n}{x^{m+1}}\frac {dt_1}{t_1} \wedge
\cdots \wedge\frac {dt_{n-1}}{t_{n-1}}
\end{align*}
We shall denote these forms simply by $\omega_l$ when $n$ and $m$
are fixed. For the definition of residue of a meromorphic form at a closed
point of a nonsingular curve used here, we refer the reader to \cite{Park}.

We also recall here the definition of the regulator maps $R^n_{1,m}$; we refer the reader to \cite{Park}, section~3,
 for details. For cycle of the form $[C] \in \TZ^{n-1}
(k, n; m)$ where $C \subset B_n$ is an irreducible curve,
we have
\begin{equation}\label{eqn:regd}
R^n_{1,m}([C]) = \sum_{P \in Supp(p^*\{x=0\})} R^n_{1,m} 
\left(C,P \right)
\end{equation}
where 
\begin{equation}\label{eqn:regd1}
R^n_{1,m} \left(C,P \right) := (-1)^{l-1} {res}_P \left 
( p^*({\omega}_l) \right) \ {\rm if} \ \ \left(C,P \right) \in
{\mathcal M}^m(t_l)
\end{equation}
for  $1 \le l \le n-1$; this definition extends to  $\TZ^{n-1}(k, n; m)$ by $\Z$-linearity.  

The following result, which explains how the residues behave with respect
to the action of $k^*$ on a smooth curve, will be one of our main tools
to prove theorem~\ref{thm:reg}.

\begin{prop}\label{res}
Assume $n=3$. Then for any $i = 1, \cdots , r$, we have for 
$(C,P_i) \in {\mathcal M}^m(t_j)$, 
\[
res_{P_i} \left (p^*({\omega}_l) \right ) = {\frac {1}{a^{m+1}}}
 res_{Q_i} \left (p^*_a({\omega}_l) \right ) \ \ {\rm for} \ l = 1, 2.
 \]
\end{prop} 

\begin{proof} We give the prooof for $\omega_1$; the other cases are
identical.  We can work locally in order to compute the residue at a given point; take $\A^1:=\P^1\setminus\{\infty\}$,  with $0\in \A^1$ corresponding to $1\in\P^1$.  Let  $A = k[x, t_1, t_2]$ $X =\A^1 \times (\A^1)^2 = \Spec(A)$. Then the action of $a\in k^*$ translates to the
$k$-algebra automorphism $f: A \rightarrow A$ given by $f(x) = a^{-1}x,
f(t_1) = t_1, f(t_2) = t_2$. 

Let $P = p(P_i) 
= (0, {\alpha}_1, {\alpha}_2)$. Then   ${\sigma}_a(P) = P$.
Let ${\mathcal O}_{X,P}$ be the local ring of $X$ at $P$. We obtain a commutative diagram of complete local rings.
\[\xymatrix{
\widehat{{\mathcal O}_{X,P}}  \ar[r]^{h_a}
\ar[d]^{\hat{f}}  & \widehat{{\mathcal O}_{{\widehat C}^N_a ,Q_i}} \ar[d]^{\phi} \\
\widehat{{\mathcal O}_{X,P}} \ar[r]_{h} & \widehat {{\mathcal O}_{{\widehat C}^N ,P_i}}
}
\]
Observe that the vertical maps are isomorphisms. Then
\begin{multline}\label{eqn:res}
{res}_{Q_i} \left ({p_a}^*{\omega}_1 \right )  =  
{res}_{Q_i} \left ({h_a}({\omega}_1) \right )   
\\
 =  {res}_{P_i}\left( {\phi} \circ h_a\left ({\omega}_1\right) \right )
 =  {res}_{P_i}\left(h\circ \hat{f}\left ({\omega}_1\right) \right ). 
 \end{multline}

Choose a parameter $t:=t_{P_i}$ of the local
ring at $P_i$ and use the parameter $t_{Q_i}:=\phi^{-1}(t)$ in ${\widehat{{\mathcal O}_{{\widehat C}^N_a ,Q_i}}} $.  There are power series 
\[
x(t), t_1(t), t_2(t)\in {\widehat{{\mathcal O}_{{\widehat C}^N ,P_i}}} \cong k(P_i)[[t]]
\]
with $h(x)=x(t)$, etc. Since $f(t_i)=t_i$ and $f(x)=a^{-1}x$, we have
\[
h \circ {\hat f} (x)=a^{-1}x(t),\ h \circ {\hat f} (t_1)=t_1(t),\ h \circ {\hat f} (t_2)=t_2(t).
\]
Thus, using the explicit description of ${\omega}_1$, we have
\[
p^*({\omega}_1) = p^* \left ({\frac {1-t_1}{x^{m+1}}} {\frac {dt_2}{t_2}}
\right ) = {\frac {1-t_1(t)}{{x(t)}^{m+1}}} {\frac {d(t_2(t))}{t_2(t)}}.
\]
Also
\begin{align*}
h \circ {\hat f}\left (\omega_1\right)&=h \circ {\hat f}\left ({\frac {1-t_1}{x^{m+1}}} 
{\frac {dt_2}{t_2}} \right )\\
 & = 
{\frac {1 - h \circ {\hat f}(t_1)}{(h \circ {\hat f} (x))^{m+1}}} 
{\frac {d\left( h \circ {\hat f}(t_2)\right )}{h \circ {\hat f}(t_2)}} \\
& =  {\frac {1 - t_1(t)}{(a^{-1}x)^{m+1}}} 
{\frac {d\left( t_2 (t)\right )}{t_2(t)}} \\
&=a^{m+1}\cdot p^*(\omega_1)
\end{align*}

 It follows directly from the definition of residue  that ${res}_{P_i}$ is $k$-linear. Thus
\[
{res}_{P_i}(h \circ {\hat f}\left (\omega_1\right))=a^{m+1}{res}_{P_i}(p^*({\omega}_1)).
\]
This together with \eqref{eqn:res} completes the proof.
\end{proof}

\begin{cor}\label{cor:reg1}
Consider the regulator map 
$$R^3_{1,m} : \TH^2(k,3;m) \rightarrow k$$
Then  for each $\alpha \in
\TH^2(k,3;m)$ and $a \in K^*$ 
\[
R^3_{1,m} (a * {\alpha}) = {\frac {1} {a^{m+1}}} R^3_{1,m} ({\alpha})
\]
\end{cor}

\begin{proof}
Since the regulator map is  $\Z$-linear, it suffices to show that the above equality holds when $\alpha$
represents the class of an irreducible curve $C \subset B_3$.
But in that case we have from  \eqref{eqn:regd}
$$R^3_{1,m} (a*[C]) = {\sum}^{r}_{i=1} R^3_{1,m} 
\left ( a*C, Q_i \right ) = {\sum}^{r}_{i=1} R^3_{1,m} 
\left ( C_{{a}^{-1}}, Q_i \right )$$ 
Now for a given $i, 1 \le i \le r$, we have by lemma~\ref{lem:mod1},
$\left ( C_{{a}^{-1}}, Q_i \right ) \in {\mathcal M}^m(t_l) \Leftrightarrow 
\left ( C, P_i \right ) \in {\mathcal M}^m(t_l)$. Applying 
proposition~\ref{res} in \eqref{eqn:regd1}, we get
$$
\begin{array}{lll}
R^3_{1,m} \left ( C_{{a}^{-1}}, Q_i \right ) & = & 
(-1)^{l-1} {res}_{Q_i} \left ( p^*_{{a}^{-1}} ({\omega}_l) \right ) \\
& = & (a^{-1})^{m+1} (-1)^{l-1} {res}_{P_i} \left ( p^* ({\omega}_l)) \right )
\\
& = & a^{-m-1} R^3_{1,m} \left ( C, P_i \right )
\end{array} 
$$  
Summing over $i$ gives the result.
\end{proof}

The following is the case $n=3$ of theorem~\ref{thm:reg}.
\begin{cor}\label{reg2}
If $k$ is algebraically closed, then $R^3_{1,m}$ is surjective.
\end{cor}

\begin{proof}
By Park \cite{Park1}, there exists $\alpha \in \TH^2(k,3;m)$ such that
\[
R^3_{1,m} (\alpha) = \beta \neq 0
\]
 in $k$. Fix any $0 \neq y \in k$.
$k= {\bar k}$ implies there are $z, {\beta}' \in k$ such that $z^{m+1} = y$
and ${({\beta}')}^{m+1} = {\beta}^{-1}$. Put $a = z{\beta}'$.
Then $a \in k^*$ and by corollary~\ref{cor:reg1}, 
$$R^3_{1,m} \left ( a^{-1} * {\alpha} \right ) = a^{m+1}{\beta}  
= (z{\beta}')^{m+1} {\beta} = y {\beta}^{-1} {\beta} = y.$$
Since $0$ is clearly in the image of $R^3_{1,m}$, we see that $R^3_{1,m}$
is surjective.
\end{proof}

Let $k$ be as before an arbitrary field of characteristic zero. Recall that by  theorem~\ref{thm:nolabel} we have the product  
\[
\boxtimes:z^r(k,*)\otimes \TZ^s(k,*;m) \to \TZ^{r+s}(k,*;m)
\]
induced by the external product. By theorem~\ref{thm:product}, the map $\boxtimes$ induces the natural
map
\begin{equation}\label{eqn:expro1} 
\mu: \CH^n(k,n) \otimes \TH^2(k,n;m) \to \TH^{n+2}(k,n+3;m)
\end{equation}

For any $n \ge 1$, a cycle in
$z^n(k,n)$ is a $\Z$-linear combination of closed points 
$b = (b_1, \cdots , b_n)\in\Box^n$, with
each $b_i $ in $\Box - \{0,\infty \}$.
For a curve $C \in B_3$, we denote by $C^b$ the curve in $B_{n+3} \cong
B_3 \times \Box^n$ given by the  product $C \times b$. We thus have
\begin{equation}\label{eqn:expro}
\mu \left ( [C] \otimes [b] \right ) = [C^b]
\end{equation}

\begin{lem}\label{lem:regprod}
Let $[C] \in \TZ^2(k,3;m)$ be the cycle defined by an irreducible curve
$C \subset B_3$ and let $b = (b_1, \cdots , b_n)\in\Box^n(k)$ be a $k$-point
defining a cycle in $z^n(k,n)$. Then
\begin{align*}
&R^{n+3}_{1,m} \left ([C^b] \right ) = \\
&\left(
\sum_{(C,P) \in {\mathcal M}^m (t_1)} R^3_{1,m} \left (C, P \right )
+ (-1)^n \sum_{(C,P) \in {\mathcal M}^m (t_2)} R^3_{1,m} 
\left (C, P \right ) \right ) \cdot \bigwedge_{i=1}^n \frac{db_i}{b_i}.
\end{align*}
\end{lem}

\begin{proof}
It is immediate that $\widehat{C^b} = \widehat{C} \times b$ and
$(\widehat{C^b})^N = \widehat{C}^N \times b$. Let
\[
p :\widehat{C}^N \to \widehat {B}_3;\ p^b : (\widehat{C^b})^N \to
\widehat{B}_{n+3}
\]
denote the maps on the normalizations, and for $P \in p^{-1} (\{x=0\})$, set
$P^b:=P\times b$. Then \\
\\
$(i) \ {(p^b)}^{-1} (\{x=0\}) = p^{-1} (\{x=0\}) \times b$. \\
$(ii)\ \text{for } l = 1, 2, \forall P \in p^{-1} (\{x=0\})$:
\[
 (C, P) \in {\mathcal M}^m(t_l)
\Leftrightarrow (C^b, P^b) \in {\mathcal M}^m(t_l). 
\]
$(iii) \  \text{for } l \ge 3, \forall P \in p^{-1}(\{x=0\})\!\!: p(P^b) \notin  F^1_{n+3,l}$.
\\
 \\
Thus
\begin{align*}
R^{n+3}_{1,m} \left ([C^b] \right ) & =  \sum_{Q \in 
(p^b)^{-1}(\{x=0\})}R^{n+3}_{1,m} \left ( C^b,Q \right ) \\
& =   \sum_{(C,P) \in {\mathcal M}^m (t_1)} R^{n+3}_{1,m} 
\left (C^b, Q \right ) + 
 \sum_{(C^b,Q) \in {\mathcal M}^m (t_2)} R^{n+3}_{1,m} 
\left (C^b, Q \right )  \\
& \hskip50pt + \sum_{l \ge 3} \sum_{(C^b,Q) \in {\mathcal M}^m (t_l)} R^{n+3}_{1,m} 
\left (C^b, Q \right )
\end{align*}
Note that the second equality holds because if  
\[
(C^b,Q) \in {\mathcal M}^m (t_l)
\cap {\mathcal M}^m (t_{l'})
\]
 for some $Q 
\in (p^b)^{-1}(\{x=0\})$, and for some $l \neq l'$, 
then $R^{n+3}_{1,m} \left ( C^b,Q \right ) = 0$ (see \cite{Park}, lemma~3.2).

Using $(i) , (ii)$ and $(iii)$ above, this gives
\begin{align*}
R^{n+3}_{1,m} \left ( [C^b] \right )
& = 
\sum_{(C^b,Q) \in {\mathcal M}^m (t_1)}
R^{n+3}_{1,m} \left (C^b, Q \right ) + 
 \sum_{(C^b,Q) \in {\mathcal M}^m (t_2)} 
R^{n+3}_{1,m} \left (C^b, Q \right ) \\
& =  \sum_{(C,P) \in {\mathcal M}^m (t_1)}
R^{n+3}_{1,m} \left (C^b, P^b \right ) + 
 \sum_{(C,P) \in {\mathcal M}^m (t_2)} 
R^{n+3}_{1,m} \left (C^b, P^b \right ) \\
& =  \sum_{(C,P) \in {\mathcal M}^m (t_1)}
{res}_{P^b} \left (  {\frac {1-t_1}{x^{m+1}}}  
\cdot \bigwedge_{j=2}^{n+2} {\frac {dt_j}{t_j}}  
\right )\\
&\hskip50pt -   \sum_{(C,P) \in {\mathcal M}^m (t_2)} {res}_{P^b} \left (
 {\frac {1-t_2}{x^{m+1}}}  \cdot \bigwedge_{j=3}^{n+2}
{\frac {dt_j}{t_j}}  \wedge {\frac {dt_1}{t_1}} \right ) \\ 
& =  \sum_{(C,P) \in {\mathcal M}^m (t_1)}
\left ( {res}_{P}  \frac {1-t_1}{x^{m+1}}\cdot
\frac {dt_2}{t_2} \right ) \cdot
\bigwedge_{j=1}^{n}
\frac {db_j}{b_j}\\
&\hskip50pt +
 (-1)^{n+1} \sum_{(C,P) \in {\mathcal M}^m (t_2)}
\left ( {res}_{P} \frac {1-t_2}{x^{m+1}}
\frac {dt_1}{t_1} \right )\cdot  \bigwedge_{j=1}^{n}
\frac {db_j}{b_j}  \\
& =  \sum_{(C,P) \in {\mathcal M}^m (t_1)} R^3_{1,m} \left (
C,P \right )\cdot\bigwedge_{j=1}^{n}\frac {db_j}{b_j}  \\
&\hskip50pt+  (-1)^n \sum_{(C,P) \in {\mathcal M}^m (t_2)}
R^3_{1,m} \left (C,P \right ) \cdot\bigwedge_{j=1}^{n}\frac {db_j}{b_j} 
\end{align*}
This completes the proof of the lemma.
\end{proof}

\begin{cor}\label{cor:even}
If $n$ is even, then 
\[
R^{n+3}_{1,m} \left ( {\mu}({\alpha} \otimes [b]) \right )
= R^3_{1,m} \left ( \alpha \right ) \cdot \bigwedge_{j=1}^{n}
{\frac {db_j}{b_j}} .
\]
\end{cor}

Next we recall some explicit $1$-cycles in $\TZ^2(k,3;m)$ constructed 
by Park \cite{Park1}. Let $a, a_1, a_2, \in k , b, b_1, b_2 \in k^*$.
Let ${\Gamma}_1, {\Gamma}_2, C^{(a_1,a_2),b}_1$, 
$C^{a,(b_1,b_2)}_2$ be parametrized $1$-cycles
in $\TZ^2(k,3;m)$ defined as follows: 
\begin{align*}
C^{(a_1,a_2),b}_1 &=  \begin{cases} 
\left \{ \left ( t, {\frac {(1-a_1 t)(1-a_2 t)}{1-(a_1 + a_2 )t}}, b 
\right ) | t \in k \right \} & {\rm if} \ a_1 a_2 (a_1 + a_2) \neq 0 \\
\left \{ \left ( t, 1- a^2 t^2 , b \right ) |\ t \in \A^1_k \right \} &
{\rm if} \ a:= a_1 = -a_2 \neq 0 \\
0 & {\rm if} \ a_1 a_2 = 0 
\end{cases}\\
\\
C^{a,(b_1,b_2)}_2 &=  \begin{cases}
\left \{ \left ( {\frac {1}{a}}, t, {\frac {b_1 t - b_1 b_2}{t - b_1 b_2}}
\right ) |\ t \in \A^1_k \right \} & {\rm if} \ a \neq 0 \\
0 & {\rm if} \ a = 0
\end{cases}\\
\\
{\Gamma}_1 &= \left \{ \left ( t, t, {\frac {(1-(1/2) t)^2}{1-t}} \right )
|\ t \in \A^1_k \right \}\\
\\
{\Gamma}_2 &= \left \{ \left ( t, 1+ {\frac {t} {6}}, 1- {\frac {t^2}{6}}
\right ) |\ t \in \A^1_k \right \}
\end{align*}

Let ${\ov {\Gamma}_1} = {\Gamma}_1 + C^{(1/2,1/2),2}_1$ and
\begin{align*}
{\ov {\Gamma}_2} & = & {\Gamma}_2 + C^{(-1/2,1/2),2/3}_1
+ 3  C^{(-1/3,-1/3),-1}_1 -  C^{(-1/6,-1/6),2}_1 \\
& & - C^{(-1/6,-1/3),2}_1 - C^{(-1/2,1/2),2}_1 + C^{(1/2,1/2),2}_1 
 \\
& & + C^{-1/6,(4,-2)}_2 -  C^{-1/6,(-2,-2)}_2 +
C^{1/2,(2/3,3/2)}_2 \\
& & - 3 C^{-1/6,(2,-1)}_2 - 3 C^{-1/3,(-1,-1)}_2 + 
C^{-1/2,(4/3,3/2)}_2 
\end{align*}
Define the cycle ${\Gamma}$ by
\begin{equation}\label{eqn:gamma}
{\Gamma} : = {\ov {\Gamma}_1} - {\ov {\Gamma}_2}
\end{equation}

\begin{lem}\label{lem:gamma1}
For each $b = (b_1, \cdots , b_n) \in {({\Box} \setminus \{0, \infty \})}^n(k)$ 
and
$a \in k^*$, we have 
\[
R^{n+3}_{1,m}\left ( \mu \left ( a* {\Gamma} \otimes [b] \right )
\right ) = (-1)^n  R^3_{1,m} \left ( a* {\Gamma} \right )\cdot
\bigwedge_{j=1}^{n}\frac {db_j}{b_j} \in \Omega^n_k.
\]
\end{lem}

\begin{proof} It is enough to show that the above formula holds for
a cycle of the form $[C] \in \left \{ \Gamma_1, \Gamma_2, 
C^{(a_1,a_2),b}_1, C^{a,(b_1,b_2)}_2 \right \}$. We do it separately
for each cycle in this set. \\
\\
$(i) \ \ [C] = {\Gamma}_1:$ By Park \cite{Park1}, lemma~1.9, for
each $p \in p^{-1}(\{x=0\})$, $(C, P)$ is in ${\mathcal M}^m (t_2)$ and hence by
lemma~\ref{lem:mod1} and lemma~\ref{lem:regprod}, we see immediately that the
formula holds for ${\Gamma}_1$. \\
\\
$(ii) \ \ [C] = {\Gamma}_2:$ In this case again by Park 
\cite{Park1}, lemma~1.10, for
each $p \in p^{-1}(\{x=0\})$, $(C, P)$ is in ${\mathcal M}^m (t_2)$ and hence
the formula follows. \\
\\
$(iii) \ \ [C] = C^{(a_1,a_2),b}_1 :$ In this case, we apply lemma~1.7 of Park
to see that for
each $p \in p^{-1}(\{x=0\})$, $(C, P)$ is in  ${\mathcal M}^m (t_1)$. 
Hence the second sum in lemma~\ref{lem:regprod} is zero which in
particular means that the first sum is same as 
$R^3_{1,m}(a*[C])$. Thus applying lemma~\ref{lem:regprod} again, we
get 
\[
R^{n+3}_{1,m}\left ( {\mu} \left ( a* [C] \otimes [b] \right ) \right )
= R^3_{1,m}(a*[C]) \cdot {\bigwedge_{j=1}^{n}} {\frac {db_j}{b_j}}. 
\]
So enough to show that the right hand term is zero.
However, in this case we have
$R^3_{1,m}(a*[C]) = {\frac {1}{a^{m+1}}} R^3_{1,m}([C])$ by 
corollary~\ref{cor:reg1}. But it is easy to see that 
$R^3_{1,m}([C]) = 0$ since $p^*(dt_2) = db = 0$. \\
\\
$(iv) \ \ [C] = C^{a,(b_1,b_2)}_2 : $
In this case again we have by lemma~1.7 of Park, for
each $p \in p^{-1}(\{x=0\})$, $(C, P)$ is in  ${\mathcal M}^m (t_2)$. We apply
lemma~\ref{lem:mod1} and lemma~\ref{lem:regprod} again
to conclude that the formula holds for in this case. 
\end{proof} 

\begin{lem}\label{lem:milnor}
Let $k$ be a field of characteristic different from $2$. Then there are 
surjective $k$-linear maps 
$$k {\otimes}_{\Z} T(k^*) \surj k {\otimes}_{\Z} {\bigwedge}(k^*)
\surj k {\otimes}_{\Z} K^M_*(k) \surj {\Omega}^*_k$$
where the $T(k^*)$ and ${\bigwedge}(k^*)$ denote the tensor and 
exterior algebras of $k^*$ and $K^M_*(k)$ is the Milnor $k$-theory
ring of $k$. The composite map is given in degree $n$ by
\[
a \otimes (b_1 \otimes \cdots \otimes b_n) \mapsto
a{\bigwedge_{j=1}^{n}} {\frac {db_j}{b_j}} \in {\Omega}^n_k.
\]
\end{lem}

\begin{proof}
The exterior algebra
of $k^*$ is the quotient of the tensor algebra by the homogeneous 
two-sided ideal $I$ generated by degree two elements 
$\left \{ b \otimes b | b \in k^* \right \}$. Let $I'$ be
the homogeneous 
two-sided ideal generated by degree two elements 
$\left \{ b_1 \otimes b_2 + b_2 \otimes b_1 | b_1, b_2 \in k^* \right \}$.
It is clear that $I' \subset I$.
However, since the characteristic of $k$ is different from $2$,
it is easy to check that $k {\otimes}_{\Z} (I/{I'}) = 0$. This shows
that $k {\otimes}_{\Z} {\bigwedge}(k^*)$ is the quotient of the 
$k$-algebra $k {\otimes}_{\Z} T(k^*)$ by the homogeneous ideal $I$.

It is known \cite{Nessuslin} that the Milnor $K$-theory ring of $k$ 
is the quotient of the tensor algebra of $k^*$ by the homogeneous 
two-sided ideal $J$ generated by degree two elements 
$\left \{ b \otimes (1-b)  | b, 1- b \in k^* \right \}$.
Thus to get the middle map of the lemma, it suffices to show that 
$I \subset J$. For this, it is enough to show that the images of the degree two 
generators of $I$ under the quotient map ${\tau} : T(k^*) \to K^M_*(k)$
are all zero. But we have 
\[
{\tau}(b_1 \otimes  b_2 + b_2 \otimes b_1 )
= \{ b_1, b_2 \} + \{b_2, b_1 \} = 0
\]
 in $K^M_2(k)$ by \cite{Nessuslin}, lemma~3.2(b).

By \cite{BlochEsnault}, lemma~4.1, there is a surjective map
of $k$-vector spaces 
\[
{\eta} : k {\otimes}_{\Z} T(k^*) \to
{\Omega}^*_k
\]
 which is given in degree $n$ by
\[
{\eta} (a \otimes (b_1 \otimes \cdots \otimes b_n)) = 
a {\bigwedge_{j=i}^{n}} {\frac {db_j}{b_j}}.
\]
To show that this factors through $k \otimes K^M_*(k)$, it suffices to
show that ${\eta} \left ( k \otimes J \right ) = 0$. A generator of 
$k \otimes J$ is of the form $a \otimes (b \otimes (1-b))$ and one has
\begin{align*}
{\eta} \left ( a \otimes (b \otimes (1-b)) \right ) &= 
{\frac {a} {b(1-b)}} db \wedge d(1-b) \\
&= - {\frac {a} {b(1-b)}}(db \wedge 
db) = 0.
\end{align*}
\end{proof}

We shall also denote the induced map $k \otimes K^M_*(k) \to  {\Omega}^*_k$
by $\eta$. Now we go back to the situation of $k$ being of characteristic
zero. Let $\Gamma \in \TH^2(k,3;m)$ be the class of the cycle defined in 
\eqref{eqn:gamma}.

\begin{prop}\label{prop:commutative}
For $n \ge 0$, there are maps 
\[
\xymatrix{
CH^n(k,n) {\otimes}_{\Z} \TH^2(k,3;m) \ar[r]^{\mu} \ar[d]^{S^n_{1,m}} &
\TH^{n+2}(k,n+3;m) \ar[d]^{R^{n+3}_{1,m}} \\
K^M_n(k) {\otimes}_{\Z} k \ar[r]^{\eta} & {\Omega}^n_k 
}
\]   
such that for $a \in k^*$ and ${\delta} \in CH^n(k,n)$, one has
$${R^{n+3}_{1,m}} \circ {\mu} \left ( a * {\Gamma} \otimes {\delta}
\right ) = (-1)^n {\eta} \circ {S^n_{1,m}} \left ( a*{\Gamma} \otimes
{\delta} \right ).$$
\end{prop}  

\begin{proof}
The map $\mu$ comes from theorem~\ref{thm:product} (part 1). ${S^n_{1,m}}$
is defined by ${S^n_{1,m}} = {\phi}^n \otimes R^3_{1,m}$ where ${\phi}^n 
: CH^n(k,n) \to K^M_n(k)$ is the Nesterenko-Suslin-Totaro isomorphism
\cite{Nessuslin}. We can write $\delta = \sum_{i=1}^{r}b^i -
\sum_{j=1}^{s}c^j$ where $b^i, c^j \in ({\Box} \setminus \{0,\infty\})^n$. 
Then
\begin{align*}
{R^{n+3}_{1,m}} \circ {\mu} \left ( a * {\Gamma} \otimes {\delta}\right ) & = 
{\sum} {R^{n+3}_{1,m}} \circ {\mu} \left ( a * {\Gamma} \otimes b^i 
\right ) \\
&\hskip80pt - {\sum} {R^{n+3}_{1,m}} \circ {\mu} \left (a * {\Gamma} \otimes c^j \right ) \\
& =  (-1)^n {\sum} {R^{3}_{1,m}} \left (   a * {\Gamma} \right ) {\bigwedge_{l=1}^{n}} {\frac {db^i_l}{b^i_l}} \\
&\hskip80pt - (-1)^n {\sum} {R^{3}_{1,m}} \left (   a * {\Gamma} \right )
{\bigwedge_{l=1}^{n}} {\frac {dc^j_l}{c^j_l}} \\
& \hskip150pt ({\rm by \ \ lemma~\ref{lem:gamma1}})  \\
& =  (-1)^n {\sum} {R^{3}_{1,m}} \left (   a * {\Gamma} \right ) 
{\phi}^n ([b^i]) \\
& \hskip80pt - (-1)^n {\sum} {R^{3}_{1,m}} \left (   a * {\Gamma} \right )
{\phi}^n ([c^j]) \\
& =  (-1)^n {R^{3}_{1,m}} \left (   a * {\Gamma} \right ) \left \{
{\sum}{\phi}^n ([b^i]) - {\sum}{\phi}^n ([d^j]) \right \}  \\
& =  (-1)^n {R^{3}_{1,m}} \left (   a * {\Gamma} \right ) \left (
{\eta} \left ( 1 \otimes {\phi}^n (\delta) \right ) \right ) \\
& =  (-1)^n {\eta} \left \{ {R^{3}_{1,m}} \left (   a * {\Gamma} \right )
\otimes {\phi}^n (\delta) \right \} \\
& =  (-1)^n {\eta} \circ S^n_{1,m} \left (  a * {\Gamma} \otimes 
\delta \right )
\end{align*}
\end{proof}

\begin{cor}\label{cor:gamma2}
Let ${\widetilde {\TH}}^2(k,3;m) \hookrightarrow \TH^2(k,3,m)$ 
be the ${\Z}[k^*]$-submodule
generated by $\Gamma$. Then for each element $\alpha$ in 
$CH^n(k,n) \otimes {\widetilde {\TH}}^2(k,3;m)$, one has 
$${R^{n+3}_{1,m}} \circ {\mu} (\alpha) = 
(-1)^n {\eta} \circ {S^n_{1,m}} \left ( \alpha \right ).$$
\end{cor}

\begin{proof}
This follows from proposition~\ref{prop:commutative} since elements of 
$CH^n(k,n) \otimes {\widetilde {\TH}}^2(k,3;m)$ are $\Z$-linear combinations
of elements of the type $a * {\Gamma} \otimes [b]$ with $a \in k^*$.
\end{proof} 

\begin{proof}[Proof of theorem~\ref{thm:reg}]
For $n \le 2$, the theorem is obvious as the term on the right is zero.
Thus it is enough to show that $R^{n+3}_{1,m}$ is surjective for $n \ge 0$.

Let  
\[
\iota : CH^n(k,n) \otimes {\widetilde {\TH}}^2(k,3;m) \to 
CH^n(k,n) \otimes {\TH}^2(k,3;m) 
\]
be the natural map.   It is enough to show that ${R^{n+3}_{1,m}} \circ {\mu} \circ \iota$
is surjective. $\eta$ is surjective by lemma~\ref{lem:milnor}. ${\phi}^n$
is isomorphism for all $n \ge 0$ by \cite{Nessuslin}.
By \cite{Park1}, proposition~1.11, $R^3_{1,m} (\Gamma) \neq 0$. Now since
$k = {\ov k}$, we repeat the argument of the proof of 
corollary~\ref{reg2} to see that the map 
\[
R^3_{1,m}:{\widetilde {\TH}}^2(k,3;m) \to k 
\]
is surjective. Combining these
implies that $\eta \circ S^n_{1,m} \circ \iota$ is surjective.
We now apply corollary~\ref{cor:gamma2} to show
$(-1)^n {R^{n+3}_{1,m}} \circ {\mu} \circ \iota$ is surjective. 
In particular, 
\[
{\Omega}^n_k \subset {\rm Image}
({R^{n+3}_{1,m}} \circ {\mu} \circ \iota).
\]
\end{proof}

School of Mathematics \\
Tata Institute of Fundamental Research \\
Homi Bhabha Road \\
Mumbai, 400005, India. \\
e-mail : amal@math.tifr.res.in \\
\\
Northeastern University \\
Department of Mathematics \\
360 Huntington Avenue \\
Boston, MA 02115, USA \\
e-mail :  marc@neu.edu \\
\end{document}